\documentclass[11pt,a4paper]{article}
\usepackage[T1]{fontenc}
\usepackage[utf8]{inputenc}
\usepackage[english]{babel}
\usepackage{lmodern,microtype}
\usepackage[margin=2.5cm]{geometry}
\usepackage{amsmath,amssymb,amsthm,mathtools}
\usepackage{enumitem,needspace}
\usepackage{graphicx,booktabs,array,multirow}
\usepackage[authoryear,round]{natbib}
\usepackage[hidelinks]{hyperref}
\hypersetup{pdftitle={Recursive estimation of the Oja median and related simplex medians},pdfauthor={Antoine Godichon-Baggioni and Guillaume Sallé}}
\newtheorem{theorem}{Theorem}[section]
\newtheorem{proposition}[theorem]{Proposition}
\newtheorem{lemma}[theorem]{Lemma}

\theoremstyle{definition}
\newtheorem{assumption}{Assumption}

\theoremstyle{remark}
\newtheorem*{remark}{Remark}

\newcommand{\R}{\mathbb R}
\newcommand{\E}{\mathbb E}
\newcommand{\Prob}{\mathbb P}
\newcommand{\F}{\mathcal F}
\newcommand{\Id}{I_d}
\newcommand{\ind}{\mathbf 1}
\newcommand{\norm}[1]{\left\lVert#1\right\rVert}
\newcommand{\abs}[1]{\left|#1\right|}
\newcommand{\scal}[2]{\left\langle #1,#2\right\rangle}
\newcommand{\as}{\text{a.s.}}
\newcommand{\op}{\mathrm{op}}
\DeclareMathOperator{\Simpl}{Simpl}
\DeclareMathOperator{\aff}{aff}
\DeclareMathOperator{\spanop}{span}
\DeclareMathOperator{\dist}{dist}
\DeclareMathOperator{\Gr}{Gr}
\DeclareMathOperator{\Var}{Var}
\DeclareMathOperator{\Cov}{Cov}
\DeclareMathOperator{\sgn}{sgn}
\setlist[enumerate]{itemsep=0.3em,topsep=0.4em}
\allowdisplaybreaks[2]
\newcommand{\doi}[1]{\href{https://doi.org/#1}{doi:\nolinkurl{#1}}}
\numberwithin{equation}{section}
\begin{document}
\title{Recursive estimation of the Oja median\protect\\ and related simplex medians}
\author{Antoine Godichon-Baggioni$^{1}$ \and Guillaume Sall\'e$^{2}$}
\date{}
\maketitle
\begin{center}
\small
$^{1}$Sorbonne Universit\'e, Universit\'e Paris Cit\'e, CNRS,\\
Laboratoire de Probabilit\'es, Statistique et Mod\'elisation, LPSM,\\
F-75005 Paris, France\par\smallskip
$^{2}$INSA Rouen Normandie, Laboratoire de math\'ematiques de l'INSA,\\
76800 Saint-\'Etienne-du-Rouvray, France
\end{center}

\begin{abstract}
We study the recursive estimation of multivariate location functionals defined by minimizing the expected volume of a random simplex. This framework encompasses both the spatial median and the affine-equivariant Oja median. We introduce a sliding-window stochastic gradient algorithm and its averaged version, and establish their strong consistency and asymptotic normality. We further show that the sliding-window approach offers a favorable trade-off between statistical accuracy and computational cost compared with a naive disjoint-block stochastic gradient algorithm and offline methods, with the statistical comparison based on asymptotic covariances. Numerical experiments illustrate the computational advantages of the proposed procedures and their finite-sample accuracy.
\end{abstract}


\section{Introduction}\label{sec:introduction}

The estimation of multivariate location is a classical problem in nonparametric statistics. The mean vector is affine-equivariant: applying a nonsingular affine transformation to a distribution transforms its mean in the same way. Nevertheless, its sensitivity to extreme observations motivates the use of multivariate medians. The spatial median, defined by minimizing expected Euclidean distance, is a natural extension of the univariate median; see \citet{Hal48,brown1983statistical,Kem87}. It is equivariant under translations, orthogonal transformations and uniform changes of scale, but not under general nonsingular affine transformations. The Oja median, introduced by \citet{oja1983descriptive}, is affine-equivariant. It is defined by minimizing the expected volume of a full-dimensional simplex whose vertices are the candidate location and independent observations from the distribution.

These two medians belong to the family of location functionals studied by \citet{durre_affine-equivariant_2022}. For $p\geq1$ and $1\leq\ell\leq d$, their population objective function is
\begin{equation*}
 \mu\longmapsto
 \E\left[m_\ell\left(\Simpl\left(X_1,\ldots,X_\ell,\mu\right)\right)^p\right],
\end{equation*}
where $X_1,\ldots,X_\ell$ are independent random vectors with a common distribution on $\R^d$, $m_\ell$ denotes $\ell$-dimensional volume, and $\Simpl$ denotes the convex hull of its arguments. For $p=1$, the choices $\ell=1$ and $\ell=d$ give the spatial and Oja medians, respectively. Intermediate values of $\ell$ give related location functionals based on lower-dimensional simplices. We call these functionals \emph{simplex medians}. We study recursive estimation in the median case $p=1$, for the full range $1\leq\ell\leq d$. In this case, the norm of the subgradient is bounded uniformly in $\mu$ by the volume of the simplex formed by the $\ell$ observations, divided by $\ell$. This property allows us to control the successive updates through moments of this simplex volume.

Given $n$ observations, the empirical objective function of the complete offline estimator is a $U$-statistic of order $\ell$, obtained by averaging over all $\binom n\ell$ subsets of $\ell$ observations \citep[Section~2.1]{durre_affine-equivariant_2022}. Evaluating all terms of this sum has a computational cost of order $n^\ell$ for fixed $d$ and $\ell$, and becomes prohibitive when $n$ is large and $\ell$ is not small, particularly for the Oja median, for which $\ell=d$ \citep{Fischer2020}. To reduce this cost, \citet[Section~2.4]{durre_affine-equivariant_2022} consider an incomplete $U$-statistic based on a random selection of tuples from the observed sample, greatly reducing the number of terms in the empirical objective function. However, these offline methods still involve repeated evaluations of an empirical objective function based on a fixed dataset. In this article, we propose an online stochastic gradient algorithm that uses a sliding window to update the estimate whenever a new observation arrives and requires a fixed amount of storage.

Stochastic approximation provides a natural method for this purpose \citep{robbins1951}. For the spatial median, \citet{CCZ2013} study the almost sure convergence and mean-square convergence rates of a stochastic gradient estimator in a Hilbert space, and establish asymptotic normality of its averaged version. Nonasymptotic confidence balls are studied by \citet{CCG2017}. Averaging the iterates follows the method of \citet{ruppert1988efficient} and \citet{PolyakJud92}; see also \citet{pelletier2000asymptotic}. Since each stochastic gradient depends on $\ell$ observations, applying stochastic gradient descent with independent data blocks requires waiting for $\ell$ new observations before each update. However, this gives only $\lfloor n/\ell\rfloor$ updates after $n$ observations. We instead use a sliding-window stochastic gradient estimator, where an update is performed with each new observation, using this observation and the $\ell-1$ preceding ones. Therefore, over a single pass through $n$ observations, the algorithm performs $n-\ell+1$ updates. As discussed below, reusing observations in overlapping windows also reduces the asymptotic variance of the averaged estimator.

The use of overlapping windows breaks the independence between the current iterate and the observations used to compute the next stochastic gradient, so that standard results based on conditionally unbiased gradients cannot be applied directly. Stochastic approximation with dependent observations has a substantial literature; see \citet{benveniste-book90,Kushner2003stochastic}, and \citet{KushnerYang1993} for averaging. Here, the difficulty is to control the variation of the stochastic gradient under assumptions on the common distribution of the observations $X_i$, even though the simplex volume is not differentiable everywhere.

In the particular case of the Oja median, the gradient of the simplex volume is constant on each side of the hyperplane determined by the observations, with opposite values on the two sides. Under suitable assumptions on the density of the observations $X_i$, we control the variation of this gradient in expectation and use the resulting bounds in the proofs of convergence rates and asymptotic normality for the recursive estimator and its averaged version. This case, corresponding to $\left(p,\ell\right)=\left(1,d\right)$, is excluded from the limit theorem for the offline estimator in \citet[Theorem~2.5]{durre_affine-equivariant_2022}. Earlier asymptotic results for the sample Oja median under multivariate normality were obtained by \citet{oja1985asymptotic}.

We first prove almost sure convergence without any density assumption, and then establish almost sure convergence rates and asymptotic normality for both the recursive estimator and its averaged version under a polynomial upper bound on the density and suitable regularity and step-size conditions. The averaged estimator has squared error $O\left(n^{-1}\log\log n\right)$ almost surely and is asymptotically normal at the $\sqrt n$ scale. For $\ell<d$, the rates and limit theorems for both estimators also hold under alternative assumptions that do not require a density. The asymptotic covariance matrix of the averaged sliding-window estimator lies between those of the complete offline estimator and the averaged disjoint-block estimator in the Loewner order, with strict inequalities when $\ell\geq2$ and the distribution has a density. The two recursions have the same order of computational cost.

The paper is organized as follows.  Section~\ref{sec::framework} gives the statistical setting and assumptions.  Section~\ref{sec::algo} defines the recursive estimators, states the convergence results and compares their limiting covariance matrices and computational costs.  Section~\ref{sec::simus} contains the simulation study.  The proofs and auxiliary results are postponed to the appendices.

\section{Framework}\label{sec::framework}\label{sec:framework}

\noindent\emph{Notation.}\enspace
We write $\norm{\cdot}$ for the Euclidean norm and
$\norm{\cdot}_{\op}$ for the induced operator norm.
For a finite tuple $\mathbf x$ of points in $\R^d$, we denote by
$\Simpl(\mathbf x)$ and $\aff(\mathbf x)$ its convex hull and
affine span, respectively.
For $j\in\{0,\ldots,d\}$, we write $m_j$ for
$j$-dimensional volume.
We write $\dist(u,E)$ for the Euclidean distance from $u$
to a set $E$.

\subsection{Simplex medians and their empirical counterparts}

Let $P$ be a probability measure on $\R^d$, with $d\geq2$,
and fix $1\leq\ell\leq d$.
Let $(X_i)_{i\geq1}$ be a sequence of independent random vectors
with common distribution $P$, and write
$\mathbf X\coloneqq\left(X_1,\ldots,X_\ell\right)$.
For $\mathbf x\coloneqq\left(x_1,\ldots,x_\ell\right)\in(\R^d)^\ell$
and $\mu\in\R^d$, we define $g(\mathbf x;\mu)$ as the
$\ell$-dimensional volume of the simplex with vertices
$x_1,\ldots,x_\ell,\mu$:
\begin{equation*}
 g\left(\mathbf x;\mu\right)\coloneqq
 m_\ell\left(\Simpl\left(x_1,\ldots,x_\ell,\mu\right)\right).
\end{equation*}
We study the minimizer $\mu^*$ of
\begin{equation*}
 G\left(\mu\right)\coloneqq \E\left[g\left(\mathbf X;\mu\right)\right].
\end{equation*}
For $\ell=1$, $G\left(\mu\right)=\E\left[\norm{X_1-\mu}\right]$ and $\mu^*$ is the spatial median. For $\ell=d$, it is the Oja median.

If the points $x_1,\ldots,x_\ell$ are affinely independent, let $\Gamma\left(\mathbf x\right)$ be the orthogonal projection onto the orthogonal complement of
$\spanop\left(x_1-x_\ell,\ldots,x_{\ell-1}-x_\ell\right)$.
For such a tuple and every $\mu\in\R^d$, we have
$\norm{\Gamma(\mathbf x)(\mu-x_\ell)}
=\dist(\mu,\aff(\mathbf x))$.
For $\ell=1$, we set $\Gamma\left(x\right)\coloneqq I_d$ and $m_0\left(\Simpl\left(x\right)\right)\coloneqq 1$. For $\ell>1$, we set $\Gamma\left(\mathbf x\right)\coloneqq 0$ on affinely dependent tuples; both $g\left(\mathbf x;\mu\right)$ and $m_{\ell-1}\left(\Simpl\left(\mathbf x\right)\right)$ are then zero. The base--height formula \citep[Lemma~S.1.1]{DPsupp} gives
\begin{equation}\label{eq:base-height}
 g\left(\mathbf x;\mu\right)=\frac{m_{\ell-1}\left(\Simpl\left(\mathbf x\right)\right)}{\ell}
                   \norm{\Gamma\left(\mathbf x\right)\left(\mu-x_\ell\right)}.
\end{equation}
This representation shows that $g\left(\mathbf x;\cdot\right)$ is convex. We use the following selection of its subgradient:
\begin{equation}\label{eq:little_g_grad}
 \nabla_\mu g\left(\mathbf x;\mu\right)\coloneqq
 \begin{cases}
 \displaystyle\frac{m_{\ell-1}\left(\Simpl\left(\mathbf x\right)\right)}{\ell}
 \frac{\Gamma\left(\mathbf x\right)\left(\mu-x_\ell\right)}
 {\norm{\Gamma\left(\mathbf x\right)\left(\mu-x_\ell\right)}},&\text{if }\Gamma\left(\mathbf x\right)\left(\mu-x_\ell\right)\ne0,\\[0.8em]
 0,&\text{if }\Gamma\left(\mathbf x\right)\left(\mu-x_\ell\right)=0.
 \end{cases}
\end{equation}
This notation denotes the usual gradient when it exists. At a point of nondifferentiability, the chosen value is a subgradient.

For $n\geq\ell$, the empirical objective function of the offline
estimator is
\begin{equation}\label{eq:empirical-criterion}
 \widehat G_n\left(u\right)\coloneqq \binom n\ell^{-1}
 \sum_{1\leq i_1<\cdots<i_\ell\leq n}g\left(X_{i_1},\ldots,X_{i_\ell};u\right).
\end{equation}
When its set of minimizers is nonempty, let $\widehat\mu_n$ denote any measurable choice of minimizer. Under the density assumption below, this set is nonempty, compact and convex almost surely for $n\geq\ell+1$. Existence and measurability are established in the proof of
Theorem~\ref{th:batch} in Appendix~\ref{sec:batch-proof}.

If $P$ is centrally symmetric about $a$ and the minimizer is unique, then $\mu^*=a$. Indeed, reflection of all the vertices about $a$ preserves the simplex volume, so $G\left(a+u\right)=G\left(a-u\right)$. For $\ell=d$, an invertible affine transformation $x\mapsto Ax+b$ multiplies every simplex volume by $|\det A|$. Thus the unique population minimizer transforms into $A\mu^*+b$; see \citet[Section~2.2]{durre_affine-equivariant_2022}. The set of empirical minimizers has the same transformation property.

\subsection{Assumptions on the distribution}

For almost sure convergence, we only require the following
condition on $P$.

\begin{assumption}
\label{ass:consistency}
The function $G$ has a unique minimizer $\mu^*\in\R^d$, and
$\E\left[m_{\ell-1}\left(\Simpl\left(\mathbf X\right)\right)^2\right]
<\infty$.
\end{assumption}

Since $G$ has a unique minimizer, it is not identically infinite,
and hence $G\left(\mu^*\right)<\infty$.
By \eqref{eq:base-height} and the moment condition in Assumption~\ref{ass:consistency}, $G$ is finite and continuous on $\R^d$.  By \eqref{eq:product-bound}, the moment condition in Assumption~\ref{ass:consistency} is satisfied whenever $\E\left[\norm{X_1}^2\right]<\infty$.

For the convergence rates and asymptotic distributions, we use the following assumptions for both $\ell<d$ and $\ell=d$.

\begin{assumption}
\label{ass:density}
The law $P$ has a Borel measurable density $f$ with respect to
Lebesgue measure, and there exist $\nu>2$ and a finite constant
$C_f>0$ such that, for every $x\in\R^d$,
\begin{equation*}
 0\leq f\left(x\right)
 \leq \frac{C_f}{\left(1+\norm{x}\right)^{d+\nu}}.
\end{equation*}
\end{assumption}

Assumption~\ref{ass:density} implies
Assumption~\ref{ass:consistency}; see
Proposition~\ref{prop:population}.
It gives $\E\left[\norm{X_1}^q\right]<\infty$
for every $0<q<\nu$.
The same proposition shows that $G$ is finite, strictly convex,
coercive and differentiable, and that its gradient is
Lipschitz.

\begin{assumption}
\label{ass:taylor}
At the minimizer $\mu^*$ defined in Assumption~\ref{ass:consistency}, the function $G$ is twice differentiable, with positive definite Hessian $H\coloneqq \nabla^2G\left(\mu^*\right)$. Moreover, there exist an exponent $\alpha$ with $0<\alpha\leq1$
and finite constants $r_{\alpha},C_{\alpha}>0$ such that,
for every $u\in\R^d$ with $\norm{u}\leq r_{\alpha}$,
the first-order Taylor remainder of $\nabla G$ at $\mu^*$ satisfies
\begin{equation*}
 \norm{\nabla G\left(\mu^*+u\right)-Hu}
 \leq C_{\alpha}\norm{u}^{1+\alpha}.
\end{equation*}
\end{assumption}

\begin{remark}
The remainder bound in Assumption~\ref{ass:taylor} can be verified from H\"older continuity of the Hessian at $\mu^*$. More precisely, suppose that $G$ is twice differentiable on a
neighbourhood of $\mu^*$ and that there exist an exponent $\alpha$
with $0<\alpha\leq1$ and finite constants $L_{\alpha},r>0$ such that,
for every $u\in\R^d$ with $\norm{u}\leq r$,
\begin{equation*}
 \norm{\nabla^2G\left(\mu^*+u\right)-H}_{\op}
 \leq L_{\alpha}\norm{u}^{\alpha}.
\end{equation*}
Since $\nabla G\left(\mu^*\right)=0$, integration along the segment from $\mu^*$ to $\mu^*+u$ gives
\begin{equation*}
 \nabla G\left(\mu^*+u\right)-Hu
 =\int_0^1\left(\nabla^2G\left(\mu^*+tu\right)-H\right)u\,dt.
\end{equation*}
The remainder is therefore bounded by $L_{\alpha}\norm{u}^{1+\alpha}/\left(1+\alpha\right)$, so the bound in Assumption~\ref{ass:taylor} holds with $r_{\alpha}=r$ and $C_{\alpha}=L_{\alpha}/\left(1+\alpha\right)$.

Local linearization with a remainder of order $\norm{u}^{1+\alpha}$ is used by \citet[Assumption~3.2]{PolyakJud92}. For $\alpha=1$, the remainder bound in Assumption~\ref{ass:taylor} is the local quadratic remainder condition of \citet[Assumption~5]{GodichonWerge2025}.
\end{remark}

\subsection{A sufficient condition for the local expansion}

The following proposition gives a sufficient condition for
Assumption~\ref{ass:taylor} that can be checked directly from
the density and its first derivative.

\begin{proposition}\label{prop:smooth-density}
Suppose that Assumption~\ref{ass:density} holds with a density
$f\in C^1\left(\R^d\right)$.
For the exponent $\nu$ in Assumption~\ref{ass:density},
suppose that there exists a finite constant $C_{\nu}>0$
such that, for every $x\in\R^d$,
\begin{equation*}
 \norm{\nabla f\left(x\right)}
 \leq \frac{C_{\nu}}{\left(1+\norm{x}\right)^{d+\nu}}.
\end{equation*}
Then, for every $1\leq\ell\leq d$, the function $G$ is twice
continuously differentiable, and its Hessian is
Lipschitz and positive definite at every point of $\R^d$.
In particular, Assumption~\ref{ass:taylor} holds with $\alpha=1$.
\end{proposition}

The proof is given in Appendix~\ref{sec:hessian}.
The derivative bound is used to control the variation of the
population Hessian through translations of the density.

These conditions permit polynomial tails.
Nondegenerate Gaussian densities and finite mixtures of such
densities satisfy the hypotheses of
Proposition~\ref{prop:smooth-density} for every $\nu>2$.
The same is true of $C^1$ densities with compact support.
For an elliptical Student density with $\nu>2$ degrees of
freedom, the hypotheses hold with this value of $\nu$.
These assertions are verified at the end of
Appendix~\ref{sec:hessian}.

\subsection{Alternative assumptions when $\ell<d$}

For $p=1$, \citet[Theorem~2.5]{durre_affine-equivariant_2022} establish the asymptotic distribution of the offline estimator for $1\leq\ell<d$ under moment and geometric conditions, without requiring a density. For the recursive estimator with $1\leq\ell<d$,
Assumption~\ref{ass:density} can be replaced by the following
conditions, with Assumptions~\ref{ass:consistency}
and~\ref{ass:taylor} imposed separately.

\begin{assumption} \label{ass:moment-alternative}
The law $P$ is $\left(\ell-1\right)$-smooth, meaning that it assigns
probability zero to every affine subspace of dimension $\ell-1$,
and there exists $q>2$ such that
$\E\left[\norm{X_1}^q\right]<\infty$.
\end{assumption}

\begin{assumption} \label{ass:inverse-alternative}
For $2\leq\ell<d$,
\begin{equation*}
 \sup_{u\in\R^d}\E\left[
 \frac{
 m_{\ell-1}\left(\Simpl\left(\mathbf X\right)\right)
 \prod_{i=1}^{\ell}\left(1+\norm{X_i}\right)}
 {\ell\norm{\Gamma\left(\mathbf X\right)\left(u-X_\ell\right)}}
 \right]<\infty.
\end{equation*}
\end{assumption}

The $\left(\ell-1\right)$-smoothness condition in Assumption~\ref{ass:moment-alternative} is that of \citet[Section~2.3]{durre_affine-equivariant_2022}. For $p=1$, their Theorem~2.5 requires a finite second moment, whereas Assumption~\ref{ass:moment-alternative} requires a moment of order $q>2$. For $2\leq\ell<d$, Assumption~\ref{ass:inverse-alternative} strengthens their inverse-distance condition (2.6), which is imposed only at $u=\mu^*$ and does not contain the factor $\prod_{i=1}^{\ell}\left(1+\norm{X_i}\right)$. Here, the additional factor and the uniform bound over $u$ are used to control the difference between gradients evaluated at iterates that depend on overlapping observations; see \eqref{eq:overlap-first-moment}.

Under Assumption~\ref{ass:density}, the bound in Assumption~\ref{ass:inverse-alternative} follows from \eqref{eq:weighted-inverse} in Lemma~\ref{lem:oscillation}.  When $\ell=1$, successive windows are independent, and Assumption~\ref{ass:inverse-alternative} is not required in Proposition~\ref{cor:nonoja-alternative}.

\section{Recursive estimators and asymptotic results}\label{sec::algo}\label{sec:results}
\subsection{The sliding-window recursion}

Let $\mu_0$ be an $\R^d$-valued random vector and let $c_\gamma$ be a strictly positive random variable. Suppose that $\left(\mu_0,c_\gamma\right)$ is independent of $\left(X_i\right)_{i\geq1}$ and that $\E\left[\norm{\mu_0}^2+c_\gamma^2\right]<\infty$. For $0\leq j<\ell$, we set $\mu_j\coloneqq \mu_0$. For $n\geq\ell$, we define
\begin{equation}\label{eq:recursion}
 \mu_n\coloneqq \mu_{n-1}
 -\gamma_n\nabla_\mu g\left(X_{n-\ell+1},\ldots,X_n;\mu_{n-1}\right),
\end{equation}
where the step sizes are given by
$\gamma_n\coloneqq c_\gamma n^{-\gamma}$,
with a fixed exponent $1/2<\gamma<1$. At each step, one observation enters the window and the oldest observation is removed. For every $n\geq\ell$, the averaged estimator is defined by
\begin{equation}\label{eq:averaged-estimator}
 \bar\mu_n\coloneqq \frac1{n-\ell+1}\sum_{k=\ell}^n\mu_k,
\end{equation}
and can be computed recursively, for every $n\geq\ell$, as
\[
 \bar\mu_{n+1}=\bar\mu_n+\frac{\mu_{n+1}-\bar\mu_n}{n-\ell+2},
 \quad\text{and}\quad\bar\mu_\ell=\mu_\ell.
\]

For $\ell\geq2$ and $n>\ell$, the observations
$X_{n-\ell+1},\ldots,X_{n-1}$ have already been used to compute
$\mu_{n-1}$.
Even when $G$ is differentiable, the conditional expectation
\[
 \E\left[
 \nabla_\mu g\left(X_{n-\ell+1},\ldots,X_n;\mu_{n-1}\right)
 \mid\mu_0,c_\gamma,X_1,\ldots,X_{n-1}
 \right]
\]
need not equal $\nabla G\left(\mu_{n-1}\right)$.
In the proofs, we compare $\mu_{n-1}$ with $\mu_{n-\ell}$,
which is independent of $\left(X_{n-\ell+1},\ldots,X_n\right)$.
The resulting error is controlled in \eqref{eq:overlap-first-moment}
and \eqref{eq:overlap-variance}.

\subsection{Convergence, rates and limit distributions}

We first state the almost sure convergence result. It uses only Assumption~\ref{ass:consistency}; in particular, no density assumption is needed.

\begin{theorem}\label{th:consistency}\label{theo::consistency}
Under Assumption~\ref{ass:consistency}, the recursion \eqref{eq:recursion} satisfies
\begin{equation*}
 \sup_n\E\left[\norm{\mu_n-\mu^*}^2\right]<\infty
 \qquad\text{and}\qquad
 \mu_n\xrightarrow[n\to\infty]{\hspace{1.5em}}\mu^*\quad\as
\end{equation*}
\end{theorem}

The proof in Appendix~\ref{sec:consistency} uses the convexity of $g\left(\mathbf x;\cdot\right)$ to establish the uniform second-moment bound before proving almost sure convergence.

We next give the almost sure convergence rates and the asymptotic distributions of both estimators.

\begin{theorem}
\label{th:main}\label{theo::av}\label{theo::rate}
Suppose that Assumptions~\ref{ass:density} and~\ref{ass:taylor}
hold and that $\gamma\nu>2$. Let $H\coloneqq \nabla^2G\left(\mu^*\right)$.
For $j\in\{0,\ldots,\ell-1\}$, we set
\begin{equation*}
 \Gamma_j\coloneqq \E\left[\nabla_\mu g\left(X_1,\ldots,X_\ell;\mu^*\right)
      \nabla_\mu g\left(X_{1+j},\ldots,X_{\ell+j};\mu^*\right)^\top\right].
\end{equation*}
We also set
\begin{align}
 \Sigma\coloneqq \Gamma_0+\sum_{j=1}^{\ell-1}\left(\Gamma_j+\Gamma_j^\top\right)
 &\quad\text{and}\quad
 V_H\coloneqq \int_0^\infty e^{-tH}\Sigma e^{-tH}\,dt.\label{eq:Sigma}\\
\intertext{Then}
 \norm{\mu_n-\mu^*}^2=O\left(n^{-\gamma}\log n\right)\quad\as
 &\quad\text{and}\quad
 \gamma_n^{-1/2}\left(\mu_n-\mu^*\right)
 \xrightarrow[n\to\infty]{\hspace{.5em}\mathcal L\hspace{.5em}}\mathcal N_d\left(0,V_H\right).\notag\\
\intertext{If, in addition, $\gamma\left(1+\alpha\right)>1$, then}
 \norm{\bar\mu_n-\mu^*}^2
 =O\left(\frac{\log\log n}{n}\right)\quad\as
 &\quad\text{and}\quad
 \sqrt n\left(\bar\mu_n-\mu^*\right)
 \xrightarrow[n\to\infty]{\hspace{.5em}\mathcal L\hspace{.5em}}
 \mathcal N_d\left(0,H^{-1}\Sigma H^{-1}\right).\notag
\end{align}
\end{theorem}

The restrictions on $\gamma$ have separate uses. The bound $\gamma>1/2$ gives summability of the squared step sizes and makes the accumulated conditional bias negligible after division by $\sqrt n$. The condition $\gamma<1$ controls the boundary term from summation by parts. The moment restriction $\gamma\nu>2$ allows a choice $2/\gamma<q<\nu$ in the almost sure rate proof. Finally, $\gamma\left(1+\alpha\right)>1$ makes the accumulated Taylor remainder negligible at the same scale.

The matrix $V_H$ is the unique symmetric solution of $HV_H+V_HH=\Sigma$. With normalization $n^{\gamma/2}$, the limiting distribution of the unaveraged iterate, conditionally on $c_\gamma$, is $\mathcal N_d\left(0,c_\gamma V_H\right)$. The limiting covariance $H^{-1}\Sigma H^{-1}$ of $\sqrt n\left(\bar\mu_n-\mu^*\right)$ is the same conditionally on $c_\gamma$ and unconditionally, and does not depend on the multiplicative factor $c_\gamma$ in the step sizes.

For the almost sure convergence rate of $\mu_n$ alone, Assumption~\ref{ass:taylor} is used only through the twice differentiability of $G$ at $\mu^*$ and the positive definiteness of its Hessian; see the proof in Appendix~\ref{sec:asymptotics}.

\begin{remark}[Affine transformations]
Suppose that $\ell=d$ and that the assumptions for $\bar\mu_n$
in Theorem~\ref{th:main} hold. Let $A$ be an invertible
$d\times d$ matrix, let $b\in\R^d$, and set $c\coloneqq|\det A|$.
The simplex-volume identity gives
\[
 g\left(Ax_1+b,\ldots,Ax_d+b;A\mu+b\right)
 =c\,g\left(x_1,\ldots,x_d;\mu\right).
\]
By the chain rule and the zero subgradient convention, the selected
gradient at the transformed points equals
$c\left(A^{-1}\right)^\top\nabla_\mu g(\mathbf x;\mu)$.
Since $c\left(A^{-1}\right)^\top$ generally differs from $A$,
the recursion \eqref{eq:recursion} with the same step sizes is not
affine-equivariant.

For the transformed law, the population minimizer is $A\mu^*+b$,
and the matrices $H$ and $\Sigma$ become, respectively,
\[
 c\left(A^{-1}\right)^\top HA^{-1}
 \quad\text{and}\quad
 c^2\left(A^{-1}\right)^\top\Sigma A^{-1}.
\]
Assumptions~\ref{ass:density} and~\ref{ass:taylor} are preserved
with the same exponents $\nu$ and $\alpha$, by change of variables
and equivalence of norms. Theorem~\ref{th:main} therefore applies
to the recursion based on $(AX_i+b)_{i\geq1}$, initialized at
$A\mu_0+b$ with the same step sizes. Its average is asymptotically
normal about $A\mu^*+b$ at the $\sqrt n$ scale, with limiting
covariance
\[
 A\left(H^{-1}\Sigma H^{-1}\right)A^\top.
\]
\end{remark}

The next proposition gives the convergence rates and asymptotic
distributions under the alternative assumptions for $\ell<d$.

\begin{proposition}
\label{cor:nonoja-alternative}
Suppose that $1\leq\ell<d$ and that Assumptions~\ref{ass:consistency}, \ref{ass:taylor} and~\ref{ass:moment-alternative} hold, with $\gamma q>2$, where $q$ is the moment exponent in Assumption~\ref{ass:moment-alternative}. If $\ell\geq2$, assume in addition that Assumption~\ref{ass:inverse-alternative} holds.
Then the almost sure rate and asymptotic normality of $\mu_n$ stated in Theorem~\ref{th:main} hold. If, in addition, $\gamma\left(1+\alpha\right)>1$, the assertions of that theorem for $\bar\mu_n$ also hold.
\end{proposition}

\subsection{Comparison with offline and disjoint-block estimators}

We compare the averaged sliding-window estimator $\bar\mu_n$
defined in \eqref{eq:averaged-estimator} with the offline
estimator $\widehat\mu_n$, a measurable minimizer of
$\widehat G_n$ in \eqref{eq:empirical-criterion}, and an averaged
stochastic gradient estimator based on disjoint blocks of size $\ell$. The covariances are compared at the $\sqrt n$ scale, where $n$ is the number of observations. We set
\begin{equation*}
 \zeta_1\coloneqq \Var\left(
 \E\left[\nabla_\mu g\left(X_1,\ldots,X_\ell;\mu^*\right)
 \mid X_1\right]\right).
\end{equation*}

For $p=1$ and $\ell<d$, \citet[Theorem~2.5]{durre_affine-equivariant_2022} establish asymptotic normality for the offline estimator, under moment and geometric conditions that do not require a density. The following theorem establishes its asymptotic normality under the assumptions used for the recursion, including the Oja case $\ell=d$.

\begin{theorem}\label{th:batch}
Suppose that Assumption~\ref{ass:taylor} holds.
Assume either that Assumption~\ref{ass:density} holds,
or that $\ell<d$ and Assumptions~\ref{ass:consistency}
and~\ref{ass:moment-alternative} hold.
Then the offline estimator $\widehat\mu_n$, defined as a
measurable minimizer of $\widehat G_n$ in
\eqref{eq:empirical-criterion}, satisfies
\begin{equation}\label{eq:batch-clt}
 \sqrt n\left(\widehat\mu_n-\mu^*\right)
 \xrightarrow[n\to\infty]{\hspace{.5em}\mathcal L\hspace{.5em}}
 \mathcal N_d\left(0,\ell^2H^{-1}\zeta_1H^{-1}\right).
\end{equation}
\end{theorem}

For comparison, we set $\mu_0^{\mathrm{SGD}}\coloneqq \mu_0$ and define, for $j\geq1$,
\begin{equation*}
 \mu_j^{\mathrm{SGD}}
 \coloneqq \mu_{j-1}^{\mathrm{SGD}}
 -\gamma_j\nabla_\mu g\left(
 X_{(j-1)\ell+1},\ldots,X_{j\ell};
 \mu_{j-1}^{\mathrm{SGD}}\right),
\end{equation*}
where $\gamma_j\coloneqq c_\gamma j^{-\gamma}$ is indexed by the number of blocks, with the same $c_\gamma$ and fixed exponent $\gamma$ as in \eqref{eq:recursion}. After $n\geq\ell$ observations, this recursion has performed $\lfloor n/\ell\rfloor$ updates, and its averaged estimator is
\begin{equation*}
 \bar\mu_n^{\mathrm{ASGD}}
 \coloneqq \frac1{\lfloor n/\ell\rfloor}
 \sum_{j=1}^{\lfloor n/\ell\rfloor}\mu_j^{\mathrm{SGD}}.
\end{equation*}

Under Assumptions~\ref{ass:density} and~\ref{ass:taylor}, with $\gamma\nu>2$ and $\gamma\left(1+\alpha\right)>1$, independence of the blocks makes each gradient conditionally unbiased given the preceding blocks and $\left(\mu_0,c_\gamma\right)$. The averaging argument used in the proof of Theorem~\ref{th:main} therefore gives limiting covariance $H^{-1}\Gamma_0H^{-1}$ at the $\sqrt{\lfloor n/\ell\rfloor}$ scale. Since $n/\lfloor n/\ell\rfloor\to\ell$, it follows that
\begin{equation*}
 \sqrt n\left(\bar\mu_n^{\mathrm{ASGD}}-\mu^*\right)
 \xrightarrow[n\to\infty]{\hspace{.5em}\mathcal L\hspace{.5em}}
 \mathcal N_d\left(0,\ell H^{-1}\Gamma_0H^{-1}\right).
\end{equation*} The same limiting covariance holds under the assumptions for $\bar\mu_n$ in Proposition~\ref{cor:nonoja-alternative}, by the averaging argument with independent blocks.
At this same scale, Theorems~\ref{th:batch} and~\ref{th:main} give the covariances $\ell^2H^{-1}\zeta_1H^{-1}$ and $H^{-1}\Sigma H^{-1}$ for $\widehat\mu_n$ and $\bar\mu_n$, respectively. Their comparison follows from the next lemma.

\begin{lemma}
\label{prop:covariance-comparison}\label{lem:loewner_bounds}
Under Assumptions~\ref{ass:consistency} and~\ref{ass:taylor},
\begin{equation*}
 \ell^2\zeta_1\preceq\Sigma\preceq\ell\Gamma_0,
\end{equation*}
where $A\preceq B$ means that $B-A$ is positive semidefinite; we write $A\prec B$ when $B-A$ is positive definite.
If, in addition, Assumption~\ref{ass:density} holds and $\ell\geq2$, then both inequalities are strict.
\end{lemma}

The proof is given in Appendix~\ref{sec:covariance-proof}.
Let $V_{\mathrm{offline}}$, $V_{\mathrm{SW}}$ and
$V_{\mathrm{ASGD}}$ denote the asymptotic covariance matrices
of $\widehat\mu_n$, $\bar\mu_n$ and
$\bar\mu_n^{\mathrm{ASGD}}$, respectively, under the common
$\sqrt n$ normalization. Consequently,
\begin{equation*}
 \underbrace{\ell^2H^{-1}\zeta_1H^{-1}}_{V_{\mathrm{offline}}}
 \preceq
 \underbrace{H^{-1}\Sigma H^{-1}}_{V_{\mathrm{SW}}}
 \preceq
 \underbrace{\ell H^{-1}\Gamma_0H^{-1}}_{V_{\mathrm{ASGD}}}.
\end{equation*}
The averaged sliding-window estimator thus provides a compromise between computational cost and asymptotic variance. Like the disjoint-block estimator, it has computational cost linear in $n$ and storage independent of $n$, for fixed $d$ and $\ell$; see Section~\ref{sec:computation}. Under Assumption~\ref{ass:density}, for $\ell\geq2$, reusing observations in overlapping windows gives a strictly smaller asymptotic variance in every direction than that of the disjoint-block estimator, without minimizing the offline empirical objective function.

\begin{remark}
For $\ell=1$, the sliding-window and disjoint-block recursions give $\mu_n=\mu_n^{\mathrm{SGD}}$ and $\bar\mu_n=\bar\mu_n^{\mathrm{ASGD}}$ for every $n\geq1$. Since $\zeta_1=\Gamma_0=\Sigma$, these averaged estimators have the same limiting covariance as the empirical spatial median.

For $\ell\geq2$, the matrices $V_{\mathrm{offline}}$,
$V_{\mathrm{SW}}$ and $V_{\mathrm{ASGD}}$, defined by the
formulas above, are pairwise distinct under
Assumptions~\ref{ass:consistency} and~\ref{ass:taylor}. More precisely, for every $v\in\R^d$ such that $v^\top V_{\mathrm{ASGD}}v>0$,
\begin{equation*}
 v^\top V_{\mathrm{offline}}v
 <
 v^\top V_{\mathrm{SW}}v
 <
 v^\top V_{\mathrm{ASGD}}v.
\end{equation*}
Moreover, their differences satisfy
\begin{equation*}
 \frac1\ell\left(V_{\mathrm{ASGD}}-V_{\mathrm{offline}}\right)
 \preceq
 V_{\mathrm{SW}}-V_{\mathrm{offline}}
 \preceq
 \frac{2\ell-1}{3\ell}
 \left(V_{\mathrm{ASGD}}-V_{\mathrm{offline}}\right).
\end{equation*}
In particular, for $\ell=2$,
$V_{\mathrm{SW}}=\frac12\left(V_{\mathrm{offline}}+V_{\mathrm{ASGD}}\right)$.
The assertions of this remark are established in the proof of Lemma~\ref{prop:covariance-comparison}.
\end{remark}

We also consider a recursion that uses each disjoint block for a fixed
number $r\geq1$ of successive updates, reevaluates the gradient after
each update, and averages all iterates.
Repeated gradient updates on the same block are also studied by
\citet{agarwal2020stochastic}, who establish excess-risk bounds
for smooth convex objectives.
Under the assumptions for the
averaged estimator in Theorem~\ref{th:main}, its limiting covariance is
still $V_{\mathrm{ASGD}}$, by Proposition~\ref{prop:repeated-blocks}.
Taking $r=\ell$ gives $\ell\lfloor n/\ell\rfloor$ gradient evaluations,
compared with $n-\ell+1$ for sliding windows; both numbers are
asymptotically equivalent to $n$. The strict improvement for $\ell\geq2$
therefore comes from using different tuples, rather than merely
repeating updates on the same tuples.

\subsection{Computational cost}\label{sec:computation}

We now describe the storage requirements and the arithmetic cost of the sliding-window recursion \eqref{eq:recursion}. After $n\geq\ell$ observations, it has performed $n-\ell+1$ gradient evaluations. Storing the current window requires $O\left(d\ell\right)$ memory. The current iterate and its average require $O\left(d\right)$ additional memory, and averaging costs $O\left(d\right)$ operations per update.

For $2\leq\ell\leq d$, the gradient can be computed from the Gram matrix of the $\ell-1$ difference vectors defining the base simplex. This matrix is nonsingular whenever the observations in the window are affinely independent, as they are almost surely under Assumption~\ref{ass:density}, or under Assumption~\ref{ass:moment-alternative} when $\ell<d$. Replacing one observation changes the Gram matrix by a matrix of rank at most two. Updating its inverse and determinant gives $O\left(d\ell\right)$ operations per update after an initial calculation costing $O\left(d\ell^2\right)$. The total cost is therefore $O\left(nd\ell\right)$, with $O\left(d\ell\right)$ storage. For $\ell=1$, each gradient evaluation costs $O\left(d\right)$ operations, giving the same bounds. The formulas and operation counts are given in Appendix~\ref{sec:gradient-computation}.

The disjoint-block recursion defined in the preceding subsection performs $\lfloor n/\ell\rfloor$ updates. Computing each block separately costs $O\left(d\ell^2\right)$ operations, giving a total cost of $O\left(nd\ell\right)$. Thus both recursions have the same order of computational cost, while the sliding-window recursion updates the estimate after each new observation. For the Oja median, $\ell=d$, these bounds give $O\left(nd^2\right)$ operations and $O\left(d^2\right)$ storage.

In contrast, the empirical objective function \eqref{eq:empirical-criterion} of the offline estimator contains $\binom n\ell$ simplex volumes. Computing these volumes separately from their Gram matrices requires $O\left(d\ell^2\binom n\ell\right)$ operations for a single evaluation of this function, and numerical minimization generally requires repeated evaluations. For the Oja median, the number of terms is $\binom nd$, so direct enumeration can become costly even for moderate sample sizes and dimensions. An incomplete offline estimator reduces the number of terms by minimizing an objective function based on a randomly selected collection of tuples \citep[Section~2.4]{durre_affine-equivariant_2022}. The simulations in Section~\ref{sec::simus} use this incomplete version.

\section{Simulation study}\label{sec::simus}

This section evaluates the finite-sample performance and computational efficiency of the proposed sliding-window estimator. While the framework of Section~\ref{sec::algo} accommodates every $1\leq\ell\leq d$, the experiments focus on the Oja median, corresponding to $\ell=d$. We consider distributions that are spherically symmetric about $\mu^*=0$. The simulations are implemented in Python, and the main numerical routines are compiled with Numba before timing. The implementation updates the determinant of the matrix of observations and uses the Sherman--Morrison formula to update its inverse; see Appendix~\ref{sec:oja-implementation}.

\subsection{Online convergence and averaging}

We compare the stochastic gradient estimator based on disjoint blocks (SGD) and the sliding-window estimator (SW), together with their averaged versions, denoted by ASGD and A-SW. SGD is updated after each block of $d$ observations, whereas SW is updated after each new observation, starting from the $d$th.

For each $d\in\{3,10,50\}$, we generate 20 independent samples of $10^7$ observations from $\mathcal N_d(0,I_d)$. All estimators are initialized at a vector of ones, and each sample is processed in a single pass. The sliding-window algorithm uses step size $c_\gamma n^{-\gamma}$ at observation $X_n$, for $n\geq d$, with $\gamma=0.75$ and $c_\gamma\coloneqq\sqrt{d!}$. This choice of $c_\gamma$ gives $c_\gamma^2\Gamma_0=I_d$, since the gradient at $\mu^*=0$ has covariance $\Gamma_0=I_d/d!$; see Appendix~\ref{sec:gaussian-calibration}. The disjoint-block algorithm uses the same sequence of step sizes, indexed by the number of updates.

\begin{figure}[htbp]
    \centering
    \includegraphics[width=1.0\textwidth]{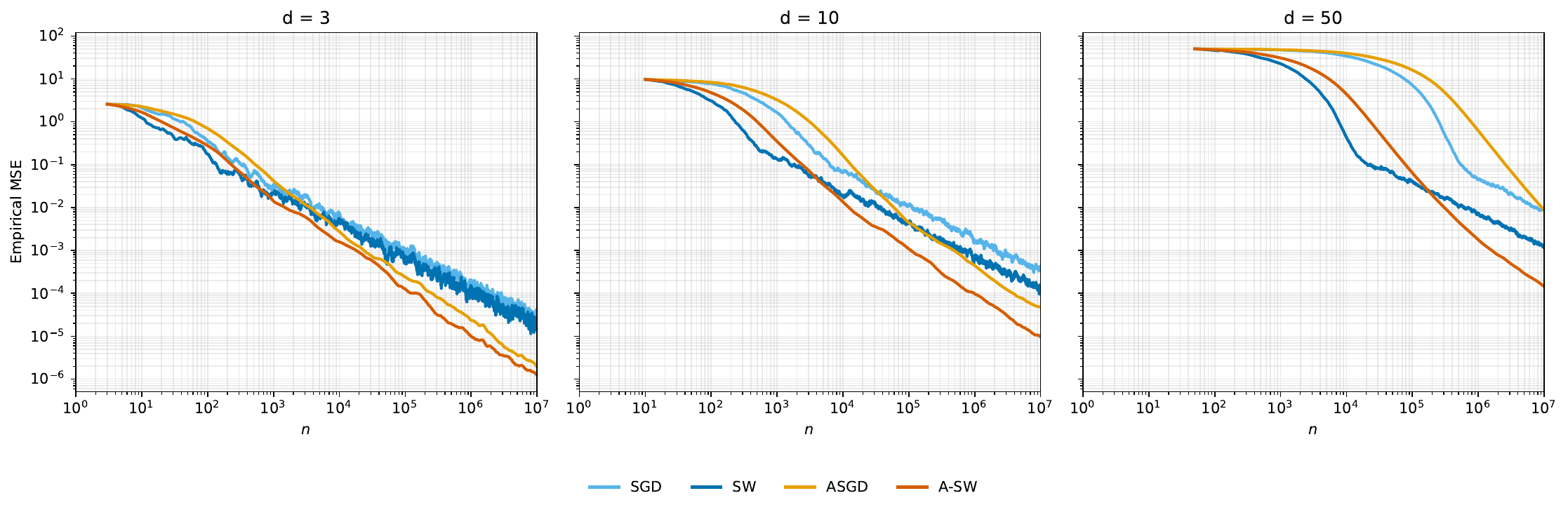}
    \caption{Mean squared errors estimated from 20 independent samples, as a function of the number $n$ of observations, for $d\in\{3,10,50\}$.}
    \label{fig:online_mse}
\end{figure}

Figure~\ref{fig:online_mse} shows the mean squared errors estimated from these samples. After $10^7$ observations, A-SW has the smallest error for each of the three values of $d$. Averaging improves the final error of SW in all three cases. For disjoint blocks, it improves the final error when $d=3$ and $d=10$, whereas SGD and ASGD give similar errors when $d=50$. For $d=50$, the disjoint-block algorithm performs $2\cdot10^5$ updates; the influence of the early iterates may still limit the benefit of averaging.

\subsection{Comparison with an incomplete offline estimator}

We compare our sliding-window estimator and its averaged version with an incomplete offline estimator on datasets of 2000 observations. As in Table~1 of \citet{durre_affine-equivariant_2022}, we consider dimensions $d\in\{2,3,4,6,10\}$ and the following distributions: the standard Gaussian distribution, standard multivariate Student distributions with degrees of freedom $\nu\in\{3,6,9,12\}$, and power-exponential distributions with shape parameter $\eta\in\{2,3,4,6\}$ and density proportional to $\exp\left(-\norm{x}^{2\eta}/2\right)$. The full sum in \eqref{eq:empirical-criterion} contains $\binom{2000}{d}$ terms. Following the incomplete construction in Section~2.4 of \citet{durre_affine-equivariant_2022}, we replace this sum by a sum of $M\coloneqq4000$ terms drawn uniformly with replacement, keeping the selected terms fixed throughout the optimization. Denoting the corresponding tuples by $\mathbf X^{(s)}$, $1\leq s\leq M$, we minimize
\begin{equation*}
 \sum_{s=1}^M g\left(\mathbf X^{(s)};u\right)
\end{equation*}
over $u\in\R^d$. Each term is the absolute value of an affine function of $u$. We precompute the coefficients of these functions, so that a joint evaluation of the objective function and its gradient costs $O(Md)$ operations.

For our sliding-window estimator and its averaged version, we use step size $c_\gamma n^{-\gamma}$, with $\gamma=0.75$, and choose the initial value $\mu_0$ and the constant $c_\gamma$ from a preliminary sample consisting of the first 100 observations. The initial value $\mu_0$ is the vector of coordinate-wise marginal medians of this sample. We draw 50 tuples uniformly with replacement from the $\binom{100}{d}$ possible tuples of $d$ distinct observations in this sample, evaluate the gradient at $\mu_0$ for each tuple, and set $c_\gamma$ to the inverse of the average gradient norm. The algorithm then processes the remaining 1900 observations, which are independent of the preliminary sample and are indexed from $n=1$.

The optimization uses BFGS with a line search satisfying the Armijo and weak Wolfe conditions, with constants $c_1=10^{-4}$ and $c_2=0.9$, respectively. It is initialized at the vector of coordinate-wise marginal medians of the full dataset. The stopping rule uses a numerical upper bound on the error in the objective value, obtained by duality from the selected tuples. The absolute tolerance is $10^{-10}$ after division of the objective function by the sum of the norms of the precomputed linear coefficients. This tolerance is met by all 4500 final estimates. All configurations use the same solver settings.

Table~\ref{tab:offline_mse} reports the mean squared errors estimated from 100 independent samples. The incomplete offline calculation is labelled ``Batch'' in the tables. Comparing the two recursive estimators, averaging gives a smaller error in 35 of the 45 configurations. Its effect is not uniform: for example, the sliding-window estimator has a smaller error than its averaged version under the Student distribution with 3 degrees of freedom in dimension $d=3$. Since the recursion performs fewer than 1900 updates, the contribution of the early iterates may still limit the benefit of averaging.

Table~\ref{tab:offline_time} reports mean total computation times. The sliding-window estimator and its averaged version are computed together, so their shared execution time is reported as ``(A)-SW''. The times include initialization and calibration for the recursive estimators, and initialization, tuple sampling, precomputation and optimization for the offline estimator. For each dimension $d$, since $\ell=d$, the recursive estimators perform the same number of updates for all distributions and have similar mean computation times.

\begin{table}[htbp]
    \centering
    \footnotesize
    \setlength{\tabcolsep}{2.5pt}

    \begin{tabular}{llccccccccc}
\hline
 & & $\nu=3$ & $\nu=6$ & $\nu=9$ & $\nu=12$ & $\mathcal{N}$ & $\eta=2$ & $\eta=3$ & $\eta=4$ & $\eta=6$ \\
\hline
\multirow{3}{*}{$d=2$} & Batch & $3.4 \cdot 10^{-3}$ & $2.8 \cdot 10^{-3}$ & $2.5 \cdot 10^{-3}$ & $2.4 \cdot 10^{-3}$ & $2.2 \cdot 10^{-3}$ & $1.1 \cdot 10^{-3}$ & $1.2 \cdot 10^{-3}$ & $9.7 \cdot 10^{-4}$ & $7.9 \cdot 10^{-4}$ \\
 & SW & $7.7 \cdot 10^{-3}$ & $7.0 \cdot 10^{-3}$ & $5.8 \cdot 10^{-3}$ & $6.6 \cdot 10^{-3}$ & $6.3 \cdot 10^{-3}$ & $3.9 \cdot 10^{-3}$ & $3.2 \cdot 10^{-3}$ & $3.7 \cdot 10^{-3}$ & $3.1 \cdot 10^{-3}$ \\
 & A-SW & $6.6 \cdot 10^{-3}$ & $3.9 \cdot 10^{-3}$ & $3.6 \cdot 10^{-3}$ & $4.4 \cdot 10^{-3}$ & $3.4 \cdot 10^{-3}$ & $1.6 \cdot 10^{-3}$ & $1.4 \cdot 10^{-3}$ & $1.3 \cdot 10^{-3}$ & $1.1 \cdot 10^{-3}$ \\
\addlinespace
\multirow{3}{*}{$d=3$} & Batch & $6.4 \cdot 10^{-3}$ & $4.5 \cdot 10^{-3}$ & $3.7 \cdot 10^{-3}$ & $3.9 \cdot 10^{-3}$ & $3.4 \cdot 10^{-3}$ & $1.4 \cdot 10^{-3}$ & $1.1 \cdot 10^{-3}$ & $1.1 \cdot 10^{-3}$ & $1.0 \cdot 10^{-3}$ \\
 & SW & $1.9 \cdot 10^{-2}$ & $1.4 \cdot 10^{-2}$ & $1.2 \cdot 10^{-2}$ & $9.5 \cdot 10^{-3}$ & $9.7 \cdot 10^{-3}$ & $6.8 \cdot 10^{-3}$ & $5.8 \cdot 10^{-3}$ & $5.4 \cdot 10^{-3}$ & $4.4 \cdot 10^{-3}$ \\
 & A-SW & $2.1 \cdot 10^{-2}$ & $9.6 \cdot 10^{-3}$ & $9.0 \cdot 10^{-3}$ & $7.8 \cdot 10^{-3}$ & $7.2 \cdot 10^{-3}$ & $2.7 \cdot 10^{-3}$ & $2.0 \cdot 10^{-3}$ & $2.2 \cdot 10^{-3}$ & $1.7 \cdot 10^{-3}$ \\
\addlinespace
\multirow{3}{*}{$d=4$} & Batch & $1.3 \cdot 10^{-2}$ & $6.9 \cdot 10^{-3}$ & $6.1 \cdot 10^{-3}$ & $5.7 \cdot 10^{-3}$ & $5.0 \cdot 10^{-3}$ & $1.6 \cdot 10^{-3}$ & $1.3 \cdot 10^{-3}$ & $1.3 \cdot 10^{-3}$ & $1.0 \cdot 10^{-3}$ \\
 & SW & $4.7 \cdot 10^{-2}$ & $2.1 \cdot 10^{-2}$ & $1.9 \cdot 10^{-2}$ & $1.6 \cdot 10^{-2}$ & $1.6 \cdot 10^{-2}$ & $7.9 \cdot 10^{-3}$ & $5.5 \cdot 10^{-3}$ & $5.4 \cdot 10^{-3}$ & $6.6 \cdot 10^{-3}$ \\
 & A-SW & $5.5 \cdot 10^{-2}$ & $1.8 \cdot 10^{-2}$ & $1.6 \cdot 10^{-2}$ & $1.5 \cdot 10^{-2}$ & $1.3 \cdot 10^{-2}$ & $3.8 \cdot 10^{-3}$ & $2.7 \cdot 10^{-3}$ & $2.6 \cdot 10^{-3}$ & $2.3 \cdot 10^{-3}$ \\
\addlinespace
\multirow{3}{*}{$d=6$} & Batch & $4.3 \cdot 10^{-2}$ & $1.3 \cdot 10^{-2}$ & $1.1 \cdot 10^{-2}$ & $1.0 \cdot 10^{-2}$ & $7.6 \cdot 10^{-3}$ & $2.0 \cdot 10^{-3}$ & $1.7 \cdot 10^{-3}$ & $1.4 \cdot 10^{-3}$ & $1.1 \cdot 10^{-3}$ \\
 & SW & $6.4 \cdot 10^{0}$ & $4.9 \cdot 10^{-2}$ & $4.0 \cdot 10^{-2}$ & $3.3 \cdot 10^{-2}$ & $2.4 \cdot 10^{-2}$ & $1.1 \cdot 10^{-2}$ & $9.5 \cdot 10^{-3}$ & $7.9 \cdot 10^{-3}$ & $7.1 \cdot 10^{-3}$ \\
 & A-SW & $7.4 \cdot 10^{0}$ & $5.0 \cdot 10^{-2}$ & $3.9 \cdot 10^{-2}$ & $3.4 \cdot 10^{-2}$ & $2.8 \cdot 10^{-2}$ & $6.2 \cdot 10^{-3}$ & $4.6 \cdot 10^{-3}$ & $3.2 \cdot 10^{-3}$ & $3.0 \cdot 10^{-3}$ \\
\addlinespace
\multirow{3}{*}{$d=10$} & Batch & $1.6 \cdot 10^{-1}$ & $4.0 \cdot 10^{-2}$ & $2.4 \cdot 10^{-2}$ & $2.3 \cdot 10^{-2}$ & $1.5 \cdot 10^{-2}$ & $3.2 \cdot 10^{-3}$ & $2.0 \cdot 10^{-3}$ & $1.6 \cdot 10^{-3}$ & $1.4 \cdot 10^{-3}$ \\
 & SW & $2.4 \cdot 10^{0}$ & $1.4 \cdot 10^{-1}$ & $1.2 \cdot 10^{-1}$ & $1.8 \cdot 10^{-1}$ & $5.5 \cdot 10^{-2}$ & $1.8 \cdot 10^{-2}$ & $1.4 \cdot 10^{-2}$ & $1.3 \cdot 10^{-2}$ & $1.1 \cdot 10^{-2}$ \\
 & A-SW & $3.9 \cdot 10^{0}$ & $1.6 \cdot 10^{-1}$ & $1.0 \cdot 10^{-1}$ & $2.9 \cdot 10^{-1}$ & $6.3 \cdot 10^{-2}$ & $1.1 \cdot 10^{-2}$ & $7.2 \cdot 10^{-3}$ & $5.7 \cdot 10^{-3}$ & $5.0 \cdot 10^{-3}$ \\
\hline
\end{tabular}

    \caption{Mean squared errors estimated from 100 independent samples for the incomplete offline calculation (``Batch''), SW and A-SW, for the dimensions and distributions considered.}
    \label{tab:offline_mse}
\end{table}

\begin{table}[htbp]
    \centering
    \footnotesize
    \setlength{\tabcolsep}{2.5pt}

    \begin{tabular}{llccccccccc}
\hline
 & & $\nu=3$ & $\nu=6$ & $\nu=9$ & $\nu=12$ & $\mathcal{N}$ & $\eta=2$ & $\eta=3$ & $\eta=4$ & $\eta=6$ \\
\hline
\multirow{2}{*}{$d=2$} & Batch & $4.9 \cdot 10^{-2}$ & $4.9 \cdot 10^{-2}$ & $4.8 \cdot 10^{-2}$ & $4.9 \cdot 10^{-2}$ & $4.8 \cdot 10^{-2}$ & $4.8 \cdot 10^{-2}$ & $4.8 \cdot 10^{-2}$ & $4.9 \cdot 10^{-2}$ & $5.3 \cdot 10^{-2}$ \\
 & (A)-SW & $2.9 \cdot 10^{-3}$ & $2.8 \cdot 10^{-3}$ & $2.8 \cdot 10^{-3}$ & $2.9 \cdot 10^{-3}$ & $2.8 \cdot 10^{-3}$ & $2.8 \cdot 10^{-3}$ & $2.8 \cdot 10^{-3}$ & $2.8 \cdot 10^{-3}$ & $3.0 \cdot 10^{-3}$ \\
\addlinespace
\multirow{2}{*}{$d=3$} & Batch & $6.1 \cdot 10^{-2}$ & $6.2 \cdot 10^{-2}$ & $6.1 \cdot 10^{-2}$ & $6.1 \cdot 10^{-2}$ & $6.1 \cdot 10^{-2}$ & $5.9 \cdot 10^{-2}$ & $6.0 \cdot 10^{-2}$ & $5.8 \cdot 10^{-2}$ & $5.8 \cdot 10^{-2}$ \\
 & (A)-SW & $3.0 \cdot 10^{-3}$ & $3.0 \cdot 10^{-3}$ & $3.0 \cdot 10^{-3}$ & $3.1 \cdot 10^{-3}$ & $3.0 \cdot 10^{-3}$ & $3.0 \cdot 10^{-3}$ & $3.0 \cdot 10^{-3}$ & $2.9 \cdot 10^{-3}$ & $2.9 \cdot 10^{-3}$ \\
\addlinespace
\multirow{2}{*}{$d=4$} & Batch & $7.2 \cdot 10^{-2}$ & $6.9 \cdot 10^{-2}$ & $6.9 \cdot 10^{-2}$ & $7.0 \cdot 10^{-2}$ & $7.0 \cdot 10^{-2}$ & $6.6 \cdot 10^{-2}$ & $6.7 \cdot 10^{-2}$ & $6.6 \cdot 10^{-2}$ & $6.6 \cdot 10^{-2}$ \\
 & (A)-SW & $2.9 \cdot 10^{-3}$ & $2.9 \cdot 10^{-3}$ & $2.9 \cdot 10^{-3}$ & $3.0 \cdot 10^{-3}$ & $3.1 \cdot 10^{-3}$ & $2.9 \cdot 10^{-3}$ & $3.0 \cdot 10^{-3}$ & $3.1 \cdot 10^{-3}$ & $3.1 \cdot 10^{-3}$ \\
\addlinespace
\multirow{2}{*}{$d=6$} & Batch & $1.0 \cdot 10^{-1}$ & $1.0 \cdot 10^{-1}$ & $1.0 \cdot 10^{-1}$ & $1.0 \cdot 10^{-1}$ & $1.0 \cdot 10^{-1}$ & $9.6 \cdot 10^{-2}$ & $9.5 \cdot 10^{-2}$ & $9.4 \cdot 10^{-2}$ & $9.5 \cdot 10^{-2}$ \\
 & (A)-SW & $3.3 \cdot 10^{-3}$ & $3.3 \cdot 10^{-3}$ & $3.3 \cdot 10^{-3}$ & $3.3 \cdot 10^{-3}$ & $3.4 \cdot 10^{-3}$ & $3.4 \cdot 10^{-3}$ & $3.5 \cdot 10^{-3}$ & $3.4 \cdot 10^{-3}$ & $3.4 \cdot 10^{-3}$ \\
\addlinespace
\multirow{2}{*}{$d=10$} & Batch & $2.2 \cdot 10^{-1}$ & $2.2 \cdot 10^{-1}$ & $2.0 \cdot 10^{-1}$ & $2.0 \cdot 10^{-1}$ & $2.1 \cdot 10^{-1}$ & $1.9 \cdot 10^{-1}$ & $1.7 \cdot 10^{-1}$ & $1.7 \cdot 10^{-1}$ & $1.8 \cdot 10^{-1}$ \\
 & (A)-SW & $4.2 \cdot 10^{-3}$ & $4.2 \cdot 10^{-3}$ & $4.2 \cdot 10^{-3}$ & $4.1 \cdot 10^{-3}$ & $4.2 \cdot 10^{-3}$ & $4.1 \cdot 10^{-3}$ & $4.2 \cdot 10^{-3}$ & $4.1 \cdot 10^{-3}$ & $4.3 \cdot 10^{-3}$ \\
\hline
\end{tabular}

    \caption{Mean total computation times in seconds over 100 independent replications for the incomplete offline calculation (``Batch'') and the joint computation of SW and A-SW.}
    \label{tab:offline_time}
\end{table}
\clearpage

\section*{Conclusion}

We have investigated a sliding-window stochastic gradient estimator for simplex medians, corresponding to $p=1$ and $1\leq\ell\leq d$. The proposed procedure updates the estimate as each new observation becomes available, requiring $O(d\ell)$ operations per update. It offers a favorable trade-off between estimation accuracy and computational cost compared with a naive disjoint-block stochastic gradient algorithm and offline estimators, as supported by both our theoretical results and numerical experiments.

Several questions remain open. A first direction is to establish convergence rates and limit theorems for the Oja median under weaker distributional assumptions. Extending our analysis to the case $p>1$ would also be of interest. Finally, developing recursive estimators of the asymptotic covariance matrix of the averaged sliding-window estimators would enable the construction of asymptotic confidence intervals and confidence ellipsoids \citep{Godichon2019variance,GodichonLu2024}.

\section*{Acknowledgements}
The authors would like to express their sincere gratitude to Davy Paindaveine for insightful and constructive discussions that inspired this work.

\appendix
\section{Proofs for the sliding-window estimators}\label{sec:proofs}\label{sec::proofs}

Throughout Appendix~\ref{sec:proofs} and
Section~\ref{sec:repeated-blocks},
we put $\F_0\coloneqq \sigma\left(\mu_0,c_\gamma\right)$ and
$\F_k\coloneqq \sigma\left(\mu_0,c_\gamma,X_1,\ldots,X_k\right)$
for $k\geq1$.
For the sliding-window recursion and $k\geq\ell$, we write
$Z_k\coloneqq \left(X_{k-\ell+1},\ldots,X_k\right)$.
In the proofs for the recursive estimators, we first work
conditionally on $\F_0$. By independence of the pair
$\left(\mu_0,c_\gamma\right)$ and the sequence $\left(X_i\right)_{i\geq1}$,
a version of the conditional law of $\left(\mu_0,c_\gamma,\left(X_i\right)_{i\geq1}\right)$ given
$\F_0$ is $\delta_{\left(\mu_0,c_\gamma\right)}\otimes P^{\otimes\infty}$,
where $\delta_z$ denotes the Dirac probability measure at $z$.

\subsection{Proof of Theorem~\ref{th:consistency}}\label{sec:consistency}

For every $\mathbf x\in\left(\R^d\right)^\ell$ and $u,v\in\R^d$, convexity, \eqref{eq:base-height} and \eqref{eq:little_g_grad} give
\begin{gather}
 g\left(\mathbf x;v\right)\geq g\left(\mathbf x;u\right)+\scal{\nabla_\mu g\left(\mathbf x;u\right)}{v-u},\qquad
 \norm{\nabla_\mu g\left(\mathbf x;u\right)}\leq\frac{m_{\ell-1}\left(\Simpl\left(\mathbf x\right)\right)}{\ell},\notag\\*
 \abs{g\left(\mathbf x;u\right)-g\left(\mathbf x;v\right)}\leq\frac{m_{\ell-1}\left(\Simpl\left(\mathbf x\right)\right)}{\ell}\norm{u-v}.\label{eq:envelope}
\end{gather}

For $k\geq\ell$ and $u\in\R^d$, we define
\[
 \xi_k\left(u\right)\coloneqq
 g\left(Z_k;u\right)-G\left(u\right).
\]
By \eqref{eq:envelope}, for every $k\geq\ell$ and $u,v\in\R^d$,
\begin{equation}\label{eq:xi-difference}
 \abs{\xi_k\left(u\right)-\xi_k\left(v\right)}
 \leq\frac1\ell\left(
 m_{\ell-1}\left(\Simpl\left(Z_k\right)\right)
 +\E\left[m_{\ell-1}\left(\Simpl\left(\mathbf X\right)\right)\right]
 \right)\norm{u-v}\quad\as
\end{equation}

\begin{proof}
We set $V_n\coloneqq \norm{\mu_n-\mu^*}^2$.
By \eqref{eq:recursion} and \eqref{eq:envelope}, induction gives,
for every $n\geq0$,
$\E\left[V_n\mid\F_0\right]<\infty$ almost surely.
For $n\geq\ell-1$, summing the successive updates gives
\begin{equation}\label{eq:increment-derived}
 \norm{\mu_n-\mu_{n+1-\ell}}
 \leq\sum_{j=\max\left(\ell,n+2-\ell\right)}^n\gamma_j\frac{m_{\ell-1}\left(\Simpl\left(Z_j\right)\right)}{\ell}.
\end{equation}
The lower limit of the sum is $\max\left(\ell,n+2-\ell\right)$ because
the first $\ell$ iterates are equal to $\mu_0$.
For $\ell=1$, the sum is empty and there is no overlap error.
Since the sequence $(X_i)_{i\geq1}$ is independent of $\F_0$ and all windows
have the same distribution, the conditional Cauchy--Schwarz
inequality gives,
for every $j\geq\ell$,
\begin{align*}
 \E\left[m_{\ell-1}\left(\Simpl\left(Z_{n+1}\right)\right)
          m_{\ell-1}\left(\Simpl\left(Z_j\right)\right)\mid\F_0\right]
 &\leq
 \E\left[m_{\ell-1}\left(\Simpl\left(\mathbf X\right)\right)^2\right]\quad\as,\\
 \E\left[m_{\ell-1}\left(\Simpl\left(Z_j\right)\right)\mid\F_0\right]
 &=\E\left[m_{\ell-1}\left(\Simpl\left(\mathbf X\right)\right)\right]\quad\as
\end{align*}
The sum in \eqref{eq:increment-derived} contains at most $\ell-1$
terms, and, for all sufficiently large $n$, each step size $\gamma_j$
in this sum is at most $2^\gamma\gamma_n$.
This implies that there exists a finite constant $C>0$ such that,
for all sufficiently large $n$,
\begin{equation}\label{eq:increment-expectation}
 \begin{aligned}
 &\frac1\ell
 \E\left[
 \left(m_{\ell-1}\left(\Simpl\left(Z_{n+1}\right)\right)
 +\E\left[m_{\ell-1}\left(\Simpl\left(\mathbf X\right)\right)\right]\right)
 \norm{\mu_n-\mu_{n+1-\ell}}
 \mid\F_0\right]\\
 &\qquad\leq
 \frac{
 \E\left[m_{\ell-1}\left(\Simpl\left(\mathbf X\right)\right)^2\right]
 +\left(\E\left[m_{\ell-1}\left(\Simpl\left(\mathbf X\right)\right)\right]\right)^2
 }{\ell^2}
 \sum_{j=\max\left(\ell,n+2-\ell\right)}^n\gamma_j\\
 &\qquad\leq C\gamma_n\quad\as
 \end{aligned}
\end{equation}

\smallskip\noindent\emph{Uniform bound on $\E\left[V_n\mid\F_0\right]$.}\enspace
The subgradient inequality in \eqref{eq:envelope} yields
\begin{equation}\label{eq:V-basic}
 V_{n+1}\leq V_n-2\gamma_{n+1}
 \left(g\left(Z_{n+1};\mu_n\right)-g\left(Z_{n+1};\mu^*\right)\right)
 +\gamma_{n+1}^2\frac{m_{\ell-1}\left(\Simpl\left(Z_{n+1}\right)\right)^2}{\ell^2}.
\end{equation}
The window $Z_{n+1}=\left(X_{n+2-\ell},\ldots,X_{n+1}\right)$
is independent of $\F_{n+1-\ell}$. Hence
\[
 \E\left[g\left(Z_{n+1};\mu_{n+1-\ell}\right)-g\left(Z_{n+1};\mu^*\right)\mid\F_0\right]
 =\E\left[G\left(\mu_{n+1-\ell}\right)-G\left(\mu^*\right)\mid\F_0\right]
 \geq0\quad\as
\]
By \eqref{eq:envelope} and this nonnegativity, for every $n\geq\ell-1$,
\[
 \E\left[g\left(Z_{n+1};\mu_n\right)-g\left(Z_{n+1};\mu^*\right)\mid\F_0\right]
 \geq-\frac1\ell\E\left[
 m_{\ell-1}\left(\Simpl\left(Z_{n+1}\right)\right)
 \norm{\mu_n-\mu_{n+1-\ell}}\mid\F_0\right]\quad\as
\]
Taking conditional expectations given $\F_0$ in \eqref{eq:V-basic}
and using this bound, \eqref{eq:increment-derived} and the preceding
Cauchy--Schwarz bound gives
\begin{equation*}
 \E\left[V_{n+1}\mid\F_0\right]
 \leq\E\left[V_n\mid\F_0\right]
 +\frac{\E\left[m_{\ell-1}\left(\Simpl\left(\mathbf X\right)\right)^2\right]}{\ell^2}
 \left(\gamma_{n+1}^2
 +2\gamma_{n+1}\sum_{j=\max\left(\ell,n+2-\ell\right)}^n\gamma_j\right)\quad\as
\end{equation*}
Using $2\gamma_{n+1}\gamma_j\leq\gamma_{n+1}^2+\gamma_j^2$ and summing over $n$ shows that
$\sup_n\E\left[V_n\mid\F_0\right]\leq V_0+Kc_\gamma^2$
almost surely, where the constant $K$ is defined by
\begin{equation*}
 K\coloneqq\frac{2\ell-1}{\ell^2}
 \E\left[m_{\ell-1}\left(\Simpl\left(\mathbf X\right)\right)^2\right]
 \sum_{k=\ell}^{\infty}k^{-2\gamma}<\infty.
\end{equation*}
Indeed, each squared step size occurs in at most $2\left(\ell-1\right)$ of the terms obtained from the double sum.

\smallskip\noindent
\emph{Decomposition of
$g\left(Z_{n+1};\mu_n\right)-g\left(Z_{n+1};\mu^*\right)$.}\enspace
For $n\geq\ell-1$, we have
\[
 g\left(Z_{n+1};\mu_n\right)-g\left(Z_{n+1};\mu^*\right)
 =G\left(\mu_n\right)-G\left(\mu^*\right)
 +\varepsilon_{n+1}+\rho_{n+1},
\]
where
\begin{align*}
 \varepsilon_{n+1}&\coloneqq
 \xi_{n+1}\left(\mu_{n+1-\ell}\right)
 -\xi_{n+1}\left(\mu^*\right),\\
 \rho_{n+1}&\coloneqq
 \xi_{n+1}\left(\mu_n\right)
 -\xi_{n+1}\left(\mu_{n+1-\ell}\right).
\end{align*}
Then $\E\left[\varepsilon_{n+1}\mid\F_{n+1-\ell}\right]=0$ almost surely.
By \eqref{eq:xi-difference} and independence of $Z_{n+1}$ and
$\F_{n+1-\ell}$, there exists
a finite constant $C_\varepsilon>0$ such that
\[
 \E\left[\varepsilon_{n+1}^2\mid\F_{n+1-\ell}\right]
 \leq C_\varepsilon V_{n+1-\ell}\quad\as
\]
Taking conditional expectations given $\F_0$ and using the
uniform bound on $\E\left[V_n\mid\F_0\right]$ yields
\[
 \sup_{n\geq\ell-1}\E\left[\varepsilon_{n+1}^2\mid\F_0\right]
 \leq C_\varepsilon\left(V_0+Kc_\gamma^2\right)<\infty\quad\as
\]
Consequently,
\[
 \sum_{n\geq\ell-1}
 \E\left[\left(\gamma_{n+1}\varepsilon_{n+1}\right)^2\mid\F_0\right]
 \leq C_\varepsilon c_\gamma^2\left(V_0+Kc_\gamma^2\right)
 \sum_{k\geq\ell}k^{-2\gamma}<\infty\quad\as
\]
Applying Lemma~\ref{lem:series}(i) conditionally on $\F_0$,
with delay $\ell$, gives
\[
 \Prob\left[\sum_{n\geq\ell-1}\gamma_{n+1}\varepsilon_{n+1}
 \text{ converges}\mid\F_0\right]=1\quad\as
\]

Equations~\eqref{eq:xi-difference} and~\eqref{eq:increment-expectation},
Assumption~\ref{ass:consistency} and
$\sum_{k\geq\ell}\gamma_k^2<\infty$ imply
\begin{equation*}
 \sum_{n\geq\ell-1}
 \E\left[\gamma_{n+1}|\rho_{n+1}|\mid\F_0\right]<\infty
 \quad\text{and}\quad
 \sum_{n\geq\ell-1}
 \E\left[
 \gamma_{n+1}^2
 \frac{m_{\ell-1}\left(\Simpl\left(Z_{n+1}\right)\right)^2}{\ell^2}
 \mid\F_0\right]<\infty\quad\as
\end{equation*}
Since the summands are nonnegative, it follows that
\begin{equation*}
 \sum_{n\geq\ell-1}\gamma_{n+1}|\rho_{n+1}|<\infty
 \quad\as,
 \qquad\text{and}\qquad
 \sum_{n\geq\ell-1}\gamma_{n+1}^2
 \frac{m_{\ell-1}\left(\Simpl\left(Z_{n+1}\right)\right)^2}{\ell^2}
 <\infty\quad\as
\end{equation*}

\smallskip\noindent
\emph{Almost sure convergence of $\left(\mu_n\right)_n$.}\enspace
For $n\geq\ell$, we set
\[
 a_n\coloneqq
 -2\gamma_n\left(\varepsilon_n+\rho_n\right)
 +\gamma_n^2
 \frac{m_{\ell-1}\left(\Simpl\left(Z_n\right)\right)^2}{\ell^2}.
\]
Substituting the decomposition above into \eqref{eq:V-basic} gives
\[
 V_n\leq V_{n-1}
 -2\gamma_n\left(G\left(\mu_{n-1}\right)-G\left(\mu^*\right)\right)
 +a_n.
\]
The convergence of the three series above gives
\[
 \Prob\left[\sum_{n\geq\ell}a_n\text{ converges}\mid\F_0\right]
 =1\quad\as
\]
Since $G(\mu_{n-1})\geq G(\mu^*)$, the inequality for $V_n$
shows that
$\left(V_n-\sum_{k=\ell}^n a_k\right)_{n\geq\ell}$ is non-increasing.
The convergence of $\sum_{k\geq\ell}a_k$ and $V_n\geq0$
show that this sequence is bounded below. Consequently,
$\left(V_n\right)_n$ converges to a finite limit almost surely.
Summing the inequality for $V_n$ and using $V_n\geq0$ gives
$\sum_{n\geq\ell}\gamma_n
\left(G\left(\mu_{n-1}\right)-G\left(\mu^*\right)\right)<\infty$
almost surely.
Since $\sum_{n\geq\ell}\gamma_n=\infty$ almost surely,
we obtain
\[
 \liminf_{n\to\infty}
 \left(G\left(\mu_n\right)-G\left(\mu^*\right)\right)=0\quad\as
\]
The convergence of $\left(V_n\right)_n$ implies that
$\left(\mu_n\right)_n$ is bounded almost surely.
By continuity of $G$ and uniqueness of its minimizer, a subsequence
converges to $\mu^*$ almost surely.
Thus $\liminf_{n\to\infty}V_n=0$ almost surely, and the convergence
of $\left(V_n\right)_n$ gives $V_n\to0$ almost surely.
We have proved
\[
 \sup_n\E\left[V_n\mid\F_0\right]
 \leq V_0+Kc_\gamma^2
 \quad\text{and}\quad
 \Prob\left[\mu_n\xrightarrow[n\to\infty]{\hspace{1.5em}}\mu^*\mid\F_0\right]=1
 \quad\as
\]
Taking expectations gives
\[
 \sup_n\E\left[V_n\right]
 \leq\E\left[V_0\right]+K\E\left[c_\gamma^2\right]<\infty
 \quad\text{and}\quad
 \Prob\left[\mu_n\xrightarrow[n\to\infty]{\hspace{1.5em}}\mu^*\right]=1.
\]
This proves Theorem~\ref{th:consistency}.
\end{proof}

\subsection{Proof of Theorem~\ref{th:main}}\label{sec:asymptotics}
The proof of Theorem~\ref{th:main} is divided into three parts.
For the functions $\xi_k$ introduced in Appendix~\ref{sec:consistency},
and for $k\geq\ell$ and $u\in\R^d$, we use the notation
\[
 \nabla\xi_k\left(u\right)\coloneqq
 \nabla_\mu g\left(Z_k;u\right)-\nabla G\left(u\right),
\]
with the convention $\nabla_\mu g\left(Z_k;u\right)=0$ at points where
$g\left(Z_k;\cdot\right)$ is not differentiable, as in
\eqref{eq:little_g_grad}.
Independence of $Z_k$ and $\F_{k-\ell}$ gives
$\E\left[\nabla\xi_k\left(\mu_{k-\ell}\right)
\mid\F_{k-\ell}\right]=0$ almost surely.
For $k\geq2\ell-1$, the recursion \eqref{eq:recursion} can be written as
\[
 \mu_k=\mu_{k-1}-\gamma_k\left(
 \nabla G\left(\mu_{k-1}\right)
 +\left(\nabla\xi_k\left(\mu_{k-1}\right)-R_k\right)+R_k
 \right),
\]
where
\[
 R_k\coloneqq
 \E\left[\nabla\xi_k\left(\mu_{k-1}\right)
 \mid\F_{k-2\ell+1}\right].
\]
The vector $R_k$ is $\F_{k-2\ell+1}$-measurable and
\[
 \E\left[\nabla\xi_k\left(\mu_{k-1}\right)-R_k
 \mid\F_{k-2\ell+1}\right]=0\quad\as
\]
Since $\nabla G\left(\mu^*\right)=0$, for every $k\geq\ell$,
$\nabla\xi_k\left(\mu^*\right)=\nabla_\mu g\left(Z_k;\mu^*\right)$.
For $k\geq0$, we define the Taylor remainder by
\[
 \delta_k\coloneqq\nabla G\left(\mu_k\right)
 -H\left(\mu_k-\mu^*\right).
\]
For $k\geq2\ell-1$, the recursion becomes
\begin{equation*}
 \mu_k-\mu^*
 =\left(\Id-\gamma_kH\right)\left(\mu_{k-1}-\mu^*\right)
 -\gamma_k\left(\left(\nabla\xi_k\left(\mu_{k-1}\right)-R_k\right)+R_k+\delta_{k-1}\right).
\end{equation*}
Let $n_H$ be the smallest positive integer such that
$\gamma_{n_H}\norm H_{\op}\leq1/2$.
This integer is $\F_0$-measurable. Since $\left(\gamma_k\right)_k$
is decreasing, $\gamma_k\norm H_{\op}\leq1/2$ for every $k\geq n_H$.
For $n\geq k\geq n_H-1$, we set
\[
 \Pi_{n,k}\coloneqq \prod_{j=k+1}^n\left(I_d-\gamma_jH\right),
 \quad\text{and}\quad \Pi_{n,n}=I_d.
\]
We set $K_0\coloneqq\max\left\{n_H,2\ell-1\right\}$.
The random integer $K_0$ and the matrices $\Pi_{n,k}$ are
$\F_0$-measurable.
For every $n\geq K_0$, iterating this recursion from $K_0$ gives
\begin{equation}\label{eq:variation}
 \mu_n-\mu^*=\Pi_{n,K_0}\left(\mu_{K_0}-\mu^*\right)
 -\sum_{k=K_0+1}^n\gamma_k\Pi_{n,k}
           \left(\left(\nabla\xi_k\left(\mu_{k-1}\right)-R_k\right)+R_k+\delta_{k-1}\right).
\end{equation}
The factors commute, since they are polynomials in $H$.
Let $\lambda_{\min}(H)>0$ denote the smallest eigenvalue of $H$.
For every $j\geq n_H$, diagonalizing $H$ in an orthonormal basis gives
$\norm{I_d-\gamma_jH}_{\op}=1-\gamma_j\lambda_{\min}(H)$.
Hence, for $n\geq k\geq n_H-1$,
\begin{equation}\label{eq:Phi-bound}
 \norm{\Pi_{n,k}}_{\op}
 \leq\prod_{j=k+1}^n\left(1-\gamma_j\lambda_{\min}(H)\right)
 \leq\exp\left(-\lambda_{\min}(H)\sum_{j=k+1}^n\gamma_j\right)
 \quad\as
\end{equation}

\begin{proof}[Proof of the almost sure convergence rate of $\left(\mu_n\right)_n$ in Theorem~\ref{th:main}]
By Proposition~\ref{prop:population},
Assumption~\ref{ass:density} implies
Assumption~\ref{ass:consistency}.
Theorem~\ref{th:consistency} therefore gives
$\mu_n\to\mu^*$ almost surely.
Since $\gamma\nu>2$, we fix $q$ such that
\begin{equation*}
 \frac2\gamma<q<\nu.
\end{equation*}

Since
$\sum_{j=K_0+1}^n\gamma_j
\sim c_\gamma n^{1-\gamma}/\left(1-\gamma\right)$,
\eqref{eq:Phi-bound} gives
\begin{equation}\label{eq:initial-product-bound}
 \norm{\Pi_{n,K_0}\left(\mu_{K_0}-\mu^*\right)}
 \leq\norm{\mu_{K_0}-\mu^*}
 \exp\left(-\lambda_{\min}(H)\sum_{j=K_0+1}^n\gamma_j\right)
 =O\left(e^{-cn^{1-\gamma}}\right)\quad\as,
\end{equation}
with $c\coloneqq\lambda_{\min}(H)c_\gamma/\left(2\left(1-\gamma\right)\right)>0$.

Suppose first that $\ell\geq2$. Let $C_1>0$ be a constant for which
the product bound for simplex volumes in \eqref{eq:product-bound}
holds, and set $C\coloneqq2^\gamma\left(\ell-1\right)C_1$.
For every $k\geq2\ell-1$, \eqref{eq:increment-derived} gives
\begin{equation}\label{eq:increment-product}
\begin{aligned}
 \norm{\mu_{k-1}-\mu_{k-\ell}}
 &\leq\sum_{j=k-\ell+1}^{k-1}\gamma_j
 \frac{m_{\ell-1}\left(\Simpl\left(Z_j\right)\right)}{\ell}\\
 &\leq C_1\sum_{j=k-\ell+1}^{k-1}\gamma_j
 \prod_{i=j-\ell+1}^{j}\left(1+\norm{X_i}\right)\\
 &\leq2^\gamma C_1\gamma_k\sum_{j=k-\ell+1}^{k-1}
 \prod_{i=k-2\ell+2}^{k}\left(1+\norm{X_i}\right)\\
 &=C\gamma_k\prod_{i=k-2\ell+2}^{k}\left(1+\norm{X_i}\right)\quad\as
\end{aligned}
\end{equation}

The vector $\mu_{k-\ell}$, the step size $\gamma_k$ and the product
$\prod_{i=k-2\ell+2}^{k-\ell}\left(1+\norm{X_i}\right)$ are
$\F_{k-\ell}$-measurable, and $Z_k$ is independent of $\F_{k-\ell}$.
We set
$t\coloneqq C\gamma_k\prod_{i=k-2\ell+2}^{k-\ell}\left(1+\norm{X_i}\right)$,
where $C$ is the constant in \eqref{eq:increment-product}.
By \eqref{eq:increment-product}, $\mu_{k-1}$ belongs almost surely to
the closed ball with centre $\mu_{k-\ell}$ and radius
$t\prod_{i=k-\ell+1}^{k}\left(1+\norm{X_i}\right)$.
We apply the gradient bound \eqref{eq:osc-majorant} of
Lemma~\ref{lem:oscillation} to $Z_k$, with $u=\mu_{k-\ell}$
and $v=\mu_{k-1}$, which gives
\[
\begin{aligned}
 &\norm{\nabla_\mu g\left(Z_k;\mu_{k-1}\right)
 -\nabla_\mu g\left(Z_k;\mu_{k-\ell}\right)}\\
 &\quad\leq
 \begin{cases}
 \displaystyle\frac{2m_{\ell-1}\left(\Simpl\left(Z_k\right)\right)}{\ell}
 \min\!\left(1,
 \frac{t\prod_{i=k-\ell+1}^{k}\left(1+\norm{X_i}\right)}
 {\dist\left(\mu_{k-\ell},\aff\left(Z_k\right)\right)}\right),
 &\text{if }\ell<d,\\[1em]
 \displaystyle\frac{2m_{d-1}\left(\Simpl\left(Z_k\right)\right)}{d}
 \ind_{\left\{\dist\left(\mu_{k-\ell},\aff\left(Z_k\right)\right)
 \leq t\prod_{i=k-\ell+1}^{k}\left(1+\norm{X_i}\right)\right\}},
 &\text{if }\ell=d,
 \end{cases}
 \quad\as
\end{aligned}
\]
The expectation bounds \eqref{eq:osc-expectation-oja} in
Lemma~\ref{lem:oscillation}, after conditioning on $\F_{k-\ell}$,
give a constant $C_2\geq0$ such that
\[
 \E\left[\norm{\nabla_\mu g\left(Z_k;\mu_{k-1}\right)
 -\nabla_\mu g\left(Z_k;\mu_{k-\ell}\right)}\mid\F_{k-\ell}\right]
 \leq C_2t\quad\as
\]
Taking conditional expectations in \eqref{eq:increment-product} gives
\[
 \E\left[\norm{\mu_{k-1}-\mu_{k-\ell}}\mid\F_{k-\ell}\right]
 \leq t\E\left[\prod_{i=k-\ell+1}^k\left(1+\norm{X_i}\right)
 \mid\F_{k-\ell}\right]
 =t\left(\E\left[1+\norm{X_1}\right]\right)^\ell\quad\as
\]
The triangle inequality and the Lipschitz continuity of
$\nabla G$ from Proposition~\ref{prop:population} now give a constant
$C\geq0$ such that, for every $k\geq2\ell-1$,
\[
 \E\left[\norm{\nabla\xi_k\left(\mu_{k-1}\right)
 -\nabla\xi_k\left(\mu_{k-\ell}\right)}\mid\F_{k-\ell}\right]
 \leq C\gamma_k\prod_{i=k-2\ell+2}^{k-\ell}\left(1+\norm{X_i}\right)\quad\as
\]
The product on the right is independent of $\F_{k-2\ell+1}$ and has
expectation $\left(\E\left[1+\norm{X_1}\right]\right)^{\ell-1}<\infty$.
Taking conditional expectations given $\F_{k-2\ell+1}$ therefore
gives a constant $C\geq0$ such that, for every $k\geq2\ell-1$,
\[
 \E\left[\norm{\nabla\xi_k\left(\mu_{k-1}\right)
 -\nabla\xi_k\left(\mu_{k-\ell}\right)}\mid\F_{k-2\ell+1}\right]
 \leq C\gamma_k\quad\as
\]
The conditional centering of $\nabla\xi_k\left(\mu_{k-\ell}\right)$
and the tower property give
\[
 \E\left[\nabla\xi_k\left(\mu_{k-\ell}\right)\mid\F_{k-2\ell+1}\right]=0\quad\as
\]
Thus, by the definition of $R_k$ and conditional Jensen's inequality,
for every $k\geq2\ell-1$,
\begin{equation}\label{eq:overlap-first-moment}
\begin{aligned}
 \norm{R_k}
 &=\norm{\E\left[\nabla\xi_k\left(\mu_{k-1}\right)
 -\nabla\xi_k\left(\mu_{k-\ell}\right)\mid\F_{k-2\ell+1}\right]}\\
 &\leq\E\left[\norm{\nabla\xi_k\left(\mu_{k-1}\right)
 -\nabla\xi_k\left(\mu_{k-\ell}\right)}\mid\F_{k-2\ell+1}\right]
 \leq C\gamma_k\quad\as
\end{aligned}
\end{equation}
For $\ell=1$, these two bounds hold since $\mu_{k-\ell}=\mu_{k-1}$
and $R_k=0$.
By conditional Jensen's inequality, the gradient bound in
\eqref{eq:envelope} and the moment bounds in \eqref{eq:moment-L},
there exists a constant $C\geq0$ such that, for every $k\geq K_0+1$,
\[
 \E\left[\norm{\nabla\xi_k\left(\mu_{k-1}\right)-R_k}^q
 \mid\F_{k-2\ell+1}\right]\leq C\quad\as
\]
Since $q>2/\gamma$, Lemma~\ref{lem:convolution}, applied
conditionally on $\F_0$ with delay $2\ell-1$, gives
\[
 \norm{\sum_{k=K_0+1}^n\gamma_k\Pi_{n,k}
 \left(\nabla\xi_k\left(\mu_{k-1}\right)-R_k\right)}^2
 =O\left(\gamma_n\log n\right)\quad\as
\]
Moreover, the bound \eqref{eq:weights-power} for weighted sums
in Lemma~\ref{lem:deterministic-weights}, with $a=\gamma$ and $b=0$,
yields
\begin{equation}\label{eq:weighted-R-bound}
 \norm{\sum_{k=K_0+1}^n\gamma_k\Pi_{n,k}R_k}
 \leq C\sum_{k=K_0+1}^n\gamma_k^2\norm{\Pi_{n,k}}_{\op}
 =O\left(\gamma_n\right)\quad\as
\end{equation}
These bounds and \eqref{eq:variation} imply that there exists
a finite random variable $C_0\geq0$ such that, for every $n>K_0$,
\[
 \norm{\mu_n-\mu^*}
 \leq C_0\sqrt{\gamma_n\log n}
 +\sum_{k=K_0+1}^n
 \gamma_k\norm{\Pi_{n,k}}_{\op}\norm{\delta_{k-1}}\quad\as
\]

We define
$\kappa_n\coloneqq
\norm{\delta_n}/\norm{\mu_n-\mu^*}$
when $\mu_n\ne\mu^*$, and $\kappa_n\coloneqq0$ otherwise.
Since $\nabla G\left(\mu^*\right)=0$, we have
$\norm{\delta_n}=\kappa_n\norm{\mu_n-\mu^*}$ in both cases.
By Assumption~\ref{ass:taylor}, $\nabla G$ is differentiable
at $\mu^*$ with derivative $H$, so
$\nabla G\left(\mu^*+u\right)-Hu=o\left(\norm{u}\right)$
as $u\to0$.
Since $\mu_n\to\mu^*$ almost surely,
$\kappa_n\to0$ almost surely.
On the event of probability one where the preceding bounds
hold and $\kappa_n\to0$, Lemma~\ref{lem:perturbed-weights}
applies with
\[
 u_n=\norm{\mu_n-\mu^*},
 \qquad
 v_n=(C_0+1)\sqrt{\gamma_n\log n},
 \qquad
 n_0=K_0.
\]
It follows that
\[
 \norm{\mu_n-\mu^*}
 =O\left(\sqrt{\gamma_n\log n}\right)\quad\as,
\]
which proves the asserted almost sure convergence rate.
\end{proof}

\begin{proof}[Proof of the assertions for the average $\bar\mu_n$ in Theorem~\ref{th:main}]
We follow the linearization and averaging argument in the proof
of \citet[Theorem~2]{PolyakJud92}, with additional estimates for
the dependence between overlapping windows.
For the exponent $q$ chosen above, we set
$\beta\coloneqq\left(q-2\right)/\left(q-1\right)$.

\smallskip\noindent\emph{Decomposition of $H\left(\bar\mu_n-\mu^*\right)$.}\enspace
Rearranging and summing \eqref{eq:recursion} gives
\begin{equation}\label{eq:average-decomposition}
 \left(n-\ell+1\right)H\left(\bar\mu_n-\mu^*\right)
 =T_{1,n}-T_{2,n}-T_{3,n},
\end{equation}
where
\begin{align*}
 T_{1,n}&\coloneqq
 \sum_{k=\ell}^n\frac{\mu_k-\mu_{k+1}}{\gamma_{k+1}},
 \qquad
 T_{2,n}\coloneqq
 \sum_{k=\ell}^n\nabla\xi_{k+1}\left(\mu_{k+1-\ell}\right),\\
 T_{3,n}&\coloneqq
 \sum_{k=\ell}^n\left(
 \delta_k+\nabla\xi_{k+1}\left(\mu_k\right)
 -\nabla\xi_{k+1}\left(\mu_{k+1-\ell}\right)\right).
\end{align*}
Summation by parts yields
\begin{equation*}
 T_{1,n}=\frac{\mu_\ell-\mu^*}{\gamma_{\ell+1}}
 -\frac{\mu_{n+1}-\mu^*}{\gamma_{n+1}}
 +\sum_{k=\ell+1}^n\left(\mu_k-\mu^*\right)
 \left(\frac1{\gamma_{k+1}}-\frac1{\gamma_k}\right).
\end{equation*}
By the mean value theorem,
$0\leq\gamma_{k+1}^{-1}-\gamma_k^{-1}\leq c_\gamma^{-1}\gamma k^{\gamma-1}$.
Since $\gamma<1$, this bound and the almost sure convergence rate
of $\left(\mu_n\right)_n$ give
\[
 \norm{T_{1,n}}^2
 =O\left(n^\gamma\log n\right)=o\left(n\right)\quad\as
\]
The bound on $R_k$ in \eqref{eq:overlap-first-moment} and
$1/2<\gamma<1$ give
\[
 \sum_{k=2\ell-1}^n\norm{R_k}
 =O\left(n^{1-\gamma}\right)=o\left(\sqrt n\right)\quad\as
\]
Assumption~\ref{ass:taylor} and the almost sure convergence rate of $\left(\mu_k\right)_k$ give
\[
 \norm{\delta_k}^2=O\left(k^{-\gamma\left(1+\alpha\right)}
 \left(\log k\right)^{1+\alpha}\right)\quad\as
\]
Since $1/2<\gamma\left(1+\alpha\right)/2<1$, it follows that
\begin{equation*}
 \left(\sum_{k=\ell}^n\norm{\delta_k}\right)^2
 =O\left(n^{2-\gamma\left(1+\alpha\right)}
 \left(\log n\right)^{1+\alpha}\right)=o\left(n\right)\quad\as
\end{equation*}

\smallskip\noindent\emph{Approximation of $T_{2,n}$ and $T_{3,n}$.}\enspace
Independence of $Z_k$ and $\F_{k-\ell}$ gives
\[
 \E\left[\nabla\xi_k\left(\mu_{k-\ell}\right)-\nabla\xi_k\left(\mu^*\right)
 \mid\F_{k-\ell}\right]=0\quad\as
\]
By the bound \eqref{eq:L2-modulus} in
Proposition~\ref{prop:population}, there exists a constant $C\geq0$
such that, for every $k\geq\ell$,
\begin{equation}\label{eq:D-variance}
\begin{aligned}
 &\E\left[\norm{\nabla\xi_k\left(\mu_{k-\ell}\right)-\nabla\xi_k\left(\mu^*\right)}^2\mid\F_{k-\ell}\right]\\
 &\quad=\E\left[\norm{\nabla_\mu g\left(Z_k;\mu_{k-\ell}\right)
 -\nabla_\mu g\left(Z_k;\mu^*\right)}^2\mid\F_{k-\ell}\right]
 -\norm{\nabla G\left(\mu_{k-\ell}\right)
 -\nabla G\left(\mu^*\right)}^2\\
 &\quad\leq C\min\left(\norm{\mu_{k-\ell}-\mu^*}^\beta,1\right)
 \quad\as
\end{aligned}
\end{equation}
The almost sure convergence rate of $\left(\mu_k\right)_k$ shows that the
right-hand side is
$O\left(k^{-\gamma\beta/2}\left(\log k\right)^{\beta/2}\right)$
almost surely.
Since $\gamma\beta>0$, it follows that
\[
 \sum_{k=\ell}^{\infty}\frac1k
 \E\left[\norm{\nabla\xi_k\left(\mu_{k-\ell}\right)-\nabla\xi_k\left(\mu^*\right)}^2\mid\F_{k-\ell}\right]
 <\infty\quad\as
\]
For $k\geq2\ell-1$, we define
\[
 \eta_k\coloneqq
 \nabla\xi_k\left(\mu_{k-1}\right)
 -\nabla\xi_k\left(\mu_{k-\ell}\right)-R_k.
\]
The definition of $R_k$ and the tower property give
$\E\left[\eta_k\mid\F_{k-2\ell+1}\right]=0$ almost surely.
By \eqref{eq:envelope}, for every $k\geq2\ell-1$,
\[
 \norm{\nabla\xi_k\left(\mu_{k-1}\right)-\nabla\xi_k\left(\mu_{k-\ell}\right)}
 \leq\frac{2m_{\ell-1}\left(\Simpl\left(Z_k\right)\right)}{\ell}
 +2\E\left[\frac{m_{\ell-1}\left(\Simpl\left(\mathbf X\right)\right)}{\ell}\right]\quad\as
\]
Since $Z_k$ is independent of $\F_{k-2\ell+1}$, the moment bounds in
\eqref{eq:moment-L} give a constant $C_q\geq0$ such that, for every
$k\geq2\ell-1$,
\[
 \E\left[\norm{\nabla\xi_k\left(\mu_{k-1}\right)
 -\nabla\xi_k\left(\mu_{k-\ell}\right)}^q
 \mid\F_{k-2\ell+1}\right]\leq C_q\quad\as
\]
The conditional variance identity, \eqref{eq:overlap-first-moment}
and conditional H\"older's inequality in the form \eqref{eq:interpolation}
given in the proof of Proposition~\ref{prop:population}
therefore give a constant $C\geq0$ such that, for every $k\geq2\ell-1$,
\begin{equation}\label{eq:overlap-variance}
 \E\left[\norm{\eta_k}^2\mid\F_{k-2\ell+1}\right]
 =\E\left[\norm{\nabla\xi_k\left(\mu_{k-1}\right)
 -\nabla\xi_k\left(\mu_{k-\ell}\right)}^2\mid\F_{k-2\ell+1}\right]
 -\norm{R_k}^2
 \leq C\gamma_k^\beta\quad\as
\end{equation}
Since $\gamma_k=c_\gamma k^{-\gamma}$ and $\gamma\beta>0$,
\eqref{eq:overlap-variance} gives
\[
 \sum_{k=2\ell-1}^{\infty}\frac1k
 \E\left[\norm{\eta_k}^2\mid\F_{k-2\ell+1}\right]
 <\infty\quad\as
\]
Applying Lemma~\ref{lem:series}(iii) under the conditional laws
given $\F_0$, with delays $\ell$ and $2\ell-1$, respectively, gives
\begin{equation*}
 \sum_{k=\ell}^n\left(\nabla\xi_k\left(\mu_{k-\ell}\right)-\nabla\xi_k\left(\mu^*\right)\right)
 =o\left(\sqrt n\right)\quad\as,
 \quad\text{and}\quad
 \sum_{k=2\ell-1}^n\eta_k
 =o\left(\sqrt n\right)\quad\as
\end{equation*}
The two estimates above and the bounds on
$\sum_{k=\ell}^n\norm{\delta_k}$ and
$\sum_{k=2\ell-1}^n\norm{R_k}$ give
\begin{equation*}
 T_{3,n}=o\left(\sqrt n\right)\quad\as,
 \quad\text{and}\quad
 T_{2,n}=\sum_{k=\ell}^n\nabla\xi_{k+1}\left(\mu^*\right)
 +o\left(\sqrt n\right)\quad\as
\end{equation*}
Substituting in \eqref{eq:average-decomposition} and multiplying
by $H^{-1}/\sqrt n$ yields
\begin{equation}\label{eq:representation-N}
 \frac{n-\ell+1}{\sqrt n}\left(\bar\mu_n-\mu^*\right)
 =-H^{-1}\frac1{\sqrt n}\sum_{k=\ell}^n\nabla\xi_{k+1}\left(\mu^*\right)+o\left(1\right)\quad\as
\end{equation}

\smallskip\noindent\emph{Asymptotic normality of $\left(\bar\mu_n\right)_n$.}\enspace
The sequence $\left(\nabla\xi_k\left(\mu^*\right)\right)_{k\geq\ell}$ is strictly stationary,
since the same measurable function is applied to each window
of independent observations with common law $P$, and it is
independent of $\F_0$.
By the definition of $\nabla\xi_k$, the gradient formula in
Proposition~\ref{prop:population} and the bound \eqref{eq:envelope},
for every $k\geq\ell$,
\[
 \E\left[\nabla\xi_k\left(\mu^*\right)\right]=0,
 \quad\text{and}\quad
 \E\left[\norm{\nabla\xi_k\left(\mu^*\right)}^2\right]
 \leq\frac1{\ell^2}\E\left[
 m_{\ell-1}\left(\Simpl\left(\mathbf X\right)\right)^2\right]<\infty.
\]
The sequence $\left(\nabla\xi_k\left(\mu^*\right)\right)_{k\geq\ell}$ is
$\left(\ell-1\right)$-dependent: for every integer $a\geq\ell$,
the sequences $\left(\nabla\xi_k\left(\mu^*\right)\right)_{\ell\leq k\leq a}$
and $\left(\nabla\xi_k\left(\mu^*\right)\right)_{k\geq a+\ell}$ depend
respectively on $\left(X_i\right)_{1\leq i\leq a}$ and
$\left(X_i\right)_{i\geq a+1}$, and are therefore independent.
For $\ell\geq2$, the central limit theorem in \citet{Duflo}, recalled in
Proposition~\ref{th:duflo}, applies with $m=\ell-1$.
For $\ell=1$, we use $m=1$; the covariance between consecutive
terms is zero by independence.
For every $k\geq\ell$ and $j\in\{0,\ldots,\ell-1\}$,
the definition of $\Gamma_j$ in Theorem~\ref{th:main} gives
$\E\left[\nabla\xi_k\left(\mu^*\right)
\nabla\xi_{k+j}\left(\mu^*\right)^\top\right]=\Gamma_j$.
Applying Proposition~\ref{th:duflo} conditionally on $\F_0$ to the
sequence $\left(\nabla\xi_{k+1}\left(\mu^*\right)\right)_{k\geq\ell}$
gives weak convergence:
\[
 \mathcal L\left(\frac1{\sqrt{n-\ell+1}}
 \sum_{k=\ell}^n\nabla\xi_{k+1}\left(\mu^*\right)\mid\F_0\right)
 \xrightarrow[n\to\infty]{\hspace{1.5em}}\mathcal N_d\left(0,\Sigma\right)\quad\as
\]
Since $\left(n-\ell+1\right)/n\to1$, \eqref{eq:representation-N}
and Slutsky's theorem prove the conditional asymptotic normality
of $\left(\bar\mu_n\right)_n$, with covariance $H^{-1}\Sigma H^{-1}$.

\smallskip\noindent\emph{Almost sure convergence rate of $\left(\bar\mu_n\right)_n$.}\enspace
For $a\in\{0,\ldots,\ell-1\}$, the sequence
$\left(\nabla\xi_{a+j\ell+1}\left(\mu^*\right)\right)_{j\geq1}$,
indexed by the residue class $a+1$ modulo $\ell$, is independent
and identically distributed, since the corresponding windows are disjoint.
Applying the law of the iterated logarithm in
\citet[Theorem~6.4.23, p.~209]{Duflo} to each coordinate and summing
over the $d$ coordinates and the $\ell$ residue classes gives
\begin{equation*}
 \norm{\sum_{k=\ell}^n\nabla\xi_{k+1}\left(\mu^*\right)}
 =O\left(\sqrt{n\log\log n}\right)\quad\as
\end{equation*}
Equation~\eqref{eq:representation-N} and the preceding bound give
\begin{equation*}
 \frac{n-\ell+1}{\sqrt n}\norm{\bar\mu_n-\mu^*}
 =O\left(\sqrt{\log\log n}\right)+o\left(1\right)
 \quad\as
\end{equation*}
Since $\left(n-\ell+1\right)/n\to1$, it follows that
\begin{equation*}
 \norm{\bar\mu_n-\mu^*}^2
 =O\left(\frac{\log\log n}{n}\right)\quad\as
\end{equation*}
This proves the asserted almost sure convergence rate of
$\left(\bar\mu_n\right)_n$.

\end{proof}

\begin{proof}[Proof of the asymptotic normality of $\left(\mu_n\right)_n$ in Theorem~\ref{th:main}]
\leavevmode\par\nobreak
\smallskip\noindent\emph{Approximation of $\frac1{\sqrt{\gamma_n}}\left(\mu_n-\mu^*\right)$.}\enspace
Adding and subtracting $\nabla\xi_k\left(\mu^*\right)$ in
\eqref{eq:variation} gives, for every $n\geq K_0$,
\begin{equation}\label{eq:raw-expanded-recursion}
\begin{aligned}
 \mu_n-\mu^*
 &=\Pi_{n,K_0}\left(\mu_{K_0}-\mu^*\right)
 -\sum_{k=K_0+1}^n\gamma_k\Pi_{n,k}\nabla\xi_k\left(\mu^*\right)
 -\sum_{k=K_0+1}^n\gamma_k\Pi_{n,k}R_k\\
 &\quad-\sum_{k=K_0+1}^n\gamma_k\Pi_{n,k}
 \left(\nabla\xi_k\left(\mu_{k-1}\right)-R_k-\nabla\xi_k\left(\mu^*\right)\right)
 -\sum_{k=K_0+1}^n\gamma_k\Pi_{n,k}\delta_{k-1}.
\end{aligned}
\end{equation}
The bounds \eqref{eq:initial-product-bound} and
\eqref{eq:weighted-R-bound} give
\[
 \norm{\Pi_{n,K_0}\left(\mu_{K_0}-\mu^*\right)}
 =o\left(\sqrt{\gamma_n}\right)
 \quad\text{and}\quad
 \norm{\sum_{k=K_0+1}^n\gamma_k\Pi_{n,k}R_k}
 =o\left(\sqrt{\gamma_n}\right)\quad\as
\]
By Assumption~\ref{ass:taylor}, the almost sure convergence rate of $\left(\mu_n\right)_n$
and the bound \eqref{eq:weights-power} in
Lemma~\ref{lem:deterministic-weights}, and since $\alpha>0$,
\begin{equation}\label{eq:raw-taylor-bound}
 \norm{\sum_{k=K_0+1}^n\gamma_k\Pi_{n,k}\delta_{k-1}}
 =O\left(\gamma_n^{\left(1+\alpha\right)/2}\left(\log n\right)^{\left(1+\alpha\right)/2}\right)
 =o\left(\sqrt{\gamma_n}\right)\quad\as
\end{equation}

For $k\geq2\ell-1$, we set
\[
 \widetilde\eta_k\coloneqq
 \nabla\xi_k\left(\mu_{k-1}\right)-R_k
 -\nabla\xi_k\left(\mu^*\right).
\]
The definition of $R_k$ and the centering of
$\nabla\xi_k\left(\mu^*\right)$ give
\[
 \E\left[\widetilde\eta_k\mid\F_{k-2\ell+1}\right]=0\quad\as
\]
By the definition of $\eta_k$,
\[
 \widetilde\eta_k
 =\eta_k+\nabla\xi_k\left(\mu_{k-\ell}\right)
 -\nabla\xi_k\left(\mu^*\right).
\]
Conditional dominated convergence applies to
$\min\left(\norm{\mu_{k-\ell}-\mu^*}^\beta,1\right)$, since it is bounded
by one and $\mu_k\to\mu^*$ almost surely.
Using $\norm{x+y}^2\leq2\norm{x}^2+2\norm{y}^2$ and taking conditional
expectations in \eqref{eq:D-variance} and \eqref{eq:overlap-variance}
given $\F_0$, we obtain a constant $C\geq0$ such that
\begin{equation}\label{eq:raw-centered-second-moment}
\begin{aligned}
 \E\left[\norm{\widetilde\eta_k}^2\mid\F_0\right]
 &\leq2\E\left[\norm{\eta_k}^2\mid\F_0\right]
 +2\E\left[\norm{\nabla\xi_k\left(\mu_{k-\ell}\right)
 -\nabla\xi_k\left(\mu^*\right)}^2\mid\F_0\right]\\
 &\leq C\gamma_k^\beta+
 C\E\left[\min\left(\norm{\mu_{k-\ell}-\mu^*}^\beta,1\right)
 \mid\F_0\right]\xrightarrow[k\to\infty]{\hspace{1.5em}}0\quad\as
\end{aligned}
\end{equation}
Let $i,j$ be integers such that $2\ell-1\leq i<j$ and $j-i$ is a
multiple of $2\ell-1$. Then $i\leq j-2\ell+1$, so
$\widetilde\eta_i$
is $\F_{j-2\ell+1}$-measurable. Conditioning first on
$\F_{j-2\ell+1}$ and then on $\F_0$ gives
\[
 \E\left[
 \widetilde\eta_i
 \widetilde\eta_j^\top
 \mid\F_0\right]=0\quad\as
\]
Since $K_0$ and $\gamma_k\Pi_{n,k}$ are $\F_0$-measurable,
splitting into residue classes and applying the Cauchy--Schwarz inequality gives
\begin{align*}
 &\frac1{\gamma_n}\E\left[\norm{\sum_{k=K_0+1}^n\gamma_k\Pi_{n,k}
 \widetilde\eta_k}^2\mid\F_0\right]\\
 &\quad\leq\frac{2\ell-1}{\gamma_n}\sum_{a=0}^{2\ell-2}
 \E\left[\norm{\sum_{\substack{K_0<k\leq n\\k\equiv a\pmod{2\ell-1}}}
 \gamma_k\Pi_{n,k}\widetilde\eta_k}^2\mid\F_0\right]\\
 &\quad=\frac{2\ell-1}{\gamma_n}\sum_{k=K_0+1}^n
 \E\left[\norm{\gamma_k\Pi_{n,k}\widetilde\eta_k}^2\mid\F_0\right]\\
 &\quad\leq\frac{2\ell-1}{\gamma_n}\sum_{k=K_0+1}^n
 \gamma_k^2\norm{\Pi_{n,k}}_{\op}^2
 \E\left[\norm{\widetilde\eta_k}^2\mid\F_0\right]\quad\as
\end{align*}
By the bound \eqref{eq:weights-square} in
Lemma~\ref{lem:deterministic-weights}, there exists a finite
$\F_0$-measurable random variable $C_{\Pi}\geq0$ such that,
for every $n\geq K_0+1$,
\[
 \frac1{\gamma_n}\sum_{k=K_0+1}^n
 \gamma_k^2\norm{\Pi_{n,k}}_{\op}^2\leq C_{\Pi}\quad\as
\]
We fix an integer $K\geq2\ell-1$.
As $n$ tends to infinity with $K$ fixed, the exponential bound
\eqref{eq:Phi-bound} gives
\begin{align*}
 \frac1{\gamma_n}\sum_{k=K_0+1}^{K-1}
 \gamma_k^2\norm{\Pi_{n,k}}_{\op}^2
 \E\left[\norm{\widetilde\eta_k}^2\mid\F_0\right]
 &=O\left(\gamma_n^{-1}e^{-\lambda_{\min}(H)c_\gamma n^{1-\gamma}}\right)\\
 &=o\left(1\right)\quad\as
\end{align*}
For this fixed $K$, splitting at $K$ the sum of nonnegative terms
in the preceding second-moment bound gives
\begin{align*}
 &\frac1{\gamma_n}
 \E\left[\norm{\sum_{k=K_0+1}^n\gamma_k\Pi_{n,k}
 \widetilde\eta_k}^2\mid\F_0\right]\\
 &\quad\leq\frac{2\ell-1}{\gamma_n}
 \sum_{k=K_0+1}^{K-1}\gamma_k^2\norm{\Pi_{n,k}}_{\op}^2
 \E\left[\norm{\widetilde\eta_k}^2\mid\F_0\right]\\
 &\qquad+\frac{2\ell-1}{\gamma_n}
 \sum_{k=\max\{K,K_0+1\}}^n\gamma_k^2\norm{\Pi_{n,k}}_{\op}^2
 \E\left[\norm{\widetilde\eta_k}^2\mid\F_0\right]\\
 &\quad=o\left(1\right)+\frac{2\ell-1}{\gamma_n}
 \sum_{k=\max\{K,K_0+1\}}^n\gamma_k^2\norm{\Pi_{n,k}}_{\op}^2
 \E\left[\norm{\widetilde\eta_k}^2\mid\F_0\right]\\
 &\quad\leq o\left(1\right)+\frac{2\ell-1}{\gamma_n}
 \left(\sum_{k=K_0+1}^n\gamma_k^2\norm{\Pi_{n,k}}_{\op}^2\right)
 \sup_{j\geq K}\E\left[\norm{\widetilde\eta_j}^2\mid\F_0\right]\\
 &\quad\leq o\left(1\right)+(2\ell-1)C_{\Pi}
 \sup_{j\geq K}\E\left[\norm{\widetilde\eta_j}^2\mid\F_0\right]\quad\as
\end{align*}
By \eqref{eq:raw-centered-second-moment}, the supremum tends to zero
almost surely as $K\to\infty$.
We first take the limit superior in $n$, then let $K\to\infty$.
Conditional Markov's inequality gives, for every $\varepsilon>0$,
\begin{equation}\label{eq:raw-gradient-probability}
\begin{aligned}
 &\Prob\left[\norm{\sum_{k=K_0+1}^n\gamma_k\Pi_{n,k}
 \widetilde\eta_k}
 >\varepsilon\sqrt{\gamma_n}\mid\F_0\right]\\
 &\quad\leq\frac1{\varepsilon^2\gamma_n}
 \E\left[\norm{\sum_{k=K_0+1}^n\gamma_k\Pi_{n,k}
 \widetilde\eta_k}^2\mid\F_0\right]
 \xrightarrow[n\to\infty]{\hspace{1.5em}}0\quad\as
\end{aligned}
\end{equation}
Combining \eqref{eq:initial-product-bound}, \eqref{eq:weighted-R-bound},
\eqref{eq:raw-taylor-bound} and \eqref{eq:raw-gradient-probability}
in \eqref{eq:raw-expanded-recursion} gives
\begin{equation}\label{eq:raw-linear-representation}
 \frac1{\sqrt{\gamma_n}}\left(\mu_n-\mu^*\right)
 =-\frac1{\sqrt{\gamma_n}}\sum_{k=K_0+1}^n\gamma_k\Pi_{n,k}\nabla\xi_k\left(\mu^*\right)
 +o_{\Prob(\cdot\mid\F_0)}\left(1\right),
\end{equation}
where $o_{\Prob(\cdot\mid\F_0)}\left(1\right)$ denotes convergence to
zero in probability under almost every conditional law given $\F_0$.

\smallskip\noindent\emph{Limit distribution.}\enspace
The sequence $\left(\nabla\xi_k\left(\mu^*\right)\right)_{k\geq\ell}$ is
strictly stationary, centered, square-integrable and $\left(\ell-1\right)$-dependent,
and is independent of $\F_0$.
With $V_H$ defined in \eqref{eq:Sigma},
Lemma~\ref{lem:matrix-mdep-clt}, applied conditionally on $\F_0$, gives
\[
 \mathcal L\left(\frac1{\sqrt{\gamma_n}}
 \sum_{k=K_0+1}^n\gamma_k\Pi_{n,k}\nabla\xi_k\left(\mu^*\right)
 \mid\F_0\right)
 \xrightarrow[n\to\infty]{\hspace{1.5em}}\mathcal N_d\left(0,V_H\right)\quad\as
\]
Equation~\eqref{eq:raw-linear-representation} and Slutsky's theorem
under the conditional laws prove the conditional asymptotic
normality of $\left(\mu_n\right)_n$ stated in Theorem~\ref{th:main}.

Since $c_\gamma$ is $\F_0$-measurable, the conditional covariance
under normalization by $\sqrt{n^\gamma}$ follows from
$\gamma_n^{-1/2}=c_\gamma^{-1/2}\sqrt{n^\gamma}$.

\smallskip\noindent\emph{Characterization of $V_H$.}\enspace
To prove the characterization stated after Theorem~\ref{th:main}, we integrate
\[
 \frac d{dt}\left(e^{-tH}\Sigma e^{-tH}\right)
 =-He^{-tH}\Sigma e^{-tH}-e^{-tH}\Sigma e^{-tH}H
\]
over $[0,\infty)$ and obtain $HV_H+V_HH=\Sigma$.
In a basis diagonalizing $H$, with corresponding eigenvalues $\lambda_1,\ldots,\lambda_d$, uniqueness follows from
$\left(\lambda_i+\lambda_j\right)\left(V_H\right)_{ij}=\Sigma_{ij}$, since all $\lambda_i+\lambda_j$ are positive.
\end{proof}

\begin{proof}[Completion of the proof of Theorem~\ref{th:main}]
We pass from the conditional laws given $\F_0$ to the original law.
For each event $E$ expressing an almost sure rate assertion in
Theorem~\ref{th:main}, the preceding proofs give
$\Prob[E\mid\F_0]=1$ almost surely;
hence $\Prob[E]=\E[\Prob[E\mid\F_0]]=1$.

For the iterate, we take $Y_n\coloneqq\gamma_n^{-1/2}(\mu_n-\mu^*)$
and $Q\coloneqq\mathcal N_d(0,V_H)$; for the average, under its additional
assumptions, we take $Y_n\coloneqq\sqrt n(\bar\mu_n-\mu^*)$ and
$Q\coloneqq\mathcal N_d(0,H^{-1}\Sigma H^{-1})$.
In both cases, $Q$ is deterministic and, for every bounded continuous
function $\varphi:\R^d\to\R$, the conditional limits give
\[
 \E[\varphi(Y_n)\mid\F_0]
 \xrightarrow[n\to\infty]{\hspace{1.5em}}
 \int_{\R^d}\varphi(x)\,Q(dx)\quad\as
\]
Since $|\E[\varphi(Y_n)\mid\F_0]|\leq\sup_{x\in\R^d}|\varphi(x)|$,
the tower property and dominated convergence, also conditionally on
$c_\gamma$, yield
\begin{align*}
 \E[\varphi(Y_n)]
 &=\E[\E[\varphi(Y_n)\mid\F_0]]
 \xrightarrow[n\to\infty]{\hspace{1.5em}}
 \int_{\R^d}\varphi(x)\,Q(dx),\\
 \E[\varphi(Y_n)\mid c_\gamma]
 &=\E[\E[\varphi(Y_n)\mid\F_0]\mid c_\gamma]
 \xrightarrow[n\to\infty]{\hspace{1.5em}}
 \int_{\R^d}\varphi(x)\,Q(dx)\quad\as
\end{align*}
Thus $(Y_n)_n$ converges in distribution to $Q$ under the original law
and under almost every conditional law given $c_\gamma$.

For the normalization $\sqrt{n^\gamma}$ discussed after
Theorem~\ref{th:main}, the conditional Gaussian limit given $\F_0$
has covariance $c_\gamma V_H$.
Since $\sigma(c_\gamma)\subset\F_0$, conditional dominated convergence
gives
\[
 \E\left[\varphi\left(\sqrt{n^\gamma}(\mu_n-\mu^*)\right)\mid c_\gamma\right]
 \xrightarrow[n\to\infty]{\hspace{1.5em}}
 \int_{\R^d}\varphi(x)\,\mathcal N_d(0,c_\gamma V_H)(dx)\quad\as
\]
Taking expectations gives
\[
 \E\left[\varphi\left(\sqrt{n^\gamma}(\mu_n-\mu^*)\right)\right]
 \xrightarrow[n\to\infty]{\hspace{1.5em}}
 \E\left[\int_{\R^d}\varphi(x)\,\mathcal N_d(0,c_\gamma V_H)(dx)\right],
\]
which identifies the unconditional limit as the mixture of these Gaussian
laws under the distribution of $c_\gamma$.
\end{proof}

\subsection{Proof of Proposition~\ref{cor:nonoja-alternative}}\label{sec:alternative-proof}

\begin{proof}
We retain the notation of the proof of Theorem~\ref{th:main}.
Let $q$ be the exponent in Assumption~\ref{ass:moment-alternative} and set
$\beta\coloneqq(q-2)/(q-1)$.
By Assumption~\ref{ass:consistency} and Theorem~\ref{th:consistency},
$\mu_n\to\mu^*$ almost surely.

\smallskip\noindent\emph{Differentiability of $G$.}\enspace
The moment condition in Assumption~\ref{ass:moment-alternative},
independence of $X_1,\ldots,X_\ell$ and the product bound
\eqref{eq:product-bound} give
$\E[m_{\ell-1}(\Simpl(\mathbf X))^q]<\infty$ and $G(0)<\infty$.
The bound \eqref{eq:envelope} then shows that $G$ is finite on $\R^d$.
For every $u\in\R^d$, Lemma~S.1.3 of~\citet{DPsupp} gives
$\Prob[u\in\aff(\mathbf X)]=0$, so $g(\mathbf X;\cdot)$ is differentiable
at $u$ almost surely. For every nonzero $v\in\R^d$,
\[
 \frac{\left|g(\mathbf X;u+v)-g(\mathbf X;u)
 -\scal{\nabla_\mu g(\mathbf X;u)}{v}\right|}{\norm v}
 \leq\frac{2m_{\ell-1}(\Simpl(\mathbf X))}{\ell}\quad\as
\]
The dominated convergence theorem applies: the right-hand side is integrable,
and differentiability of $g(\mathbf X;\cdot)$ at $u$ gives convergence of the
left-hand side to zero almost surely as $v\to0$. Hence
\[
 \frac{G(u+v)-G(u)-\scal{\E[\nabla_\mu g(\mathbf X;u)]}{v}}{\norm v}
 \xrightarrow[v\to0]{\hspace{1.5em}}0.
\]
Thus $G$ is differentiable and, for every $u\in\R^d$,
$\nabla G(u)=\E[\nabla_\mu g(\mathbf X;u)]$.

\smallskip\noindent\emph{The case $\ell=1$.}\enspace
Here $g(x;u)=\norm{x-u}$, so $\norm{\nabla_\mu g(x;u)}\leq1$.
Consequently, $R_k=\eta_k=0$ for every $k\geq1$, and
$\left(\nabla\xi_k(\mu_{k-1})\right)_{k\geq1}$ is a martingale
difference sequence with respect to $\left(\F_k\right)_{k\geq0}$.
For every $k\geq1$, $\norm{\nabla\xi_k(\mu_{k-1})}\leq2$ almost surely.
The proof of the almost sure convergence rate of $\left(\mu_n\right)_n$ applies
with Assumption~\ref{ass:taylor} and Lemma~\ref{lem:convolution}.

Let $e_1,\ldots,e_d$ denote the standard basis of $\R^d$.
For $t\ne0$, convexity gives
\[
 0\leq\sum_{j=1}^d
 \frac{\norm{X_1-\mu^*+te_j}+\norm{X_1-\mu^*-te_j}
 -2\norm{X_1-\mu^*}}{t^2}
 \xrightarrow[t\to0]{\hspace{1.5em}}\frac{d-1}{\norm{X_1-\mu^*}}\quad\as,
\]
since $P$ has no atoms. Fatou's lemma and
Assumption~\ref{ass:taylor} yield
\[
 (d-1)\E\left[\frac1{\norm{X_1-\mu^*}}\right]
 \leq\liminf_{t\to0}\sum_{j=1}^d
 \frac{G(\mu^*+te_j)+G(\mu^*-te_j)-2G(\mu^*)}{t^2}
 =\operatorname{tr}(H)<\infty.
\]
Since $d-1>0$, $\E[1/\norm{X_1-\mu^*}]\leq\operatorname{tr}(H)/(d-1)<\infty$.
For nonzero vectors $a,b$, the triangle inequality gives
$\norm{\frac{a}{\norm a}-\frac{b}{\norm b}}
\leq2\min\left(1,\frac{\norm{a-b}}{\norm a}\right)$.
For every $u\in\R^d$, the vectors $\mu^*-X_1$ and $u-X_1$
are nonzero almost surely because $P$ has no atoms.
Applying the inequality with $a\coloneqq\mu^*-X_1$ and $b\coloneqq u-X_1$
and taking expectations gives
\[
 \E\left[\norm{\nabla_\mu g(X_1;u)-\nabla_\mu g(X_1;\mu^*)}^2\right]
 \leq4\min\left(\norm{u-\mu^*}
 \E\left[\frac1{\norm{X_1-\mu^*}}\right],1\right).
\]
For every $k\geq1$, the vector $\mu_{k-1}$ is $\F_{k-1}$-measurable
and $X_k$ is independent of $\F_{k-1}$.
The conditional variance identity therefore gives a constant $C\geq0$
such that
\begin{align*}
 &\E\left[\norm{\nabla\xi_k(\mu_{k-1})-\nabla\xi_k(\mu^*)}^2\mid\F_{k-1}\right]\\
 &\quad=\E\left[\norm{\nabla_\mu g(X_k;\mu_{k-1})
 -\nabla_\mu g(X_k;\mu^*)}^2\mid\F_{k-1}\right]
 -\norm{\nabla G(\mu_{k-1})-\nabla G(\mu^*)}^2\\
 &\quad\leq C\min\left(\norm{\mu_{k-1}-\mu^*}^{\beta},1\right)\quad\as
\end{align*}
This proves \eqref{eq:D-variance} for $\ell=1$.
With $R_k=\eta_k=0$, the proofs of the assertions for the average and
of the asymptotic normality of $\left(\mu_n\right)_n$ apply.

\smallskip\noindent\emph{The case $2\leq\ell<d$.}\enspace
Lemma~\ref{lem:oscillation} applies under
Assumptions~\ref{ass:moment-alternative} and~\ref{ass:inverse-alternative}.
Since $\prod_{i=1}^{\ell}(1+\norm{X_i})\geq1$, its expectation bound
\eqref{eq:osc-expectation-oja} gives a constant $C\geq0$ such that,
for every $u,v\in\R^d$,
\[
 \norm{\nabla G(u)-\nabla G(v)}
 \leq\E\left[\norm{\nabla_\mu g(\mathbf X;u)-\nabla_\mu g(\mathbf X;v)}\right]
 \leq C\norm{u-v}.
\]
Thus $\nabla G$ is Lipschitz. The bound \eqref{eq:envelope} also gives
\[
 \E\left[\norm{\nabla_\mu g(\mathbf X;u)-\nabla_\mu g(\mathbf X;v)}^q\right]
 \leq\left(\frac2\ell\right)^q\E\left[m_{\ell-1}(\Simpl(\mathbf X))^q\right]<\infty.
\]
H\"older's inequality yields
\begin{align*}
 &\E\left[\norm{\nabla_\mu g(\mathbf X;u)-\nabla_\mu g(\mathbf X;v)}^2\right]\\
 &\quad\leq
 \left(\E\left[\norm{\nabla_\mu g(\mathbf X;u)-\nabla_\mu g(\mathbf X;v)}\right]\right)^\beta
 \left(\E\left[\norm{\nabla_\mu g(\mathbf X;u)-\nabla_\mu g(\mathbf X;v)}^q\right]\right)^{1/(q-1)}\\
 &\quad\leq C\norm{u-v}^{\beta}.
\end{align*}
Combining this with the second-moment bound from \eqref{eq:envelope}
gives \eqref{eq:L2-modulus}. Conditioning on $\F_{k-\ell}$ and using
the independence of $Z_k$ and $\F_{k-\ell}$ then gives
\eqref{eq:D-variance} by the conditional variance identity.

For every $k\geq2\ell-1$, the increment bound
\eqref{eq:increment-product}, the pointwise bound
\eqref{eq:osc-majorant} and the expectation bound
\eqref{eq:osc-expectation-oja} in Lemma~\ref{lem:oscillation},
applied conditionally on $\F_{k-\ell}$, give
\[
 \E\left[\norm{\nabla\xi_k(\mu_{k-1})
 -\nabla\xi_k(\mu_{k-\ell})}\mid\F_{k-\ell}\right]
 \leq C\gamma_k\prod_{i=k-2\ell+2}^{k-\ell}(1+\norm{X_i})\quad\as
\]
The product on the right is independent of $\F_{k-2\ell+1}$ and
has expectation $(\E[1+\norm{X_1}])^{\ell-1}$.
Taking conditional expectations given $\F_{k-2\ell+1}$ and using
the definition of $R_k$ yields
\[
 \norm{R_k}\leq
 \E\left[\norm{\nabla\xi_k(\mu_{k-1})
 -\nabla\xi_k(\mu_{k-\ell})}\mid\F_{k-2\ell+1}\right]
 \leq C\gamma_k\quad\as
\]
The bound \eqref{eq:envelope} and the moment of order $q$ give
\[
 \E\left[\norm{\nabla\xi_k(\mu_{k-1})
 -\nabla\xi_k(\mu_{k-\ell})}^q\mid\F_{k-2\ell+1}\right]\leq C\quad\as
\]
The conditional variance identity and conditional H\"older's inequality
therefore give
\[
 \E\left[\norm{\eta_k}^2\mid\F_{k-2\ell+1}\right]
 \leq C\gamma_k^\beta\quad\as
\]
The same envelope and conditional Jensen's inequality give
\[
 \E\left[\norm{\nabla\xi_k(\mu_{k-1})-R_k}^q\mid\F_{k-2\ell+1}\right]
 \leq C\quad\as
\]
We have thus obtained \eqref{eq:overlap-first-moment} and
\eqref{eq:overlap-variance}, together with the conditional moment
bound used in the proof of the almost sure rate.
Since $\gamma q>2$, the three parts of the proof of
Theorem~\ref{th:main} give the conditional assertions of the proposition.
The argument at the end of Appendix~\ref{sec:asymptotics}
gives these assertions under the original law.
\end{proof}

\section{Offline estimation and covariance comparison}\label{sec:batch}
Under Assumptions~\ref{ass:consistency} and~\ref{ass:taylor},
the bound in \eqref{eq:envelope} gives
$\E[\norm{\nabla_\mu g(\mathbf X;\mu^*)}^2]<\infty$.
Since $G$ is differentiable at its minimizer, taking expectations
in the subgradient inequality in \eqref{eq:envelope} gives
$\E[\nabla_\mu g(\mathbf X;\mu^*)]=0$.

By the symmetry of $g$ and the convention in
\eqref{eq:little_g_grad}, the function
$\mathbf x\mapsto\nabla_\mu g(\mathbf x;\mu^*)$
is symmetric in $x_1,\ldots,x_\ell$.
For $1\leq j\leq\ell$, Lemma~\ref{lem:hoeffding}, applied to
$h=\nabla_\mu g(\cdot;\mu^*)$, gives a measurable symmetric
Hoeffding projection $\psi_j:(\R^d)^j\to\R^d$
\citep{Hoeffding1961} satisfying
\[
 \psi_j\left(X_1,\ldots,X_j\right)
 =\sum_{A\subseteq\{1,\ldots,j\}}\left(-1\right)^{j-|A|}
       \E\left[\nabla_\mu g\left(\mathbf X;\mu^*\right)
       \mid X_i,\ i\in A\right]\quad\as
\]
In particular,
\[
 \psi_1(X_1)
 =\E[\nabla_\mu g(\mathbf X;\mu^*)\mid X_1]\quad\as,
 \qquad\text{and}\qquad
 \E[\psi_1(X_1)\psi_1(X_1)^\top]=\zeta_1.
\]

For a finite set $A=\{i_1<\cdots<i_j\}$ of positive integers,
we write $X_A\coloneqq\left(X_{i_1},\ldots,X_{i_j}\right)$.

\subsection{Proof of Theorem~\ref{th:batch}}\label{sec:batch-proof}

\begin{proof}
Under Assumption~\ref{ass:density},
Proposition~\ref{prop:population} gives
Assumption~\ref{ass:consistency}; in the alternative case,
Assumption~\ref{ass:consistency} is imposed directly.
In both cases, $P$ is $(\ell-1)$-smooth and has finite
second moments, by Assumption~\ref{ass:density}
or~\ref{ass:moment-alternative}, respectively.

For a measurable symmetric function $J:(\R^d)^\ell\to\R$,
we define, for $n\geq\ell$, the $U$-statistic
\citep{Hoeffding1948}
\[
 U_n\left(J\right)\coloneqq \binom n\ell^{-1}
 \sum_{1\leq i_1<\cdots<i_\ell\leq n}
 J\left(X_{i_1},\ldots,X_{i_\ell}\right).
\]
In the terminology of $U$-statistics, such a measurable symmetric
function is called a kernel.
Thus the empirical objective function in \eqref{eq:empirical-criterion}
satisfies $\widehat G_n\left(u\right)=U_n\left(g\left(\cdot;u\right)\right)$.

For any integers $1\leq i_1<\cdots<i_\ell$,
Lemma~S.1.3(ii) of~\citet{DPsupp} implies that $\mu^*$
almost surely does not belong to
$\aff(X_{i_1},\ldots,X_{i_\ell})$.
Since there are countably many choices of indices,
with probability one every summand in \eqref{eq:empirical-criterion}
is differentiable at $\mu^*$ for every $n\geq\ell$, and
\[
 \nabla\widehat G_n\left(\mu^*\right)
 =\binom n\ell^{-1}\sum_{1\leq i_1<\cdots<i_\ell\leq n}
     \nabla_\mu g\left(X_{i_1},\ldots,X_{i_\ell};\mu^*\right).
\]
\smallskip\noindent\emph{Existence and measurable selection.}\enspace
Let $n\geq\ell+1$. By Lemma~S.1.3(i) of~\citet{DPsupp},
the points $X_1,\ldots,X_n$ are almost surely not contained
in an affine subspace of dimension $\ell-1$.
Theorems~2.2--2.3 of~\citet{durre_affine-equivariant_2022}
therefore imply that $\widehat G_n$ is coercive and has
a nonempty compact convex set of minimizers.
For every $u\in\R^d$, $\widehat G_n(u)$ is
$\sigma(X_1,\ldots,X_n)$-measurable, and $\widehat G_n$
is continuous.
Theorem~14.37 of~\citet{rockafellar_wets_1998} therefore
gives a $\sigma(X_1,\ldots,X_n)$-measurable random vector
that minimizes $\widehat G_n$ almost surely.
Let $\widehat\mu_n$ be any measurable minimizer.

\smallskip\noindent\emph{Local expansion of the objective function.}\enspace
For $u\in\R^d$, we consider the increment of $\widehat G_n$ from $\mu^*$ to
$\mu^*+u/\sqrt n$.
As in~\citet[equation~(S.2.14)]{DPsupp}, we subtract the linear term
in each summand and multiply the remainder by $n$:
\[
 J_{n,u}\left(\mathbf x\right)\coloneqq n\left(
 g\left(\mathbf x;\mu^*+u/\sqrt n\right)-g\left(\mathbf x;\mu^*\right)
 -\frac{u^\top \nabla_\mu g\left(\mathbf x;\mu^*\right)}{\sqrt n}\right).
\]
Averaging gives
\[
 n\left(\widehat G_n\left(\mu^*+u/\sqrt n\right)
 -\widehat G_n\left(\mu^*\right)\right)
 =u^\top\sqrt n\,\nabla\widehat G_n\left(\mu^*\right)
 +U_n\left(J_{n,u}\right).
\]
We prove that the sequence $\left(U_n(J_{n,u})\right)_n$ converges
in probability to $\tfrac12u^\top Hu$, uniformly for $u$ in bounded sets.

Let $u\in\R^d$. Convexity gives $J_{n,u}\geq0$.
Differentiability of $g(\mathbf X;\cdot)$ at $\mu^*$
and the envelope bound \eqref{eq:envelope} give
\[
 \frac{J_{n,u}(\mathbf X)}{\sqrt n}
 \xrightarrow[n\to\infty]{\hspace{1.5em}}0\quad\as,
 \quad\text{and}\quad
 0\leq\frac{J_{n,u}(\mathbf X)}{\sqrt n}
 \leq\frac{2m_{\ell-1}(\Simpl(\mathbf X))}{\ell}\norm u.
\]
The majorant has a finite second moment by Assumption~\ref{ass:consistency}.
The dominated convergence theorem therefore yields
\[
 \E[J_{n,u}(\mathbf X)^2]=o(n).
\]

Integrating the bound in Assumption~\ref{ass:taylor} along the segment
from $\mu^*$ to $\mu^*+u/\sqrt n$ gives
\[
 \E\left[J_{n,u}\left(\mathbf X\right)\right]
 =n\left(G\left(\mu^*+u/\sqrt n\right)-G\left(\mu^*\right)\right)
 =\tfrac12u^\top Hu+o\left(1\right).
\]

For a measurable symmetric function $J:(\R^d)^\ell\to\R$
with $\E[J(\mathbf X)^2]<\infty$, let $I$ and $B$ range
over the subsets of $\{1,\ldots,n\}$ of size $\ell$.
If $I\cap B=\varnothing$, the random variables $J(X_I)$ and
$J(X_B)$ are independent. If $I\cap B\ne\varnothing$,
the Cauchy--Schwarz inequality gives
$|\operatorname{Cov}(J(X_I),J(X_B))|\leq\E[J(\mathbf X)^2]$.
For each $I$, the number of sets $B$ with $I\cap B\ne\varnothing$ is at most
$\ell\binom{n-1}{\ell-1}$. Therefore,
\begin{equation}\label{eq:U-variance-bound}
 \Var\left(U_n\left(J\right)\right)
 =\binom n\ell^{-2}\sum_{|I|=|B|=\ell}
      \operatorname{Cov}(J(X_I),J(X_B))
 \leq\frac{\ell\binom{n-1}{\ell-1}}{\binom n\ell}
       \E[J(\mathbf X)^2]
 =\frac{\ell^2}{n}\E[J(\mathbf X)^2].
\end{equation}

Applying \eqref{eq:U-variance-bound} with $J=J_{n,u}$ and using
the two moment estimates gives, for every fixed $u\in\R^d$,
\[
 U_n\left(J_{n,u}\right)
 \xrightarrow[n\to\infty]{\hspace{.5em}\Prob\hspace{.5em}}
 \tfrac12u^\top Hu.
\]
The random function $u\mapsto U_n(J_{n,u})$ is finite and convex.
By Theorem~II.1 of~\citet{AndersenGill1982}, applied to the concave
function $u\mapsto-U_n(J_{n,u})$, the pointwise convergence in
probability is uniform on compact subsets of $\R^d$.
Thus, for every $R>0$,
\begin{equation}\label{eq:batch-local-uniform}
 \Delta_n\left(R\right)\coloneqq
 \sup_{\norm u\leq R}\left|U_n\left(J_{n,u}\right)-\tfrac12u^\top Hu\right|
 \xrightarrow[n\to\infty]{\hspace{.5em}\Prob\hspace{.5em}}0.
\end{equation}

\smallskip\noindent\emph{Localization of the minimizer.}\enspace
Since $\E[\nabla_\mu g(\mathbf X;\mu^*)]=0$, applying
\eqref{eq:U-variance-bound} to the scalar functions
$\mathbf x\mapsto\left(\nabla_\mu g(\mathbf x;\mu^*)\right)_r$,
$1\leq r\leq d$, and summing over $r$ gives
\[
 \sup_{n\geq\ell+1}
 \E\left[\norm{\sqrt n\,\nabla\widehat G_n(\mu^*)}^2\right]
 \leq\ell^2\E\left[\norm{\nabla_\mu g(\mathbf X;\mu^*)}^2\right]
 <\infty.
\]
We compare the minimizer of the local objective function with that
of its quadratic approximation. Set
\[
 \widehat u_n\coloneqq\sqrt n\left(\widehat\mu_n-\mu^*\right)
 \quad\text{and}\quad
 v_n\coloneqq-H^{-1}\sqrt n\,\nabla\widehat G_n\left(\mu^*\right).
\]
For a fixed $M>0$, let $R>0$ be sufficiently large that
\[
 \tfrac12\lambda_{\min}(H)R^2-MR-1>0,
 \quad\text{and}\quad R>\norm{H^{-1}}_{\op}M.
\]
On the event on which $\norm{\sqrt n\,\nabla\widehat G_n(\mu^*)}\leq M$
and $\Delta_n(R)\leq1$, the function
$u\mapsto U_n(J_{n,u})+u^\top\sqrt n\,\nabla\widehat G_n(\mu^*)$
vanishes at zero and is strictly positive on the sphere of radius $R$.
Convexity therefore implies that this function is strictly positive
for $\norm u\geq R$.
The minimum is at most zero, hence $\norm{\widehat u_n}<R$.
Moreover, $\norm{v_n}\leq\norm{H^{-1}}_{\op}M<R$.

The difference of the quadratic function
$\tfrac12u^\top Hu+u^\top \sqrt n\,\nabla\widehat G_n\left(\mu^*\right)$ at $\widehat u_n$ and $v_n$ equals
$\tfrac12\left(\widehat u_n-v_n\right)^\top H\left(\widehat u_n-v_n\right)$.
Since $\widehat u_n$ minimizes $u\mapsto U_n\left(J_{n,u}\right)+u^\top\sqrt n\,\nabla\widehat G_n\left(\mu^*\right)$, \eqref{eq:batch-local-uniform} gives
\[
 \tfrac12\lambda_{\min}(H)\norm{\widehat u_n-v_n}^2\leq2\Delta_n\left(R\right).
\]
For every $\varepsilon>0$, the preceding comparison gives
\[
 \Prob[\norm{\widehat u_n-v_n}>\varepsilon]
 \leq\Prob[\norm{\sqrt n\,\nabla\widehat G_n(\mu^*)}>M]
 +\Prob\left[\Delta_n(R)>
       \min\left(1,\frac{\lambda_{\min}(H)\varepsilon^2}{4}\right)\right].
\]
For fixed $M$ and $R$ chosen above, the second probability tends to
zero as $n\to\infty$ by \eqref{eq:batch-local-uniform}.
Markov's inequality and the $L^2$ bound give
\[
 \limsup_{n\to\infty}
 \Prob[\norm{\widehat u_n-v_n}>\varepsilon]
 \leq\frac{\ell^2}{M^2}
      \E[\norm{\nabla_\mu g(\mathbf X;\mu^*)}^2].
\]
Letting $M$ tend to infinity and writing $o_{\Prob}(1)$ for a sequence of
random vectors converging to zero in probability as $n\to\infty$,
we obtain
\begin{equation}\label{eq:batch-Bahadur}
 \sqrt n\left(\widehat\mu_n-\mu^*\right)
 =-H^{-1}\sqrt n\,\nabla\widehat G_n\left(\mu^*\right)
   +o_{\Prob}\left(1\right).
\end{equation}

\smallskip\noindent\emph{Hoeffding decomposition and asymptotic normality.}\enspace
We apply Lemma~\ref{lem:hoeffding} with
$h=\nabla_\mu g(\cdot;\mu^*)$.
For $1\leq j\leq\ell$, each subset of $\{1,\ldots,n\}$ of size $j$
is contained in $\binom{n-j}{\ell-j}$ subsets of size $\ell$.
Averaging \eqref{eq:hoeffding-expansion} over all subsets of
$\{1,\ldots,n\}$ of size $\ell$ therefore gives
\[
 \nabla\widehat G_n\left(\mu^*\right)=\sum_{j=1}^{\ell}\frac{\binom\ell j}{\binom nj}
 \sum_{\substack{A\subseteq\{1,\ldots,n\}\\|A|=j}}
       \psi_j\left(X_A\right).
\]
Subtracting the term with $j=1$ and using the orthogonality in
Lemma~\ref{lem:hoeffding}, we obtain
\[
 \E\left[\norm{\nabla\widehat G_n(\mu^*)-
                \frac\ell n\sum_{i=1}^n\psi_1(X_i)}^2\right]
 =\sum_{j=2}^{\ell}\frac{\binom\ell j^2}{\binom nj}
       \E\left[\norm{\psi_j(X_1,\ldots,X_j)}^2\right]
 =O(n^{-2}).
\]
Consequently,
\[
 \sqrt n\,\nabla\widehat G_n(\mu^*)
 =\frac\ell{\sqrt n}\sum_{i=1}^n\psi_1(X_i)+o_{\Prob}(1).
\]
The random vectors $\psi_1(X_i)$, $i\geq1$, are independent and
identically distributed, with mean zero, covariance $\zeta_1$
and finite second moments.
The central limit theorem of \citet{Duflo}, recalled in
Proposition~\ref{th:duflo}, applies with $m=1$.
Together with Slutsky's theorem, it gives
\[
 \sqrt n\,\nabla\widehat G_n(\mu^*)
 \xrightarrow[n\to\infty]{\hspace{.5em}\mathcal L\hspace{.5em}}
 \mathcal N_d(0,\ell^2\zeta_1).
\]
Equation~\eqref{eq:batch-Bahadur} and Slutsky's theorem then give
\eqref{eq:batch-clt}.
\end{proof}

\subsection{Proof of Lemma~\ref{prop:covariance-comparison}}\label{sec:covariance-proof}

\begin{proof}
For distinct nonempty sets $A,B$ of positive integers,
each of size at most $\ell$, Lemma~\ref{lem:hoeffding} gives
\[
 \Cov\left(\psi_{|A|}(X_A),\psi_{|B|}(X_B)\right)=0.
\]
For $k\in\{0,\ldots,\ell-1\}$, decomposing
$\nabla_\mu g(X_1,\ldots,X_\ell;\mu^*)$ and
$\nabla_\mu g(X_{1+k},\ldots,X_{\ell+k};\mu^*)$ gives
\begin{equation}\label{eq:Gamma0-Hoeffding}
 \Gamma_k
 =\sum_{\varnothing\ne A\subseteq\{k+1,\ldots,\ell\}}
       \Var\left(\psi_{|A|}(X_A)\right)
 =\sum_{j=1}^{\ell-k}\binom{\ell-k}{j}
       \Var\left(\psi_j(X_1,\ldots,X_j)\right).
\end{equation}

To sum \eqref{eq:Gamma0-Hoeffding} over $k$, we use, for
$1\leq j\leq\ell$, the binomial identity
$\sum_{k=1}^{\ell-j}\binom{\ell-k}{j}=\binom\ell{j+1}$,
with the convention $\binom\ell{\ell+1}=0$.
Since each $\Gamma_k$ is symmetric, \eqref{eq:Sigma} gives
\begin{align}
 \Sigma
 &=\Gamma_0+2\sum_{k=1}^{\ell-1}\Gamma_k
 =\sum_{j=1}^{\ell}
   \left(\binom\ell j+
         2\sum_{k=1}^{\ell-j}\binom{\ell-k}{j}\right)
   \Var\left(\psi_j(X_1,\ldots,X_j)\right)\notag\\
 &=\sum_{j=1}^{\ell}
   \left(\binom\ell j+2\binom\ell{j+1}\right)
   \Var\left(\psi_j(X_1,\ldots,X_j)\right).
   \label{eq:Sigma-Hoeffding}
\end{align}
For $1\leq j\leq\ell$, we have
\[
 \binom\ell j+2\binom\ell{j+1}
 =\binom\ell j\,\frac{2\ell-j+1}{j+1}.
\]
Since $\Var(\psi_1(X_1))=\zeta_1$,
\eqref{eq:Gamma0-Hoeffding} and \eqref{eq:Sigma-Hoeffding} yield
\[
 \Sigma-\ell^2\zeta_1
 =\sum_{j=2}^{\ell}\binom\ell j\,
       \frac{2\ell-j+1}{j+1}
       \Var\left(\psi_j(X_1,\ldots,X_j)\right)
 \succeq0
\]
and
\[
 \ell\Gamma_0-\Sigma
 =\sum_{j=2}^{\ell}\binom\ell j\,
       \frac{(\ell+1)(j-1)}{j+1}
       \Var\left(\psi_j(X_1,\ldots,X_j)\right)
 \succeq0.
\]
Thus
\[
 \ell^2\zeta_1\preceq\Sigma\preceq\ell\Gamma_0.
\]

Suppose that $\ell\geq2$. We compare
$\Sigma-\ell^2\zeta_1$ with $\ell\Gamma_0-\ell^2\zeta_1$.
Equation~\eqref{eq:Gamma0-Hoeffding} with $k=0$ gives
\[
 \ell\Gamma_0-\ell^2\zeta_1
 =\ell\sum_{j=2}^{\ell}\binom\ell j\,
       \Var\left(\psi_j(X_1,\ldots,X_j)\right).
\]
For $2\leq j\leq\ell$, we have $3\leq j+1\leq\ell+1$, hence
\[
 1\leq\frac{2(\ell+1)}{j+1}-1
 =\frac{2\ell-j+1}{j+1}
 \leq\frac{2\ell-1}{3}.
\]
The matrices $\Var(\psi_j(X_1,\ldots,X_j))$ are positive semidefinite, so
\[
 \frac1\ell\left(\ell\Gamma_0-\ell^2\zeta_1\right)
 \preceq\Sigma-\ell^2\zeta_1
 \preceq\frac{2\ell-1}{3\ell}
       \left(\ell\Gamma_0-\ell^2\zeta_1\right).
\]
Multiplication on the left and on the right by $H^{-1}$ gives
the bounds stated in the remark following the lemma:
\[
 \frac1\ell\left(V_{\mathrm{ASGD}}-V_{\mathrm{offline}}\right)
 \preceq V_{\mathrm{SW}}-V_{\mathrm{offline}}
 \preceq\frac{2\ell-1}{3\ell}
       \left(V_{\mathrm{ASGD}}-V_{\mathrm{offline}}\right).
\]
For $\ell=2$, both coefficients are $1/2$, and therefore
$V_{\mathrm{SW}}=\frac12(V_{\mathrm{offline}}+V_{\mathrm{ASGD}})$.

\smallskip\noindent\emph{Equality cases.}\enspace
Suppose that $\ell\geq2$. Let $w\in\R^d$ be a fixed unit vector. By \eqref{eq:Gamma0-Hoeffding} and \eqref{eq:Sigma-Hoeffding},
\begin{equation*}
 w^\top\left(\ell\Gamma_0-\Sigma\right)w
 =
 \sum_{j=2}^{\ell}
 \binom\ell j\frac{\left(\ell+1\right)\left(j-1\right)}{j+1}
 \E\left[\left(w^\top\psi_j\left(X_1,\ldots,X_j\right)\right)^2\right].
\end{equation*}
The expansion of $w^\top\left(\Sigma-\ell^2\zeta_1\right)w$ contains the same terms, also with strictly positive coefficients. Thus either quadratic form is zero if and only if $w^\top\psi_j\left(X_1,\ldots,X_j\right)=0$ almost surely for every $2\leq j\leq\ell$.
In that case, we set $a\left(x\right)\coloneqq w^\top\psi_1\left(x\right)$.
The Hoeffding decomposition \eqref{eq:hoeffding-expansion} gives
\begin{equation*}
 w^\top\nabla_\mu g\left(\mathbf X;\mu^*\right)
 =\sum_{i=1}^{\ell}a\left(X_i\right)
 \qquad\as
\end{equation*}

We show that this identity implies $a\left(X_1\right)=0$ almost surely.
If $x_1,\ldots,x_\ell$ belong to a ball of radius $r$, the vectors
$x_i-x_\ell$, $1\leq i\leq\ell-1$, have norm at most $2r$.
The Gram determinant formula for the volume of the simplex
$\Simpl\left(\mathbf x\right)$, Hadamard's inequality and
\eqref{eq:envelope} therefore give
\begin{equation*}
 \abs{w^\top\nabla_\mu g\left(\mathbf x;\mu^*\right)}
 \leq \frac{m_{\ell-1}\left(\Simpl\left(\mathbf x\right)\right)}{\ell}
 \leq \frac1{\ell!}\prod_{i=1}^{\ell-1}\norm{x_i-x_\ell}
 \leq \frac{\left(2r\right)^{\ell-1}}{\ell!}.
\end{equation*}
Let $\varepsilon>0$ and choose $r>0$ such that
$\left(2r\right)^{\ell-1}/\ell!<\ell\varepsilon$.
For every ball $B$ of radius $r$, we have
$\sum_{i=1}^{\ell}a(X_i)
 =w^\top\nabla_\mu g(\mathbf X;\mu^*)
 <\ell\varepsilon$
almost surely on the event that $X_1,\ldots,X_\ell$ belong to $B$.
By independence,
\[
 \Prob\left[X_1\in B,\ a(X_1)>\varepsilon\right]^\ell
 =\Prob\left[X_i\in B,\ a(X_i)>\varepsilon
       \text{ for every }1\leq i\leq\ell\right]
 =0.
\]
Since $\R^d$ is a countable union of balls of radius $r$,
$\Prob[a(X_1)>\varepsilon]=0$.
As $\varepsilon>0$ is arbitrary, ${a(X_1)\leq0}$ almost surely.
Since $\E[a(X_1)]=0$, we have $a(X_1)=0$ almost surely.
Hence
$w^\top\nabla_\mu g\left(\mathbf X;\mu^*\right)=0$ almost surely, or equivalently $w^\top\Gamma_0w=0$.
Conversely, if $w^\top\Gamma_0w=0$, then
$w^\top\nabla_\mu g(\mathbf X;\mu^*)=0$ almost surely.
Hence, for every $A\subseteq\{1,\ldots,\ell\}$, we have
$\E[w^\top\nabla_\mu g(\mathbf X;\mu^*)\mid X_i,\ i\in A]=0$
almost surely.
The definition of $\psi_j$ then gives
$w^\top\psi_j(X_1,\ldots,X_j)=0$ almost surely for every
$1\leq j\leq\ell$.
We have therefore proved
\begin{equation}\label{eq:covariance-kernels}
 \ker\left(\Sigma-\ell^2\zeta_1\right)
 =\ker\left(\ell\Gamma_0-\Sigma\right)
 =\ker\left(\Gamma_0\right).
\end{equation}
For $v\in\R^d$ such that $v^\top V_{\mathrm{ASGD}}v>0$,
we have $H^{-1}v\notin\ker(\Gamma_0)$.
Equation~\eqref{eq:covariance-kernels} therefore gives
\[
 v^\top V_{\mathrm{offline}}v
 <
 v^\top V_{\mathrm{SW}}v
 <
 v^\top V_{\mathrm{ASGD}}v.
\]

\smallskip\noindent\emph{Distinct covariance matrices.}\enspace
We prove that $\Gamma_0\ne0$ by contradiction.
Suppose that $\Gamma_0=0$.
By the definition of $\Gamma_0$,
$\nabla_\mu g(\mathbf X;\mu^*)=0$ almost surely.
Equation~\eqref{eq:little_g_grad} gives
$g(\mathbf X;\mu^*)=0$ almost surely, and hence $G(\mu^*)=0$
in this case.
Let $u\in\R^d\setminus\{0\}$.
The simplex-volume formula~\eqref{eq:base-height} gives,
for every $t\in\R$,
\[
 G(\mu^*+tu)=|t|G(\mu^*+u).
\]
Since $G(\mu^*)=0$, dividing by $t>0$ and letting $t$ tend to
zero gives
$G(\mu^*+u)=u^\top\nabla G(\mu^*)=0$,
contrary to the uniqueness of the minimizer in
Assumption~\ref{ass:consistency}.
Thus $\Gamma_0\ne0$.
By \eqref{eq:covariance-kernels}, the two positive semidefinite
matrices $\Sigma-\ell^2\zeta_1$ and $\ell\Gamma_0-\Sigma$ are
nonzero.
The matrices $V_{\mathrm{offline}}$, $V_{\mathrm{SW}}$ and
$V_{\mathrm{ASGD}}$ are therefore pairwise distinct.

\smallskip\noindent
\emph{Positive definiteness under the density assumption.}\enspace
Suppose that Assumption~\ref{ass:density} holds.
We prove that $\Gamma_0$ is positive definite.
The Gram determinant formula for simplex volumes gives,
for every $\mathbf x$ and $\mu$,
\begin{equation}\label{eq:simplex-gram-volume}
 g(\mathbf x;\mu)^2
 =\frac1{(\ell!)^2}
 \det\!\left(
   \left((x_i-\mu)^\top(x_j-\mu)\right)_{i,j=1}^{\ell}
 \right).
\end{equation}
Each entry of the Gram matrix is polynomial in the coordinates
of $\mathbf x$ and $\mu$, hence so is $g(\mathbf x;\mu)^2$.
For a unit vector $w$, the directional derivative
\[
 Q_w(\mathbf x)\coloneqq
 \left.\frac{d}{dt}g(\mathbf x;\mu^*+tw)^2\right|_{t=0}
\]
is therefore a polynomial in the coordinates of $\mathbf x$.
Differentiating the square of the simplex-volume formula
\eqref{eq:base-height} and using \eqref{eq:little_g_grad} gives
\[
 Q_w(\mathbf x)
 =2\left(\frac{m_{\ell-1}(\Simpl(\mathbf x))}{\ell}\right)^2
   w^\top\Gamma(\mathbf x)(\mu^*-x_\ell)
 =2g(\mathbf x;\mu^*)\,w^\top\nabla_\mu g(\mathbf x;\mu^*).
\]
We choose $\ell$ affinely independent points $x_1,\ldots,x_\ell$
in the hyperplane
$\{x\in\R^d:w^\top(x-\mu^*)=1\}$.
For this tuple, $\Gamma(\mathbf x)w=w$ and
$w^\top\Gamma(\mathbf x)(\mu^*-x_\ell)=-1$, so
\[
 Q_w(\mathbf x)
 =-2\left(\frac{m_{\ell-1}(\Simpl(\mathbf x))}{\ell}\right)^2<0.
\]
Thus $Q_w$ is a nonzero polynomial, so its zero set has
Lebesgue measure zero.
Since the joint law of $\mathbf X$ has a Lebesgue density,
$Q_w(\mathbf X)\ne0$ almost surely, and therefore
\[
 w^\top\Gamma_0w
 =\E\!\left[
   \left(w^\top\nabla_\mu g(\mathbf X;\mu^*)\right)^2
 \right]>0.
\]
Hence $\Gamma_0$ is positive definite.
Equation~\eqref{eq:covariance-kernels} gives
$\Sigma-\ell^2\zeta_1\succ0$ and $\ell\Gamma_0-\Sigma\succ0$.
\end{proof}

\subsection{Repeated updates on disjoint blocks}\label{sec:repeated-blocks}
For a fixed integer $r\geq1$, we perform $r$ successive gradient
updates as soon as a new block of $\ell$ consecutive random vectors
is available. For $m\geq1$, write
$B_m\coloneqq(X_{(m-1)\ell+1},\ldots,X_{m\ell})$.
Starting from $\mu_0^{(r)}\coloneqq\mu_0$, define, for $j\geq1$,
\[
 \mu_j^{(r)}
 \coloneqq\mu_{j-1}^{(r)}
 -\gamma_{\lceil j/r\rceil}
 \nabla_\mu g\left(B_{\lceil j/r\rceil};\mu_{j-1}^{(r)}\right),
\]
where $\gamma_m\coloneqq c_\gamma m^{-\gamma}$, with $c_\gamma$
and the fixed exponent $\gamma$ as in Section~\ref{sec::algo}.
The gradient is reevaluated after each update.
After $n\geq\ell$ observations, this recursion has performed
$r\lfloor n/\ell\rfloor$ updates, and its averaged estimator is
\[
 \bar\mu_n^{(r)}
 \coloneqq\frac1{r\lfloor n/\ell\rfloor}
 \sum_{j=1}^{r\lfloor n/\ell\rfloor}\mu_j^{(r)}.
\]

\begin{proposition}\label{prop:repeated-blocks}
Under Assumption~\ref{ass:consistency}, the sequence
$\left(\mu_j^{(r)}\right)_j$ converges almost surely to $\mu^*$.
If Assumptions~\ref{ass:density} and~\ref{ass:taylor} hold, with
$\gamma\nu>2$ and $\gamma\left(1+\alpha\right)>1$, then
\[
 \sqrt n\left(\bar\mu_n^{(r)}-\mu^*\right)
 \xrightarrow[n\to\infty]{\hspace{.5em}\mathcal L\hspace{.5em}}
 \mathcal N_d\left(0,\ell H^{-1}\Gamma_0H^{-1}\right).
\]
The same conclusions hold if the step size
$\gamma_{\lceil j/r\rceil}$ is replaced by $\gamma_j$.
\end{proposition}

\begin{proof}
Throughout the proof, we omit the superscript $(r)$.
We work conditionally on $\F_0$, as explained at the beginning of
Appendix~\ref{sec:proofs}.
For asymptotic normality, we use the averaging method of
\citet[Theorem~2]{PolyakJud92}, as in the proof of the assertions
for $\bar\mu_n$ in Theorem~\ref{th:main}.
Here we also control the error caused by the repeated use of each block.
We first consider the step size $\gamma_{\lceil j/r\rceil}$ at update $j$.
For every block index $m\geq1$ and every $1\leq s\leq r$, we have
\[
 \mu_{r(m-1)+s}
 =\mu_{r(m-1)+s-1}
 -\gamma_m\nabla_\mu g(B_m;\mu_{r(m-1)+s-1}).
\]
For every $m\geq1$ and $0\leq s\leq r$, the envelope
\eqref{eq:envelope} gives
\begin{equation}\label{eq:repeated-increments}
 \norm{\mu_{r(m-1)+s}-\mu_{r(m-1)}}
 \leq\frac{s\gamma_m}{\ell}m_{\ell-1}(\Simpl(B_m)).
\end{equation}

\smallskip\noindent\emph{Almost sure convergence.}\enspace
Summing the subgradient inequalities over the $r$ updates and
using the Lipschitz bound in \eqref{eq:envelope} gives
\[
\begin{aligned}
 \norm{\mu_{rm}-\mu^*}^2
 \leq{}&\norm{\mu_{r(m-1)}-\mu^*}^2
 -2r\gamma_m\bigl(g(B_m;\mu_{r(m-1)})-g(B_m;\mu^*)\bigr)\\
 &+\frac{r^2\gamma_m^2}{\ell^2}m_{\ell-1}(\Simpl(B_m))^2.
\end{aligned}
\]
By induction, each iterate is square-integrable under the conditional law.
Since $B_m$ is independent of $\F_{(m-1)\ell}$, conditional
expectation yields
\begin{equation}\label{eq:repeated-descent}
\begin{aligned}
 \E[\norm{\mu_{rm}-\mu^*}^2\mid\F_{(m-1)\ell}]
 \leq{}&\norm{\mu_{r(m-1)}-\mu^*}^2
 -2r\gamma_m\bigl(G(\mu_{r(m-1)})-G(\mu^*)\bigr)\\
 &+\frac{r^2\gamma_m^2}{\ell^2}
 \E[m_{\ell-1}(\Simpl(\mathbf X))^2].
\end{aligned}
\end{equation}
Assumption~\ref{ass:consistency} gives
$\E[m_{\ell-1}(\Simpl(\mathbf X))^2]<\infty$.
Since $\sum_m\gamma_m=\infty$ and $\sum_m\gamma_m^2<\infty$,
the Robbins--Siegmund theorem
\citep[Theorem~1, p.~233]{RobbinsSiegmund1971} shows that
$\left(\norm{\mu_{rm}-\mu^*}^2\right)_m$ converges almost surely and
$\sum_m\gamma_m(G(\mu_{r(m-1)})-G(\mu^*))<\infty$ almost surely.
Consequently, the sequence $\left(\mu_{rm}\right)_m$ is bounded and
$\liminf_mG(\mu_{rm})=G(\mu^*)$ almost surely.
By continuity of $G$ and uniqueness of its minimizer,
a subsequence converges to $\mu^*$.
The convergence of the squared distances therefore implies that
$\left(\mu_{rm}\right)_m$ converges almost surely to $\mu^*$.
The moment condition in Assumption~\ref{ass:consistency} and
$\sum_m\gamma_m^2<\infty$ give, by conditional expectation given $\F_0$,
$\sum_m\gamma_m^2m_{\ell-1}(\Simpl(B_m))^2<\infty$ almost surely.
Equation~\eqref{eq:repeated-increments} therefore implies that
$\left(\mu_j\right)_j$ converges almost surely to $\mu^*$.

\smallskip\noindent\emph{Error from repeated updates and convergence rate.}\enspace
Suppose now that Assumptions~\ref{ass:density} and~\ref{ass:taylor}
hold, with $\gamma\nu>2$ and $\gamma(1+\alpha)>1$.
We compare the sum of the $r$ updates with a single update of size
$r\gamma_m$ taken at $\mu_{r(m-1)}$. Define
\[
 \varepsilon_m\coloneqq\frac1r\sum_{s=1}^r
 \bigl(\nabla_\mu g(B_m;\mu_{r(m-1)+s-1})
       -\nabla_\mu g(B_m;\mu_{r(m-1)})\bigr).
\]
Then
\[
 \mu_{rm}=\mu_{r(m-1)}
 -r\gamma_m\bigl(\nabla_\mu g(B_m;\mu_{r(m-1)})+\varepsilon_m\bigr).
\]
Equations~\eqref{eq:repeated-increments} and~\eqref{eq:product-bound} give
\[
 \max_{0\leq s\leq r}\norm{\mu_{r(m-1)+s}-\mu_{r(m-1)}}
 \leq C\gamma_m\prod_{i=(m-1)\ell+1}^{m\ell}(1+\norm{X_i}).
\]
Conditionally on $\F_{(m-1)\ell}$, the initial point
$\mu_{r(m-1)}$ is fixed and $B_m$ has law $P^{\otimes\ell}$.
The pointwise bound \eqref{eq:osc-majorant} and the expectation
bounds in Lemma~\ref{lem:oscillation},
with $t=C\gamma_m$, therefore give
\begin{equation}\label{eq:repeated-error-bound}
 \E[\norm{\varepsilon_m}\mid\F_{(m-1)\ell}]
 \leq C\gamma_m.
\end{equation}
Let $q\in(2/\gamma,\nu)$. We define
\[
 Y_m\coloneqq
 \nabla_\mu g(B_m;\mu_{r(m-1)})
 -\nabla G(\mu_{r(m-1)})+\varepsilon_m
 -\E[\varepsilon_m\mid\F_{(m-1)\ell}].
\]
The sequence $\left(Y_m\right)_m$ is a martingale difference
sequence for the filtration $\left(\F_{m\ell}\right)_m$.
By the definition of $\varepsilon_m$ and the envelope
\eqref{eq:envelope},
\[
 \norm{\nabla_\mu g(B_m;\mu_{r(m-1)})+\varepsilon_m}
 =\norm{\frac1r\sum_{s=1}^r
 \nabla_\mu g(B_m;\mu_{r(m-1)+s-1})}
 \leq\frac{m_{\ell-1}(\Simpl(B_m))}{\ell}.
\]
The vector $Y_m$ equals
$\nabla_\mu g(B_m;\mu_{r(m-1)})+\varepsilon_m$ minus its conditional
expectation given $\F_{(m-1)\ell}$.
The triangle inequality, conditional Jensen's inequality and
\eqref{eq:moment-L} therefore give, for every $m\geq1$,
\[
 \E[\norm{Y_m}^q\mid\F_{(m-1)\ell}]
 \leq\frac{2^q}{\ell^q}
 \E[m_{\ell-1}(\Simpl(\mathbf X))^q]<\infty\quad\as
\]
For $k\geq0$, we write
$\delta_k\coloneqq\nabla G(\mu_k)-H(\mu_k-\mu^*)$
for the Taylor remainder.
The linearized recursion is
\[
 \mu_{rm}-\mu^*
 =(\Id-r\gamma_mH)(\mu_{r(m-1)}-\mu^*)
 -r\gamma_m\bigl(
 Y_m+\E[\varepsilon_m\mid\F_{(m-1)\ell}]+\delta_{r(m-1)}\bigr).
\]
Let $m_0$ be the smallest positive integer such that
$r\gamma_{m_0}\norm H_{\op}\leq1/2$.
This integer is $\F_0$-measurable. Since $\left(\gamma_m\right)_m$
is decreasing, $r\gamma_m\norm H_{\op}\leq1/2$ for every $m\geq m_0$.
For $m\geq k\geq m_0$, set
\[
 \Pi_{m,k}\coloneqq\prod_{j=k+1}^m(\Id-r\gamma_jH).
\]
For every $m\geq m_0$, iterating the recursion from $m_0$ gives
\[
 \mu_{rm}-\mu^*
 =\Pi_{m,m_0}(\mu_{rm_0}-\mu^*)
 -\sum_{k=m_0+1}^m r\gamma_k\Pi_{m,k}
 \bigl(Y_k+\E[\varepsilon_k\mid\F_{(k-1)\ell}]+\delta_{r(k-1)}\bigr).
\]
Lemma~\ref{lem:convolution}, with delay one and step sizes $r\gamma_k$,
and Lemma~\ref{lem:deterministic-weights}, with
$\norm{\E[\varepsilon_k\mid\F_{(k-1)\ell}]}\leq C\gamma_k$, give,
for every $m>m_0$,
\[
 \norm{\mu_{rm}-\mu^*}
 \leq C_0\sqrt{\gamma_m\log m}
 +\sum_{k=m_0+1}^m r\gamma_k\norm{\Pi_{m,k}}_{\op}\norm{\delta_{r(k-1)}}
 \quad\as
\]
for a finite random variable $C_0\geq0$.
Assumption~\ref{ass:taylor} and almost sure convergence give
$\norm{\delta_{rm}}=o(\norm{\mu_{rm}-\mu^*})$ almost surely.
Lemma~\ref{lem:perturbed-weights}, applied with step sizes
$r\gamma_m$, $n_0=m_0$ and
$v_m=(C_0+1)\sqrt{\gamma_m\log m}$, therefore yields
\[
 \norm{\mu_{rm}-\mu^*}^2=O(m^{-\gamma}\log m)\quad\as
\]

\smallskip\noindent\emph{Asymptotic normality.}\enspace
We set $\xi_m^*\coloneqq\nabla_\mu g(B_m;\mu^*)$ and
\[
 D_m\coloneqq\nabla_\mu g(B_m;\mu_{r(m-1)})
              -\nabla G(\mu_{r(m-1)})-\xi_m^*.
\]
The sequence $\left(D_m\right)_m$ is a martingale difference
sequence for $\left(\F_{m\ell}\right)_m$.
By \eqref{eq:L2-modulus}, with $\beta=(q-2)/(q-1)$, almost sure
convergence of $\left(\mu_{rm}\right)_m$ and conditional dominated
convergence give $\E[\norm{D_m}^2\mid\F_0]\to0$ almost surely.
For every $m\geq1$, orthogonality under the conditional laws gives
\[
 \E\!\left[\norm{\frac1{\sqrt m}\sum_{k=1}^m D_k}^2\mid\F_0\right]
 =\frac1m\sum_{k=1}^m\E[\norm{D_k}^2\mid\F_0]
 \xrightarrow[m\to\infty]{\hspace{1.5em}}0\quad\as
\]
The bound \eqref{eq:repeated-error-bound} and the tower property give
\[
 \E\!\left[\norm{\frac1{\sqrt m}\sum_{k=1}^m\varepsilon_k}\mid\F_0\right]
 \leq \frac C{\sqrt m}\sum_{k=1}^m\gamma_k
 =O(m^{1/2-\gamma})
 \xrightarrow[m\to\infty]{\hspace{1.5em}}0\quad\as
\]
The convergence rate, Assumption~\ref{ass:taylor} and
$\gamma(1+\alpha)>1$ imply
\[
 \frac1{\sqrt m}\sum_{k=1}^m\norm{\delta_{r(k-1)}}
 =O\!\left(m^{(1-\gamma(1+\alpha))/2}
               (\log m)^{(1+\alpha)/2}\right)
 \xrightarrow[m\to\infty]{\hspace{1.5em}}0\quad\as
\]
Summation by parts gives
\[
 \sum_{k=1}^m\frac{\mu_{r(k-1)}-\mu_{rk}}{r\gamma_k}
 =\frac{\mu_0-\mu^*}{r\gamma_1}
 -\frac{\mu_{rm}-\mu^*}{r\gamma_m}
 +\sum_{k=1}^{m-1}(\mu_{rk}-\mu^*)
       \left(\frac1{r\gamma_{k+1}}-\frac1{r\gamma_k}\right).
\]
Since $\gamma_{k+1}^{-1}-\gamma_k^{-1}=O(k^{\gamma-1})$,
\[
 \norm{\frac1{\sqrt m}\sum_{k=1}^m
        \frac{\mu_{r(k-1)}-\mu_{rk}}{r\gamma_k}}
 =O\!\left(m^{(\gamma-1)/2}\sqrt{\log m}\right)
 \xrightarrow[m\to\infty]{\hspace{1.5em}}0\quad\as
\]
Summing the recursion after division by $r\gamma_k$ gives
\[
 \frac H{\sqrt m}\sum_{k=1}^m(\mu_{r(k-1)}-\mu^*)
 =\frac1{\sqrt m}\sum_{k=1}^m
        \frac{\mu_{r(k-1)}-\mu_{rk}}{r\gamma_k}
 -\frac1{\sqrt m}\sum_{k=1}^m\xi_k^*
 -\frac1{\sqrt m}\sum_{k=1}^m(D_k+\varepsilon_k+\delta_{r(k-1)}).
\]
Moreover, \eqref{eq:repeated-increments} gives
\[
 \sqrt m\,\E\!\left[\norm{\bar\mu_{\ell m}
           -\frac1m\sum_{k=1}^m\mu_{r(k-1)}}\mid\F_0\right]
 \leq\frac{C\E[m_{\ell-1}(\Simpl(\mathbf X))]}{\ell\sqrt m}\sum_{k=1}^m\gamma_k
 \xrightarrow[m\to\infty]{\hspace{1.5em}}0\quad\as
\]
Consequently,
\[
 \sqrt m(\bar\mu_{\ell m}-\mu^*)
 =-\frac{H^{-1}}{\sqrt m}\sum_{k=1}^m\xi_k^*+o_{\Prob(\cdot\mid\F_0)}(1).
\]
The random vectors $\xi_k^*$ are independent and identically
distributed, with mean zero and covariance $\Gamma_0$, and the
sequence $\left(\xi_k^*\right)_k$ is independent of $\F_0$.
The central limit theorem of \citet{Duflo}, recalled in
Proposition~\ref{th:duflo}, applies with dependence parameter one.
Together with Slutsky's theorem, it gives weak convergence:
\[
 \mathcal L\left(\sqrt m(\bar\mu_{\ell m}-\mu^*)\mid\F_0\right)
 \xrightarrow[m\to\infty]{\hspace{1.5em}}
 \mathcal N_d(0,H^{-1}\Gamma_0H^{-1})\quad\as
\]
Since the limiting law is deterministic, the argument with bounded
continuous functions at the end of Appendix~\ref{sec:asymptotics}
gives the same convergence under the original law.

For $n\geq\ell$, set $m=\lfloor n/\ell\rfloor$.
By definition, $\bar\mu_n=\bar\mu_{\ell m}$.
Since $n/m\to\ell$, Slutsky's theorem gives the stated convergence.

\smallskip\noindent\emph{Step sizes indexed by updates.}\enspace
We now use the step size $\gamma_j$ at update $j$, again denoting
the iterates by $\mu_j$.
For every $m\geq1$ and $1\leq s\leq r$,
\[
 0\leq\frac{\gamma_{r(m-1)+s}}{\gamma_{rm}}-1
 =\left(1-\frac{r-s}{rm}\right)^{-\gamma}-1
 \leq\frac Cm.
\]
The envelope \eqref{eq:envelope} gives \eqref{eq:repeated-increments}
with $\gamma_m$ replaced by $C\gamma_{rm}$.
Since $\sum_{s=1}^r\gamma_{r(m-1)+s}\geq r\gamma_{rm}$ and
$G(\mu_{r(m-1)})\geq G(\mu^*)$, the subgradient calculation also
gives \eqref{eq:repeated-descent} with $\gamma_m$ replaced by
$\gamma_{rm}$ and $r^2$ by $C$.
The almost sure convergence argument therefore applies under
Assumption~\ref{ass:consistency}.

The block recursion holds with $\gamma_{rm}$ in place of
$\gamma_m$ if we set
\[
 \varepsilon_m\coloneqq\frac1r\sum_{s=1}^r
 \left(
 \frac{\gamma_{r(m-1)+s}}{\gamma_{rm}}
 \nabla_\mu g(B_m;\mu_{r(m-1)+s-1})
 -\nabla_\mu g(B_m;\mu_{r(m-1)})
 \right).
\]
Replacing the coefficients $\gamma_{r(m-1)+s}/\gamma_{rm}$ by one
changes this remainder by a vector of norm at most
$C m_{\ell-1}(\Simpl(B_m))/(m\ell)$.
Under the assumptions for asymptotic normality,
Lemma~\ref{lem:oscillation} and the envelope therefore give
\[
 \E[\norm{\varepsilon_m}\mid\F_{(m-1)\ell}]
 \leq C(\gamma_{rm}+m^{-1})=O(\gamma_{rm}),
 \quad\text{and}\quad\norm{\varepsilon_m}\leq\frac C\ell m_{\ell-1}(\Simpl(B_m)).
\]
The envelope and \eqref{eq:moment-L} give uniformly bounded
conditional moments of order $q$ for $\left(Y_m\right)_m$.
Since $\gamma_{rm}=r^{-\gamma}\gamma_m$ is a fixed positive
multiple of $\gamma_m$, the rate and averaging
arguments apply with $\gamma_{rm}$ in place of $\gamma_m$.
\end{proof}

\section{Geometric estimates}\label{sec:preliminary}\label{sec:auxiliary}\label{sec:geometry}

By Lemma S.1.2 of~\citet{DPsupp}, there exists a constant $C>0$
such that, for every $\mathbf x\in(\R^d)^\ell$,
\begin{equation}\label{eq:product-bound}
 \frac1\ell m_{\ell-1}\left(\Simpl\left(\mathbf x\right)\right)+g\left(\mathbf x;0\right)
 \leq C\prod_{i=1}^{\ell}\left(1+\norm{x_i}\right).
\end{equation}
Under Assumption~\ref{ass:density}, independence gives, for $0<q<\nu$,
\begin{equation}\label{eq:moment-L}
 \E\left[m_{\ell-1}\left(\Simpl\left(\mathbf X\right)\right)^q\right]<\infty,
 \quad\text{and}\quad
 \E\left[\prod_{i=1}^{\ell}\left(1+\norm{X_i}\right)^q\right]
 =\left(\E\left[\left(1+\norm{X_1}\right)^q\right]\right)^{\ell}<\infty.
\end{equation}

Under Assumption~\ref{ass:density}, $P$ assigns zero mass to every
affine subspace of dimension at most $\ell-1$.
Lemma S.1.3 of~\citet{DPsupp} therefore gives, for every
$u\in\R^d$,
\[
 \dist(u,\aff(\mathbf X))>0
 \quad\as
\]

We next bound weighted probabilities that
$\dist(u,\aff(\mathbf X))$ is small, uniformly in $u$.

Let $\Gr(\ell-1,d)$ be the set of all $(\ell-1)$-dimensional
linear subspaces of $\R^d$, called the Grassmannian.
We denote by $\vartheta$ its unique rotation-invariant
Borel probability measure
\citep[Theorem~13.2.11]{SchneiderWeil2008}.
Every affine subspace of dimension $\ell-1$ has a unique
representation $y+V$, with $V\in\Gr(\ell-1,d)$ and $y\in V^\perp$.
By the affine Blaschke--Petkantschin formula
\citep[eq.~(27), p.~361]{Miles1971}, there exists a constant $c>0$
such that, for every nonnegative measurable function $F$ on
$(\R^d)^\ell$,
\begin{equation}\label{eq:BP}
\begin{aligned}
 &\int_{(\R^d)^\ell}F(\mathbf x)\,d\mathbf x\\
 &\quad=c\int_{\Gr(\ell-1,d)}\int_{V^\perp}\int_{V^\ell}
 F(y+z_1,\ldots,y+z_\ell)\,
 m_{\ell-1}(\Simpl(\mathbf z))^{d-\ell+1}
 \,d\mathbf z\,dy\,d\vartheta(V),
\end{aligned}
\end{equation}
where $\mathbf z\coloneqq(z_1,\ldots,z_\ell)$ and
$d\mathbf z=dz_1\cdots dz_\ell$.
The measures $dy$ and $dz_i$ are the Euclidean Lebesgue measures
on $V^\perp$ and $V$, of dimensions $d-\ell+1$ and $\ell-1$,
respectively.

For any linear subspace $V$ of $\R^d$, we denote by $\Gamma_V$
the orthogonal projection onto $V$.
For $V\in\Gr(\ell-1,d)$, $y\in V^\perp$ and
$\mathbf z\in V^\ell$, we set
$\mathbf x\coloneqq(y+z_1,\ldots,y+z_\ell)$.
For almost every $\mathbf z\in V^\ell$, we have
$\aff(\mathbf z)=V$ and therefore $\aff(\mathbf x)=y+V$.
Hence, for every $u\in\R^d$,
\[
 \Gamma(\mathbf x)=\Gamma_{V^\perp},
 \quad\text{and}\quad
 \dist(u,\aff(\mathbf x))=\norm{\Gamma_{V^\perp}u-y}.
\]

\begin{lemma}\label{lem:tube}
Under Assumption~\ref{ass:density}, for every $0\leq a<\nu-1$
there exists a finite constant $C_a>0$ such that, for every $t>0$,
\[
 \sup_{u\in\R^d}\E\left[
 \frac{m_{\ell-1}\left(\Simpl\left(\mathbf X\right)\right)}{\ell}
 \prod_{i=1}^{\ell}\left(1+\norm{X_i}\right)^a
 \ind_{\{\dist\left(u,\aff\left(\mathbf X\right)\right)\leq t\}}\right]
 \leq C_a t^{d-\ell+1}.
\]
\end{lemma}
\begin{proof}
For $V\in\Gr(\ell-1,d)$ and $y\in V^\perp$, we set
\begin{equation}\label{eq:a-density}
 b_a\left(V,y\right)\coloneqq \frac{c}{\ell}\int_{V^\ell}
 m_{\ell-1}\left(\Simpl\left(\mathbf z\right)\right)^{d-\ell+2}
 \prod_{i=1}^{\ell}\left(\left(1+\norm{y+z_i}\right)^a f\left(y+z_i\right)\right)\,d\mathbf z.
\end{equation}
Since $y\perp V$, we have $\norm{y+z_i}\geq\norm{z_i}$ for
every $z_i\in V$. Assumption~\ref{ass:density} therefore gives
\[
 (1+\norm{y+z_i})^a f(y+z_i)
 \leq C_f(1+\norm{y+z_i})^{a-d-\nu}
 \leq C_f(1+\norm{z_i})^{a-d-\nu}.
\]
We identify $V$ isometrically with $\R^{\ell-1}$
using an orthonormal basis.
Combining this bound with \eqref{eq:product-bound}, applied to
$\mathbf z$, we obtain a constant $C>0$ such that, for every
$V\in\Gr(\ell-1,d)$ and $y\in V^\perp$,
\begin{align*}
 b_a(V,y)
 &\leq C\int_{V^\ell}\prod_{i=1}^{\ell}
 (1+\norm{z_i})^{-(\ell+\nu-2-a)}\,d\mathbf z\\
 &=C\left(\int_{\R^{\ell-1}}
 (1+\norm w)^{-(\ell+\nu-2-a)}\,dw\right)^\ell.
\end{align*}
Since $a<\nu-1$, we have $\ell+\nu-2-a>\ell-1$, so the integral
over $\R^{\ell-1}$ is finite.
The resulting bound is independent of $V$ and $y$, hence
\[
 \sup_{V\in\Gr(\ell-1,d)}\sup_{y\in V^\perp}b_a(V,y)<\infty.
\]
For every $u\in\R^d$ and $t>0$, translation invariance of the
simplex volume and the affine integration formula \eqref{eq:BP} give
\[
\begin{aligned}
 &\E\left[\frac{m_{\ell-1}(\Simpl(\mathbf X))}{\ell}
    \prod_{i=1}^{\ell}(1+\norm{X_i})^a
    \ind_{\{\dist(u,\aff(\mathbf X))\leq t\}}\right]\\
 &\quad=\int_{\Gr(\ell-1,d)}\int_{V^\perp}
    \ind_{\{\norm{\Gamma_{V^\perp}u-y}\leq t\}}b_a(V,y)\,dy\,d\vartheta(V).
\end{aligned}
\]
The set $\{y\in V^\perp:\norm{\Gamma_{V^\perp}u-y}\leq t\}$
is a ball of radius $t$ in $V^\perp$, of volume proportional
to $t^{d-\ell+1}$.
Integrating the uniform bound over this ball, then with respect
to the probability measure $\vartheta$, gives
\[
 \int_{\Gr(\ell-1,d)}\int_{V^\perp}
 \ind_{\{\norm{\Gamma_{V^\perp}u-y}\leq t\}}b_a(V,y)\,dy\,d\vartheta(V)
 \leq Ct^{d-\ell+1}\int_{\Gr(\ell-1,d)}d\vartheta(V)
 =Ct^{d-\ell+1}.
\]
This proves the lemma.
\end{proof}

\begin{lemma}\label{lem:oscillation}
Let $\mathbf x$ be an affinely independent tuple.
For every $u\in\R^d\setminus\aff(\mathbf x)$, $v\in\R^d$ and $t>0$
such that $\norm{v-u}\leq t\prod_{i=1}^{\ell}(1+\norm{x_i})$,
we have
\begin{equation}\label{eq:osc-majorant}
 \norm{\nabla_\mu g\left(\mathbf x;v\right)-\nabla_\mu g\left(\mathbf x;u\right)}
 \leq
 \begin{cases}
 \displaystyle\frac{2m_{\ell-1}\left(\Simpl\left(\mathbf x\right)\right)}{\ell}
 \min\!\left(1,\frac{t\prod_{i=1}^{\ell}\left(1+\norm{x_i}\right)}
 {\dist\left(u,\aff\left(\mathbf x\right)\right)}\right),&\text{if }\ell<d,\\[1em]
 \displaystyle\frac{2m_{d-1}\left(\Simpl\left(\mathbf x\right)\right)}{d}
 \ind_{\{\dist\left(u,\aff\left(\mathbf x\right)\right)\leq t\prod_{i=1}^{d}\left(1+\norm{x_i}\right)\}},&\text{if }\ell=d.
 \end{cases}
\end{equation}
Under Assumption~\ref{ass:density}, or under
Assumptions~\ref{ass:moment-alternative}
and~\ref{ass:inverse-alternative} when $2\leq\ell<d$,
there exists a constant $C\geq0$ such that, for every $t>0$,
\begin{equation}\label{eq:osc-expectation-oja}
\begin{aligned}
 \sup_{u\in\R^d}\E\left[\frac{m_{\ell-1}(\Simpl(\mathbf X))}{\ell}
 \min\!\left(1,\frac{t\prod_{i=1}^{\ell}(1+\norm{X_i})}
 {\dist(u,\aff(\mathbf X))}\right)\right]&\leq Ct
 &&\text{if }\ell<d,\\
 \sup_{u\in\R^d}\E\left[\frac{m_{d-1}(\Simpl(\mathbf X))}{d}
 \ind_{\{\dist(u,\aff(\mathbf X))\leq
 t\prod_{i=1}^{d}(1+\norm{X_i})\}}\right]&\leq Ct
 &&\text{if }\ell=d.
\end{aligned}
\end{equation}
For $\ell<d$, we also have
\begin{equation}\label{eq:weighted-inverse}
 \sup_{u\in\R^d}\E\left[
 \frac{m_{\ell-1}\left(\Simpl\left(\mathbf X\right)\right)\prod_{i=1}^{\ell}\left(1+\norm{X_i}\right)}
 {\ell\dist\left(u,\aff\left(\mathbf X\right)\right)}\right]<\infty.
\end{equation}
\end{lemma}
\begin{proof}
We first consider the case $\ell<d$.
For nonzero vectors $a,b$, adding and subtracting $b/\norm a$,
using $|\norm a-\norm b|\leq\norm{a-b}$ and the bound by two gives
$\norm{\frac{a}{\norm a}-\frac{b}{\norm b}}\leq2\min\left(1,\frac{\norm{a-b}}{\norm a}\right)$.
We set
$a\coloneqq\Gamma(\mathbf x)(u-x_\ell)$ and
$b\coloneqq\Gamma(\mathbf x)(v-x_\ell)$.
Since $u\notin\aff(\mathbf x)$, we have
$\norm a=\dist(u,\aff(\mathbf x))>0$.
The projection $\Gamma(\mathbf x)$ is orthogonal, so
$\norm{a-b}\leq\norm{u-v}$.
For $b\ne0$, the inequality for normalized vectors and
\eqref{eq:little_g_grad} give the first line of
\eqref{eq:osc-majorant}.
If $b=0$, the chosen subgradient satisfies
$\nabla_\mu g(\mathbf x;v)=0$. In this case,
\[
 \dist(u,\aff(\mathbf x))=\norm a
 =\norm{\Gamma(\mathbf x)(u-v)}
 \leq\norm{u-v}
 \leq t\prod_{i=1}^{\ell}(1+\norm{x_i}).
\]
The minimum in the first line of \eqref{eq:osc-majorant} is
therefore one, and the left-hand side of that inequality equals
$m_{\ell-1}(\Simpl(\mathbf x))/\ell$.
This proves the first line also when $b=0$.

For $\ell=d$, there exists a unit vector $N\in\R^d$ normal to
the hyperplane $\aff(\mathbf x)$, and $\Gamma(\mathbf x)=NN^\top$.
With $\sgn(0)\coloneqq0$, \eqref{eq:little_g_grad} gives, for every
$\mu\in\R^d$,
\[
 \nabla_\mu g(\mathbf x;\mu)
 =\frac{m_{d-1}(\Simpl(\mathbf x))}{d}N\sgn(N^\top(\mu-x_d)).
\]
If $\dist(u,\aff(\mathbf x))>t\prod_{i=1}^{d}(1+\norm{x_i})$,
the assumed bound on $\norm{v-u}$ implies that $u$ and $v$ lie
in the same open half-space bounded by $\aff(\mathbf x)$, so
$\nabla_\mu g(\mathbf x;u)=\nabla_\mu g(\mathbf x;v)$.
The norm of the difference is always at most
$2m_{d-1}(\Simpl(\mathbf x))/d$.
This proves the second line of \eqref{eq:osc-majorant}.

For $\ell<d$, we first prove \eqref{eq:weighted-inverse}
and then deduce the first inequality in \eqref{eq:osc-expectation-oja}.
The bound \eqref{eq:weighted-inverse} is
Assumption~\ref{ass:inverse-alternative} written in terms of
the distance to $\aff(\mathbf X)$.
Under Assumption~\ref{ass:density}, \eqref{eq:moment-L} bounds
the contribution of $\{\dist(u,\aff(\mathbf X))\geq1\}$
to the expectation in \eqref{eq:weighted-inverse}, uniformly in $u$.
Lemma~\ref{lem:tube}, with $a=1$, gives a constant $C\geq0$ such that,
for every $u\in\R^d$ and $j\geq0$, the part for which
$2^{-j-1}<\dist(u,\aff(\mathbf X))\leq2^{-j}$ is bounded by
\[
 2^{j+1}\E\left[
 \frac{m_{\ell-1}(\Simpl(\mathbf X))}{\ell}
 \prod_{i=1}^{\ell}(1+\norm{X_i})
 \ind_{\{\dist(u,\aff(\mathbf X))\leq2^{-j}\}}\right]
 \leq C2^{-j(d-\ell)}.
\]
Since $d-\ell\geq1$, summing over $j\geq0$ proves
\eqref{eq:weighted-inverse} under Assumption~\ref{ass:density}.
In both cases, bounding the minimum by its second argument gives
the first inequality in \eqref{eq:osc-expectation-oja}.

For $\ell=d$, Assumption~\ref{ass:density} holds.
We fix $a\in(1,\nu-1)$ and partition according to
$2^j\leq\prod_{i=1}^d(1+\norm{X_i})<2^{j+1}$, for $j\geq0$.
On each of these events, the radius $t\prod_{i=1}^d(1+\norm{X_i})$
is at most $t2^{j+1}$. Applying Lemma~\ref{lem:tube}
with $t$ replaced by $t2^{j+1}$ gives constants $C,C'\geq0$ such that
\begin{align*}
 &\E\left[\frac{m_{d-1}(\Simpl(\mathbf X))}{d}
 \ind_{\{\dist(u,\aff(\mathbf X))\leq t\prod_{i=1}^d(1+\norm{X_i})\}}\right]\\
 &\quad\leq\sum_{j=0}^{\infty}2^{-ja}
 \E\left[\frac{m_{d-1}(\Simpl(\mathbf X))}{d}
 \prod_{i=1}^d(1+\norm{X_i})^a
 \ind_{\{\dist(u,\aff(\mathbf X))\leq t2^{j+1}\}}\right]\\
 &\quad\leq Ct\sum_{j=0}^{\infty}2^{-j(a-1)}\leq C't.
\end{align*}
This proves the second inequality in \eqref{eq:osc-expectation-oja}.
\end{proof}

\section{Probabilistic and matrix estimates}\label{sec:prob-tools}
In this section, $\left(\F_k\right)_{k\geq0}$ is an arbitrary filtration.
\subsection{Martingale differences with a fixed delay}
\begin{lemma}\label{lem:series}
Let $d$ and $m$ be positive integers, and let
$\left(U_k\right)_{k\geq1}$ be a sequence of square-integrable random
vectors in $\R^d$, adapted to a filtration
$\left(\F_k\right)_{k\geq0}$.
Suppose that, for every $k\geq m$,
\[
 \E\left[U_k\mid\F_{k-m}\right]=0\quad\as
\]
Then the following statements hold.
\begin{enumerate}[label=(\roman*)]
\item If $\sum_{k=1}^{\infty}\E\left[\norm{U_k}^2\right]<\infty$,
then $\sum_{k=1}^{\infty}U_k$ converges almost surely and in $L^2$.
\item If
$\sum_{k=m}^{\infty}\E\left[\norm{U_k}^2\mid\F_{k-m}\right]<\infty$
almost surely, then $\sum_{k=1}^{\infty}U_k$ converges almost surely.
\item If
\[
 \sum_{k=m}^{\infty}\frac1k
 \E\left[\norm{U_k}^2\mid\F_{k-m}\right]<\infty\quad\as,
\]
then $\sum_{k=1}^nU_k=o\left(\sqrt n\right)$ almost surely.
\end{enumerate}
\end{lemma}
\begin{proof}
For (i) and (ii), we fix $a\in\{0,\ldots,m-1\}$ and set
$k_j\coloneqq a+jm$ for $j\geq0$.
Since $k_j-m=k_{j-1}$ for $j\geq1$, the partial sums of
$\sum_{j\geq1}U_{k_j}$ form a square-integrable martingale with
respect to $\left(\F_{k_j}\right)_{j\geq0}$.

Under the assumption of (i), orthogonality gives, for every $N\geq1$,
\[
 \E\left[\norm{\sum_{j=1}^NU_{k_j}}^2\right]
 =\sum_{j=1}^N\E\left[\norm{U_{k_j}}^2\right]
 \leq\sum_{k=1}^{\infty}\E\left[\norm{U_k}^2\right]<\infty.
\]
The partial sums are therefore bounded in $L^2$.
Applying \citet[Proposition~1.3.7(2), p.~16]{Duflo} to each
coordinate gives convergence almost surely and in $L^2$.

For (ii), the assumption gives
\[
 \sum_{j=1}^{\infty}
 \E\left[\norm{U_{k_j}}^2\mid\F_{k_{j-1}}\right]
 \leq\sum_{k=m}^{\infty}
 \E\left[\norm{U_k}^2\mid\F_{k-m}\right]<\infty\quad\as
\]
Applying \citet[Theorem~1.3.15(1), p.~20]{Duflo} to each
coordinate of the martingale gives almost sure convergence of
$\sum_{j\geq1}U_{k_j}$.

Summing over the $m$ residue classes gives convergence of
$\sum_{k=m}^{\infty}U_k$, almost surely and also in $L^2$ in case (i).
Adding $\sum_{k=1}^{m-1}U_k$ proves (i) and (ii).

Under the assumption in (iii), assertion (ii) applied to
$\left(U_k/\sqrt k\right)_{k\geq1}$ shows that
$\sum_{k=1}^{\infty}U_k/\sqrt k$ converges almost surely.
Kronecker's lemma \citep[Lemma~1.3.14, p.~19]{Duflo}, applied to each
coordinate, gives (iii).
\end{proof}

\subsection{Products of matrices}\label{sec:weights}
In this subsection, $H$ is a symmetric positive definite
$d\times d$ matrix, $c_\gamma$ is a positive real number,
and $0<\gamma<1$.
We set $\gamma_k\coloneqq c_\gamma k^{-\gamma}$ for $k\geq1$
and $\lambda\coloneqq\lambda_{\min}(H)$.
As in Appendix~\ref{sec:asymptotics}, $n_H$ is the smallest
positive integer such that $\gamma_{n_H}\norm H_{\op}\leq1/2$.
For $n\geq k\geq n_H-1$, we set
\[
 \Pi_{n,k}\coloneqq\prod_{j=k+1}^n(I_d-\gamma_jH)
 \qquad\text{and}\qquad
 p_{n,k}\coloneqq\prod_{j=k+1}^n(1-\lambda\gamma_j).
\]

\begin{lemma}\label{lem:deterministic-weights}
For every $n\geq n_H$,
\[
 \sum_{k=n_H}^n\gamma_k\norm{\Pi_{n,k}}_{\op}\leq\frac1{\lambda}.
\]
For every $a\geq0$ and $b\in\R$,
\begin{align}
 \sum_{k=\max\{n_H,2\}}^n
 \gamma_k\norm{\Pi_{n,k}}_{\op}\,k^{-a}\left(\log k\right)^b
 &=O\left(n^{-a}\left(\log n\right)^b\right),\label{eq:weights-power}\\
 \sum_{k=n_H}^n\gamma_k^2\norm{\Pi_{n,k}}_{\op}^2
 &=O\left(\gamma_n\right).\label{eq:weights-square}
\end{align}
For every fixed integer $K\geq n_H-1$, there exists $c>0$ such that
\[
 \max_{n_H-1\leq k\leq K}\norm{\Pi_{n,k}}_{\op}
 =O\left(e^{-c n^{1-\gamma}}\right)
 \quad\text{as }n\to\infty.
\]
\end{lemma}
\begin{proof}
Diagonalizing $H$ gives, for $n\geq k\geq n_H-1$,
\begin{equation}\label{eq:product-contraction}
 \norm{\Pi_{n,k}}_{\op}\leq p_{n,k}
 \leq\exp\left(-\lambda\sum_{j=k+1}^n\gamma_j\right)\leq1.
\end{equation}
Summing the identity
$\lambda\gamma_kp_{n,k}=p_{n,k}-p_{n,k-1}$ gives, by telescoping,
\[
 \lambda\sum_{k=n_H}^n\gamma_kp_{n,k}=1-p_{n,n_H-1}\leq1,
\]
proving the first inequality.
For $k\leq n/2$, we have
$\sum_{j=k+1}^n\gamma_j\geq(n-k)c_\gamma n^{-\gamma}
\geq(c_\gamma/2)n^{1-\gamma}$.
Hence $\norm{\Pi_{n,k}}_{\op}\leq
\exp\left(-\lambda c_\gamma n^{1-\gamma}/2\right)$.
The other factors in each summand of \eqref{eq:weights-power}
grow at most polynomially, so the sum over $k\leq n/2$ is
$o\left(n^{-a}\left(\log n\right)^b\right)$.
For $k>n/2$ in the sum \eqref{eq:weights-power}, we have
the following uniform bounds for $n\geq4$:
\[
 2^{-|b|}
 \leq\frac{k^{-a}(\log k)^b}{n^{-a}(\log n)^b}
 =\left(\frac nk\right)^a
  \left(\frac{\log k}{\log n}\right)^b
 \leq2^{a+|b|}.
\]
The first inequality of the lemma gives
\begin{align*}
 \sum_{\substack{n_H\leq k\leq n\\k>n/2}}
 \gamma_k\norm{\Pi_{n,k}}_{\op}\,k^{-a}(\log k)^b
 &\leq2^{a+|b|}n^{-a}(\log n)^b
 \sum_{k=n_H}^n\gamma_k\norm{\Pi_{n,k}}_{\op}\\
 &\leq\frac{2^{a+|b|}}{\lambda}\,n^{-a}(\log n)^b.
\end{align*}
Together with the bound for $k\leq n/2$, this proves
\eqref{eq:weights-power}.

To prove \eqref{eq:weights-square}, the sum over $k\leq n/2$ is again exponentially small. For $k>n/2$, $\gamma_k\leq2^\gamma\gamma_n$ and $\norm{\Pi_{n,k}}_{\op}\leq1$, so
\[
 \sum_{k>n/2}\gamma_k^2\norm{\Pi_{n,k}}_{\op}^2
 \leq2^\gamma\gamma_n\sum_{k>n/2}\gamma_k\norm{\Pi_{n,k}}_{\op}
 \leq\frac{2^\gamma}{\lambda}\gamma_n.
\]
Finally, for a fixed integer $K\geq n_H-1$, every index
$n_H-1\leq k\leq K$ satisfies $k\leq n/2$ for all sufficiently
large $n$. The exponential bound for $k\leq n/2$ therefore
proves the last assertion.
\end{proof}

\begin{lemma}
\label{lem:perturbed-weights}
Let $n_0\geq n_H$ be an integer, let
$\left(u_n\right)_n$ and $\left(\kappa_n\right)_n$ be sequences
of nonnegative real numbers, and let $\left(v_n\right)_{n>n_0}$
be a sequence of positive real numbers.
Suppose that, for every $n>n_0$,
\[
 u_n\leq v_n+
 \sum_{k=n_0+1}^n
 \gamma_k\norm{\Pi_{n,k}}_{\op}\kappa_{k-1}u_{k-1}.
\]
If $\left(\kappa_n\right)_n$ converges to zero and
\[
 \frac{v_{n-1}}{v_n}=1+o(\gamma_n),
\]
then $u_n=O(v_n)$.
The ratio condition holds in particular for
$v_n=c n^{-a}(\log n)^b$, with $c>0$, $a\geq0$ and $b\in\R$.
\end{lemma}

\begin{proof}
We set
\[
 s_n\coloneqq
 \sum_{k=n_0+1}^n
 \gamma_kp_{n,k}\kappa_{k-1}u_{k-1}.
\]
For $n>n_0$, \eqref{eq:product-contraction} gives
$u_n\leq v_n+s_n$, hence
\begin{align*}
 s_{n+1}
 &= (1-\lambda\gamma_{n+1})s_n
   +\gamma_{n+1}\kappa_nu_n\\
 &\leq
 \bigl(1-(\lambda-\kappa_n)\gamma_{n+1}\bigr)s_n
 +\gamma_{n+1}\kappa_nv_n.
\end{align*}
Since $v_n/v_{n+1}=1+o(\gamma_{n+1})$ and $\kappa_n\to0$,
there exists an integer $N>n_0$ such that, for every $n\geq N$,
\begin{align*}
 \frac{s_{n+1}}{v_{n+1}}
 &\leq\frac{v_n}{v_{n+1}}
 \left(
 \bigl(1-(\lambda-\kappa_n)\gamma_{n+1}\bigr)
       \frac{s_n}{v_n}
 +\gamma_{n+1}\kappa_n
 \right)\\
 &\leq
 \left(1-\frac{\lambda}{2}\gamma_{n+1}\right)
       \frac{s_n}{v_n}
 +\frac{\lambda}{2}\gamma_{n+1}.
\end{align*}
By induction, for every $n\geq N$,
\[
 \frac{s_n}{v_n}
 \leq \max\left\{\frac{s_N}{v_N},1\right\}.
\]
The conclusion follows from $u_n\leq v_n+s_n$.
For $v_n=c n^{-a}(\log n)^b$, we have
$v_{n-1}/v_n=1+O(n^{-1})=1+o(\gamma_n)$ since $\gamma<1$.
\end{proof}

The following lemma applies a bound for weighted martingale sums
from \citet[Corollary~6.4.25, p.~211]{Duflo} to conditional
centering with a fixed delay.

\begin{lemma}\label{lem:convolution}
Suppose that $1/2<\gamma<1$.
Let $m\geq1$ be an integer, and let $\left(U_k\right)_{k\geq1}$ be a sequence of $\F_k$-measurable random vectors in $\R^d$. Suppose that there exist $q>2/\gamma$ and a constant $C\geq0$ such that, almost surely for every $k\geq m$,
\[
 \E\left[U_k\mid\F_{k-m}\right]=0,
 \quad\text{and}\quad \E\left[\norm{U_k}^q\mid\F_{k-m}\right]\leq C.
\]
Then
\begin{equation*}
 \norm{\sum_{k=n_H}^n\gamma_k\Pi_{n,k}U_k}
 =O\left(\sqrt{\gamma_n\log n}\right)\quad\as
\end{equation*}
\end{lemma}
\begin{proof}
We set $n_0\coloneqq\max\{n_H,m\}$.
By Lemma~\ref{lem:deterministic-weights}, there exists $c>0$ such that
\[
 \norm{\sum_{k=n_H}^{n_0-1}\gamma_k\Pi_{n,k}U_k}
 =O\left(e^{-c n^{1-\gamma}}\right)\quad\as
\]
It remains to bound the sum over $n_0\leq k\leq n$.
To apply Duflo's corollary, we decrease $q$ and increase $C$
if necessary so that
$2/\gamma<q<4$; the conditional moment bound is preserved
by Jensen's inequality.

For $a\in\{0,\ldots,m-1\}$ and $k\geq n_0$, we set
$U_k^{(a)}\coloneqq U_k\ind_{\{k\equiv a\pmod m\}}$.
For $k\geq n_0-1$, we denote by $\rho_a(k)$ the largest positive integer
$j\leq k$ such that $j\equiv a\pmod m$, or zero if there is none.
For $k\geq n_0$, we have $\rho_a(k-1)=k-m$ if
$k\equiv a\pmod m$, and $U_k^{(a)}=0$ otherwise.
Thus $\left(U_k^{(a)}\right)_{k\geq n_0}$ is adapted to
$\left(\F_{\rho_a(k)}\right)_{k\geq n_0-1}$ and, for every $k\geq n_0$,
\[
 \E[U_k^{(a)}\mid\F_{\rho_a(k-1)}]=0,
 \quad\text{and}\quad
 \E[\norm{U_k^{(a)}}^q\mid\F_{\rho_a(k-1)}]\leq C.
\]
Conditional Jensen's inequality also gives
$\E[\norm{U_k^{(a)}}^2\mid\F_{\rho_a(k-1)}]\leq C^{2/q}$,
as required by Duflo's corollary.

We fix a unit eigenvector $e$ of $H$, with eigenvalue $h>0$,
and set $b_n\coloneqq\prod_{i=n_0}^n(1-h\gamma_i)>0$ for $n\geq n_0-1$.
For $n\geq n_0$, we define
\[
 M_n\coloneqq\sum_{k=n_0}^n\frac{\gamma_k}{b_k}e^\top U_k^{(a)},
 \qquad
 A_n\coloneqq\sum_{k=n_0}^n\frac{\gamma_k^2}{b_k^2}.
\]
Since $\sum_n\gamma_n^2<\infty$, expanding $\log(1-h\gamma_i)$ gives
$-\log b_n\sim hc_\gamma n^{1-\gamma}/(1-\gamma)$.
Moreover, $\gamma_{n-1}/\gamma_n=1+O(n^{-1})=1+o(\gamma_n)$, so
\[
 \frac{\gamma_n}{b_n^2}-\frac{\gamma_{n-1}}{b_{n-1}^2}
 =\frac{\gamma_n-\gamma_{n-1}(1-h\gamma_n)^2}{b_n^2}
 \sim2h\frac{\gamma_n^2}{b_n^2}.
\]
The Stolz--Ces\`aro theorem gives
$A_n\sim\gamma_n/(2hb_n^2)$. Taking logarithms gives
\[
 \log A_n=-2\log b_n+\log\gamma_n-\log(2h)+o(1)
 \sim\frac{2hc_\gamma}{1-\gamma}n^{1-\gamma}.
\]
Thus $A_n\to\infty$ and $\log\log A_n=O(\log n)$.
Since $\gamma q/2>1$, we have
\[
 \sum_{n=n_0}^{\infty}
 \left(\frac{\gamma_n^2}{b_n^2A_n}\right)^{q/2}<\infty.
\]
We apply Corollary~6.4.25 of \citet[p.~211]{Duflo} to
$\left(e^\top U_k^{(a)}\right)_{k\geq n_0}$, with
$\alpha=q/2-1$. The weights $\Phi$ and their bounds $T$
in that corollary are both $\gamma_k/b_k$, and the sum of
their squares is $A_n$. Hence
\[
 M_n=O\left(\sqrt{A_n\log\log A_n}\right)\quad\as
\]
Since $e^\top\Pi_{n,k}=(b_n/b_k)e^\top$, it follows that
\[
 \left|e^\top\sum_{k=n_0}^n\gamma_k\Pi_{n,k}U_k^{(a)}\right|
 =b_n|M_n|
 =O\left(\sqrt{\gamma_n\log n}\right)\quad\as
\]
Summing over the $m$ residue classes and using an orthonormal
eigenbasis of $H$ proves the lemma.
\end{proof}

\subsection{Strictly stationary sequences with finite dependence}\label{sec:finite-clt}
A sequence $\left(Y_i\right)_{i\geq1}$ is called $m$-dependent,
for an integer $m\geq0$, if for every $a\geq1$ the sigma-fields
generated by $\left(Y_i\right)_{i\leq a}$ and
$\left(Y_i\right)_{i\geq a+m+1}$ are independent.

\begin{proposition}[{\citet[Proposition~2.1.7, p.~43]{Duflo}}]\label{th:duflo}
Let $\left(Y_i\right)_{i\geq1}$ be a strictly stationary, centered,
square-integrable sequence of random vectors in $\R^d$.
Suppose that it is $m$-dependent for some integer $m\geq1$.
For $0\leq j\leq m$, we set
\[
 \Lambda_j\coloneqq\E\left[Y_1Y_{1+j}^\top\right]
 \quad\text{and}\quad
 \Sigma_Y\coloneqq\Lambda_0+
 \sum_{j=1}^m\left(\Lambda_j+\Lambda_j^\top\right).
\]
Then
\[
 \frac1{\sqrt n}\sum_{i=1}^nY_i
 \xrightarrow[n\to\infty]{\hspace{.5em}\mathcal L\hspace{.5em}}
 \mathcal N_d\left(0,\Sigma_Y\right).
\]
\end{proposition}

In the independent case, the proposition applies with $m=1$
and $\Sigma_Y=\Lambda_0$.

\begin{lemma}[Central limit theorem for sums with matrix coefficients]
\label{lem:matrix-mdep-clt}
Let $d,n_0\geq1$ and $m\geq0$ be integers, and let
$\left(Y_k\right)_{k\geq n_0}$ be a strictly stationary, centered, square-integrable,
$m$-dependent sequence of random vectors in $\R^d$.
Let $H$ be a symmetric positive definite $d\times d$ matrix, let
$c_\gamma>0$ and $0<\gamma<1$, and set
$\gamma_k\coloneqq c_\gamma k^{-\gamma}$ for $k\geq1$.
Suppose that $\gamma_{n_0}\norm H_{\op}\leq1/2$.
For $n\geq k\geq n_0$, define
$\Pi_{n,k}\coloneqq\prod_{r=k+1}^n(I_d-\gamma_rH)$, with $\Pi_{n,n}=I_d$.
For $0\leq j\leq m$, set
$\Lambda_j\coloneqq\E\left[Y_{n_0}Y_{n_0+j}^\top\right]$, and set
\[
 \Sigma_Y\coloneqq\Lambda_0+
 \sum_{j=1}^m\left(\Lambda_j+\Lambda_j^\top\right).
\]
Then
\[
 \frac1{\sqrt{\gamma_n}}\sum_{k=n_0}^n\gamma_k\Pi_{n,k}Y_k
 \xrightarrow[n\to\infty]{\hspace{.5em}\mathcal L\hspace{.5em}}\mathcal N_d\left(0,
 \int_0^\infty e^{-tH}\Sigma_Ye^{-tH}\,dt\right).
\]
\end{lemma}
\begin{proof}
We apply Orey's central limit theorem for $m$-dependent triangular
arrays, in the form stated in
\citet[Theorem~13.1, pp.~67--68]{Rosen1967}.

\smallskip\noindent\emph{Convergence of the covariance.}\enspace
For $n\geq n_0$, we set
\[
 Z_n\coloneqq\sum_{k=n_0}^n\gamma_k\Pi_{n,k}Y_k
 \quad\text{and}\quad Q_n\coloneqq\Var(Z_n).
\]
For $n>n_0+m$, the recursion
$Z_n=(I_d-\gamma_nH)Z_{n-1}+\gamma_nY_n$,
strict stationarity and $m$-dependence give
\[
 Q_n=(I_d-\gamma_nH)Q_{n-1}(I_d-\gamma_nH)+\gamma_n^2\Lambda_0
 +\gamma_n\sum_{j=1}^m\gamma_{n-j}
 \left(\Pi_{n,n-j}\Lambda_j+\Lambda_j^\top\Pi_{n,n-j}\right).
\]
For $1\leq j\leq m$, we have
$\gamma_{n-j}/\gamma_n=1+O(n^{-1})$ and
$\Pi_{n,n-j}=I_d+O(\gamma_n)$ in operator norm.
Since $n^{-1}=o(\gamma_n)$, it follows that
\[
 Q_n=(I_d-\gamma_nH)Q_{n-1}(I_d-\gamma_nH)
       +\gamma_n^2\Sigma_Y+O(\gamma_n^3).
\]
We set $V\coloneqq\int_0^\infty e^{-tH}\Sigma_Ye^{-tH}\,dt$.
Since $H$ is positive definite, the integral is finite, and integrating
the derivative of $e^{-tH}\Sigma_Ye^{-tH}$ gives $HV+VH=\Sigma_Y$.
Together with $\gamma_{n-1}-\gamma_n=O(\gamma_n/n)$, this yields
\[
 Q_n-\gamma_nV
 =(I_d-\gamma_nH)(Q_{n-1}-\gamma_{n-1}V)(I_d-\gamma_nH)
 +O(\gamma_n^3+\gamma_n/n).
\]
We set $N\coloneqq n_0+m$.
The contraction bound \eqref{eq:product-contraction} and
Lemma~\ref{lem:deterministic-weights}, with $a=2\gamma$ and $a=1$
in \eqref{eq:weights-power}, give, by iteration, for some constant
$C\geq0$ and every $n>N$,
\[
\begin{aligned}
 \norm{Q_n-\gamma_nV}_{\op}
 &\leq\norm{\Pi_{n,N}}_{\op}^2\norm{Q_N-\gamma_NV}_{\op}
 +C\sum_{k=N+1}^n\gamma_k\norm{\Pi_{n,k}}_{\op}^2
                 (\gamma_k^2+k^{-1})\\
 &=O(\gamma_n^2+n^{-1}).
\end{aligned}
\]
Consequently,
\[
 \Var\left(\frac{Z_n}{\sqrt{\gamma_n}}\right)
 =\frac{Q_n}{\gamma_n}
 =V+O\left(\gamma_n+\frac1{n\gamma_n}\right)
 \xrightarrow[n\to\infty]{\hspace{1.5em}}V.
\]

\smallskip\noindent\emph{Lindeberg condition and limit distribution.}\enspace
For $n_0\leq k\leq n$, we set
\[
 W_{n,k}\coloneqq\frac{\gamma_k}{\sqrt{\gamma_n}}\Pi_{n,k}
 \quad\text{and}\quad
 \delta_n\coloneqq\max_{n_0\leq k\leq n}\norm{W_{n,k}}_{\op}.
\]
The bound \eqref{eq:weights-square} gives
$\sum_{k=n_0}^n\norm{W_{n,k}}_{\op}^2=O(1)$.
Applying \eqref{eq:product-contraction} separately to $k\leq n/2$
and $k>n/2$ gives
\[
 \delta_n\leq
 \frac{\gamma_{n_0}}{\sqrt{\gamma_n}}
 \exp\left(-\frac{c_\gamma\lambda_{\min}(H)}2n^{1-\gamma}\right)
 +2^\gamma\sqrt{\gamma_n}
 \xrightarrow[n\to\infty]{\hspace{1.5em}}0.
\]
We fix a unit vector $v\in\R^d$ and set
$X_{n,k}\coloneqq v^\top W_{n,k}Y_k$ for $n_0\leq k\leq n$.
These variables are centered and form an $m$-dependent triangular
array. Strict stationarity and square integrability give, for every
$\varepsilon>0$,
\[
 \sum_{k=n_0}^n\E\left[X_{n,k}^2\ind_{\{|X_{n,k}|>\varepsilon\}}\right]
 \leq\left(\sum_{k=n_0}^n\norm{W_{n,k}}_{\op}^2\right)
 \E\left[\norm{Y_{n_0}}^2
 \ind_{\{\norm{Y_{n_0}}>\varepsilon/\delta_n\}}\right]
 \xrightarrow[n\to\infty]{\hspace{1.5em}}0.
\]
Moreover,
\[
 \sum_{k=n_0}^n\Var(X_{n,k})
 \leq\E\left[\norm{Y_{n_0}}^2\right]
      \sum_{k=n_0}^n\norm{W_{n,k}}_{\op}^2=O(1).
\]
If $v^\top Vv>0$, the variances of the row sums
$\sigma_n^2\coloneqq v^\top Q_nv/\gamma_n$ converge to $v^\top Vv$.
For all sufficiently large $n$, the normalized array
$\left(X_{n,k}/\sigma_n\right)_{n_0\leq k\leq n}$ satisfies the hypotheses of
\citet[Theorem~13.1]{Rosen1967}, and
\[
 v^\top\frac{Z_n}{\sqrt{\gamma_n}}
 =\sum_{k=n_0}^nX_{n,k}
 \xrightarrow[n\to\infty]{\hspace{.5em}\mathcal L\hspace{.5em}}
 \mathcal N\left(0,v^\top Vv\right).
\]
If $v^\top Vv=0$, the same projection converges to zero in $L^2$.
Since $v$ is arbitrary, the Cram\'er--Wold theorem gives the stated
vector convergence.
\end{proof}

\section{Hoeffding decomposition}\label{sec:batch-auxiliary}

We recall the classical Hoeffding decomposition \citep{Hoeffding1961}
in the form used in Appendix~\ref{sec:batch}.

\begin{lemma}\label{lem:hoeffding}
Let $\left(X_i\right)_{i\geq1}$ be independent and identically
distributed random vectors in $\R^d$, and let $\ell\geq1$ be an integer.
We set $\mathbf X\coloneqq(X_1,\ldots,X_\ell)$.
Let $h:(\R^d)^\ell\to\R^d$ be measurable and symmetric, with
$\E[h(\mathbf X)]=0$ and $\E[\norm{h(\mathbf X)}^2]<\infty$.
For every $1\leq j\leq\ell$, there exists a measurable symmetric
function $\psi_j:(\R^d)^j\to\R^d$ satisfying
$\E[\norm{\psi_j(X_1,\ldots,X_j)}^2]<\infty$ and
\begin{equation}\label{eq:canonical-kernels}
 \psi_j(X_1,\ldots,X_j)
 =\sum_{A\subseteq\{1,\ldots,j\}}(-1)^{j-|A|}
       \E[h(\mathbf X)\mid X_i,\ i\in A]\quad\as
\end{equation}
For a finite set $A=\{i_1<\cdots<i_j\}$ of positive integers,
we set $X_A\coloneqq(X_{i_1},\ldots,X_{i_j})$. Then
\begin{equation}\label{eq:hoeffding-expansion}
 h(\mathbf X)=\sum_{\varnothing\ne A\subseteq\{1,\ldots,\ell\}}
                  \psi_{|A|}(X_A)\quad\as
\end{equation}
For every finite nonempty set $A$ of positive integers with
$|A|\leq\ell$ and every $k\in A$, we have
\[
 \E[\psi_{|A|}(X_A)\mid X_i,\ i\in A\setminus\{k\}]=0\quad\as
\]
When $A=\{k\}$, this identity gives $\E[\psi_1(X_k)]=0$.
Moreover, for distinct finite nonempty sets $A,B$ of positive
integers with $|A|,|B|\leq\ell$, we have
\[
 \E[\psi_{|A|}(X_A)\psi_{|B|}(X_B)^\top]=0.
\]
\end{lemma}
\begin{proof}
We denote the common law of the observations by $P$.
For $1\leq r<\ell$ and $\mathbf x\in(\R^d)^r$, we define
\[
 h_r(\mathbf x)\coloneqq
 \begin{cases}
 \displaystyle\int h(\mathbf x,\mathbf y)\,
       P^{\otimes(\ell-r)}(d\mathbf y),
 &\text{if }\displaystyle\int\norm{h(\mathbf x,\mathbf y)}\,
       P^{\otimes(\ell-r)}(d\mathbf y)<\infty,\\[0.8em]
 0,&\text{otherwise}.
 \end{cases}
\]
We also set $h_\ell\coloneqq h$.
The functions $h_r$ are measurable by Fubini's theorem.
The symmetry of $h$ implies that each $h_r$ is symmetric.
For every $\varnothing\ne A\subseteq\{1,\ldots,\ell\}$,
independence and symmetry of $h$ give
\[
 h_{|A|}(X_A)
 =\E[h(\mathbf X)\mid X_i,\ i\in A]\quad\as
\]
For $1\leq j\leq\ell$ and $(x_1,\ldots,x_j)\in(\R^d)^j$,
we write $x_A\coloneqq(x_{i_1},\ldots,x_{i_r})$ when
$A=\{i_1<\cdots<i_r\}\subseteq\{1,\ldots,j\}$, and define
\[
 \psi_j(x_1,\ldots,x_j)\coloneqq
 \sum_{\varnothing\ne A\subseteq\{1,\ldots,j\}}
       (-1)^{j-|A|}h_{|A|}(x_A).
\]
The functions $\psi_j$ are measurable and symmetric.
Since $\E[h(\mathbf X)]=0$, this definition gives
\eqref{eq:canonical-kernels}.
Conditional Jensen's inequality and the finite sum give
$\E[\norm{\psi_j(X_1,\ldots,X_j)}^2]<\infty$.

For every $\varnothing\ne B\subseteq\{1,\ldots,\ell\}$,
the binomial identity gives
\[
 \sum_{B\subseteq A\subseteq\{1,\ldots,\ell\}}
       (-1)^{|A|-|B|}
 =\sum_{r=0}^{\ell-|B|}(-1)^r\binom{\ell-|B|}{r}
 =\begin{cases}
   1,&\text{if }B=\{1,\ldots,\ell\},\\
   0,&\text{if }B\ne\{1,\ldots,\ell\}.
  \end{cases}
\]
Consequently,
\[
\begin{aligned}
 \sum_{\varnothing\ne A\subseteq\{1,\ldots,\ell\}}
       \psi_{|A|}(X_A)
 &=
 \sum_{\varnothing\ne A\subseteq\{1,\ldots,\ell\}}
       \sum_{\varnothing\ne B\subseteq A}
       (-1)^{|A|-|B|}h_{|B|}(X_B)\\
 &=
 \sum_{\varnothing\ne B\subseteq\{1,\ldots,\ell\}}
       \left(
       \sum_{B\subseteq A\subseteq\{1,\ldots,\ell\}}
       (-1)^{|A|-|B|}
       \right)h_{|B|}(X_B)\\
 &=h_\ell(\mathbf X)=h(\mathbf X).
\end{aligned}
\]
This proves \eqref{eq:hoeffding-expansion}.

For $1\leq j\leq\ell$ and $B\subseteq\{1,\ldots,j-1\}$,
independence of $X_j$ from $(X_1,\ldots,X_{j-1})$ gives
\[
 \E\Bigl[\E[h(\mathbf X)\mid X_i,\ i\in B\cup\{j\}]
       \mid X_1,\ldots,X_{j-1}\Bigr]
 =\E[h(\mathbf X)\mid X_i,\ i\in B]\quad\as
\]
The terms indexed by $B$ and $B\cup\{j\}$ in
\eqref{eq:canonical-kernels} have opposite signs. Hence
\[
 \E[\psi_j(X_1,\ldots,X_j)\mid X_1,\ldots,X_{j-1}]=0\quad\as
\]
Symmetry gives the same identity after conditioning on any $j-1$
of the $j$ arguments. Identical distribution then gives the stated
conditional centering for every set $A$ and every $k\in A$.

For the last assertion, interchanging $A$ and $B$ and taking
the transpose if necessary, we may assume
$A\setminus B\ne\varnothing$.
We fix $k\in A\setminus B$ and set
$\mathcal G\coloneqq\sigma(X_i:i\in(A\cup B)\setminus\{k\})$.
The random vector $X_A$ is independent of
$(X_i)_{i\in B\setminus A}$.
Thus the conditional centering just proved gives
$\E[\psi_{|A|}(X_A)\mid\mathcal G]=0$ almost surely.
Since $\psi_{|B|}(X_B)$ is $\mathcal G$-measurable,
\[
 \E[\psi_{|A|}(X_A)\psi_{|B|}(X_B)^\top]
 =\E\left[\E[\psi_{|A|}(X_A)\mid\mathcal G]\,
             \psi_{|B|}(X_B)^\top\right]=0.
\]
\end{proof}

\section{Properties of the objective functions}\label{sec:verification}

Under Assumption~\ref{ass:density}, we establish the regularity
of $G$ and verify Assumption~\ref{ass:consistency}.
We then prove Proposition~\ref{prop:smooth-density}, which gives
a sufficient condition for the local expansion in
Assumption~\ref{ass:taylor}.

\subsection{Regularity of \texorpdfstring{$g$ and $G$}{g and G}}

\begin{proposition}\label{prop:population}
Under Assumption~\ref{ass:density}, $G$ is finite, strictly
convex, coercive, and differentiable on $\R^d$, with
\begin{equation*}
 \nabla G\left(u\right)=\E\left[\nabla_\mu g\left(\mathbf X;u\right)\right],\quad\text{and}\quad \norm{\nabla G\left(u\right)}\leq\frac1\ell\E\left[m_{\ell-1}\left(\Simpl\left(\mathbf X\right)\right)\right].
\end{equation*}
The gradient $\nabla G$ is Lipschitz.
Moreover, for every $q\in(2,\nu)$, with
$\beta\coloneqq(q-2)/(q-1)\in(0,1)$, there exists a constant
$C\geq0$ such that, for every $u,v\in\R^d$,
\begin{equation}\label{eq:L2-modulus}
 \E\left[\norm{\nabla_\mu g\left(\mathbf X;u\right)-\nabla_\mu g\left(\mathbf X;v\right)}^2\right]
 \leq C\min\left(\norm{u-v}^\beta,1\right).
\end{equation}
In particular, Assumption~\ref{ass:consistency} holds.
\end{proposition}
\begin{proof}
Finiteness and convexity of $G$ follow from
\eqref{eq:product-bound} and \eqref{eq:base-height}.
For every $u\in\R^d$, we have
$u\notin\aff(\mathbf X)$ almost surely, as recalled at the
beginning of Appendix~\ref{sec:geometry}.
The base--height formula \eqref{eq:base-height} therefore
shows that $g(\mathbf X;\cdot)$ is differentiable at $u$
almost surely.
The envelope \eqref{eq:envelope} gives, for $z\ne0$,
\[
 \frac{\abs{g(\mathbf X;u+z)-g(\mathbf X;u)
          -\scal{\nabla_\mu g(\mathbf X;u)}{z}}}{\norm z}
 \leq\frac{2m_{\ell-1}(\Simpl(\mathbf X))}{\ell}.
\]
The dominated convergence theorem applies: the majorant is integrable,
and the quotient tends to zero almost surely as $z\to0$.
Consequently,
\[
 \frac{G(u+z)-G(u)-\scal{\E[\nabla_\mu g(\mathbf X;u)]}{z}}{\norm z}
 \xrightarrow[z\to0]{\hspace{1.5em}}0.
\]
Thus $G$ is differentiable and has the stated gradient.
Since $P$ has a density and a finite first moment,
Theorems~2.2(ii$_a$) and~2.3(i$_a$)
of~\citet{durre_affine-equivariant_2022} give
coercivity of $G$ and existence of a minimizer.

We next prove strict convexity. Let $u,v\in\R^d$ be distinct,
and let $t\in\left(0,1\right)$.
Suppose first that $\ell<d$.
For each $1\leq i\leq\ell$, conditionally on
$X_1,\ldots,X_{i-1}$, the vector $X_i$ has law $P$.
The affine space $\aff(u,v,X_1,\ldots,X_{i-1})$ has dimension
at most $i<d$, so $X_i$ belongs to this space with conditional
probability zero.
Thus the points $u,v,X_1,\ldots,X_\ell$ are affinely independent
almost surely.
Consequently,
$\Gamma(\mathbf X)(u-X_\ell)$ and
$\Gamma(\mathbf X)(v-X_\ell)$ are not collinear almost surely.
Equality in the triangle inequality is therefore
impossible. Together with \eqref{eq:base-height}, this gives
\[
 g\left(\mathbf X;\left(1-t\right)u+tv\right)
 <
 \left(1-t\right)g\left(\mathbf X;u\right)
 +t g\left(\mathbf X;v\right)
 \quad\as
\]
Taking expectations proves strict convexity when $\ell<d$.

Suppose now that $\ell=d$.
To prove strict convexity, we first show that $\aff(\mathbf X)$
separates $u$ and $v$ strictly with positive probability.
The set $E\coloneqq\{x\in\R^d:f(x)>0\}$
has positive Lebesgue measure.
Since $d\geq2$, the line $\aff(u,v)$ has Lebesgue
measure zero.
By the Lebesgue density theorem, there exists a density point
$x_0\in E\setminus\aff(u,v)$ of $E$, so that
\[
 m_d(B(x_0,r)\setminus E)
 =o(r^d)
 \qquad\text{as }r\downarrow0.
\]
The vectors $u-x_0$ and $v-x_0$ are linearly independent.
There exists a unit vector $e_0$ such that
\[
 e_0^\top(u-x_0)>0,
 \quad\text{and}\quad
 e_0^\top(v-x_0)<0.
\]
We choose affinely independent points
$z_1,\ldots,z_d\in e_0^\perp\cap B(0,1/2)$.
Continuity of the inverse Gram matrix shows that $\Gamma$
is continuous on the set of affinely independent tuples.
Whenever $\Gamma(\mathbf x)e_0\ne0$, we set
$e(\mathbf x)\coloneqq\Gamma(\mathbf x)e_0/\norm{\Gamma(\mathbf x)e_0}$.
Since $\Gamma(z_1,\ldots,z_d)e_0=e_0$, the function $e$
gives a continuous unit normal near $(z_1,\ldots,z_d)$.
There therefore exist constants $\varepsilon>0$ and
$0<\delta<1/2$ such that every tuple
$\mathbf x\in\prod_{i=1}^d B(z_i,\delta)$ is affinely
independent, $e(\mathbf x)$ is defined, and
\[
 e(\mathbf x)^\top(u-x_0)>\varepsilon,
 \quad\text{and}\quad
 e(\mathbf x)^\top(v-x_0)<-\varepsilon.
\]
For $0<r<\varepsilon$, we set
\[
 B_i(r)\coloneqq B(x_0+rz_i,\delta r),
 \qquad 1\leq i\leq d.
\]
These balls are contained in $B(x_0,r)$.
For every $\mathbf x\in\prod_{i=1}^d B_i(r)$, we have
$(x_i-x_0)/r\in B(z_i,\delta)$ for every $i$.
Translation and dilation preserve affine independence and
leave $\Gamma$, and hence $e$, unchanged.
Thus $\mathbf x$ is affinely independent and $e(\mathbf x)$
is defined.
Since $\norm{x_d-x_0}<r$, we obtain
\[
 e(\mathbf x)^\top(u-x_d)>\varepsilon-r>0,
 \quad\text{and}\quad
 e(\mathbf x)^\top(v-x_d)<-\varepsilon+r<0.
\]
Thus $\aff(\mathbf x)$ separates $u$ and $v$ strictly.

Since $m_d(B(x_0,r)\setminus E)=o(r^d)$ as $r\downarrow0$,
we fix $r\in(0,\varepsilon)$ such that
\[
 m_d(B(x_0,r)\setminus E)<r^d m_d(B(0,\delta)).
\]
For every $1\leq i\leq d$, we then have
\[
 \begin{aligned}
 m_d(E\cap B_i(r))
 &\geq m_d(B_i(r))-m_d(B(x_0,r)\setminus E)\\
 &=r^d m_d(B(0,\delta))-m_d(B(x_0,r)\setminus E)>0.
 \end{aligned}
\]
Since $f>0$ on $E$, we have $P(B_i(r))>0$ for every $i$.
By independence, there is positive probability that
$X_i\in B_i(r)$ for all $i$.
On this event, $X_1,\ldots,X_d$ are affinely independent,
so $m_{d-1}(\Simpl(\mathbf X))>0$.
Moreover, $e(\mathbf X)^\top(u-X_d)>0$ and
$e(\mathbf X)^\top(v-X_d)<0$.
The base--height formula \eqref{eq:base-height} therefore gives
\[
 \begin{aligned}
 g(\mathbf X;(1-t)u+tv)
 &=\frac{m_{d-1}(\Simpl(\mathbf X))}{d}
   \abs{(1-t)e(\mathbf X)^\top(u-X_d)
          +t e(\mathbf X)^\top(v-X_d)}\\
 &<(1-t)g(\mathbf X;u)+t g(\mathbf X;v).
 \end{aligned}
\]
For every $\mathbf x\in(\R^d)^d$, the function
$g(\mathbf x;\cdot)$ is convex.
Since the strict inequality holds on an event of positive
probability, taking expectations proves strict convexity
of $G$ in the case $\ell=d$.

Strict convexity has now been proved for every
$1\leq\ell\leq d$, so $G$ has a unique minimizer.
Since $\nu>2$, \eqref{eq:moment-L} with $q=2$ gives the
second-moment condition in Assumption~\ref{ass:consistency}.

Since $\prod_{i=1}^{\ell}(1+\norm{X_i})\geq1$,
the bounds \eqref{eq:osc-majorant} and
\eqref{eq:osc-expectation-oja} in
Lemma~\ref{lem:oscillation}, with $t=\norm{u-v}$,
give a constant $C\geq0$ such that, for every $u,v\in\R^d$,
\[
 \E[\norm{\nabla_\mu g(\mathbf X;u)-\nabla_\mu g(\mathbf X;v)}]
 \leq C\norm{u-v}.
\]
The formula for $\nabla G$ therefore shows that $\nabla G$ is Lipschitz.
To prove \eqref{eq:L2-modulus}, we fix $q\in(2,\nu)$ and set
$U\coloneqq\norm{\nabla_\mu g(\mathbf X;u)-\nabla_\mu g(\mathbf X;v)}$.
By the envelope \eqref{eq:envelope} and the moment bound \eqref{eq:moment-L},
\[
 \E[U^q]\leq\left(\frac2\ell\right)^q
                 \E[m_{\ell-1}(\Simpl(\mathbf X))^q]<\infty.
\]
H\"older's inequality, with $\beta=(q-2)/(q-1)$, gives
\begin{equation}\label{eq:interpolation}
 \E[U^2]\leq\left(\E[U]\right)^\beta
                  \left(\E[U^q]\right)^{1/(q-1)}.
\end{equation}
Combining these bounds with
$\E[U^2]\leq4\E[m_{\ell-1}(\Simpl(\mathbf X))^2]/\ell^2$
proves \eqref{eq:L2-modulus}.
\end{proof}

\subsection{Hessians and the local expansion}\label{sec:hessian}
\begin{proposition}\label{prop:hessian-nonoja}
Under Assumption~\ref{ass:density}, if $\ell<d$, then $\nabla G$ is differentiable at every $u\in\R^d$, and
\begin{equation*}
 \nabla^2G\left(u\right)=\frac1\ell\E\left[
 \frac{m_{\ell-1}\left(\Simpl\left(\mathbf X\right)\right)}{\norm{\Gamma\left(\mathbf X\right)\left(u-X_\ell\right)}}
 \left(\Gamma\left(\mathbf X\right)-
 \frac{\Gamma\left(\mathbf X\right)\left(u-X_\ell\right)\left(u-X_\ell\right)^\top\Gamma\left(\mathbf X\right)}
 {\norm{\Gamma\left(\mathbf X\right)\left(u-X_\ell\right)}^2}\right)\right].
\end{equation*}
For every $u\in\R^d$, the denominator is nonzero almost surely.
\end{proposition}
\begin{proof}
Let $u\in\R^d$.
Under Assumption~\ref{ass:density}, $P$ has a finite first moment
and assigns probability zero to every affine subspace of
dimension $\ell-1$.
Lemma~S.1.3(ii) of~\citet{DPsupp} therefore gives
$\norm{\Gamma(\mathbf X)(u-X_\ell)}>0$ almost surely.
Since $\prod_{i=1}^{\ell}(1+\norm{X_i})\geq1$,
\eqref{eq:weighted-inverse} gives
\[
 \E\left[
 \frac{m_{\ell-1}(\Simpl(\mathbf X))}
 {\norm{\Gamma(\mathbf X)(u-X_\ell)}}
 \right]<\infty.
\]
Lemma~S.2.3 of~\citet{DPsupp}, applied with $p=1$ and
$\mu_0=u$, now gives differentiability of $\nabla G$ at $u$
and the stated formula, which is its equation~(S.2.4).
\end{proof}

\begin{proof}[Proof of Proposition~\ref{prop:smooth-density}]
By the definition in \eqref{eq:a-density} with $a=0$, for
$V\in\Gr(\ell-1,d)$ and $y\in V^\perp$, we have
\[
 b_0(V,y)=\frac c\ell\int_{V^\ell}
 m_{\ell-1}(\Simpl(\mathbf z))^{d-\ell+2}
 \prod_{i=1}^{\ell}f(y+z_i)\,d\mathbf z.
\]
Formula~\eqref{eq:BP} gives
\begin{equation}\label{eq:G-flat}
 G\left(u\right)=\int_{\Gr\left(\ell-1,d\right)}\int_{V^\perp}
          b_0\left(V,y\right)\norm{\Gamma_{V^\perp}u-y}\,dy\,d\vartheta\left(V\right).
\end{equation}
\smallskip\noindent\emph{Bounds for $b_0$.}\enspace
We first construct a bounded and integrable function to control
the variation of $b_0$ in $y$.
For $V\in\Gr(\ell-1,d)$ and $y\in V^\perp$, we define
\[
 \widetilde b(V,y)\coloneqq\frac c\ell\int_{V^\ell}
 m_{\ell-1}(\Simpl(\mathbf z))^{d-\ell+2}
 \prod_{i=1}^{\ell}(1+\norm{y+z_i})^{-d-\nu}\,d\mathbf z.
\]
By Assumption~\ref{ass:density},
$0\leq b_0(V,y)\leq C_f^\ell\widetilde b(V,y)$.
Since $y\perp V$, we have $\norm{y+z_i}\geq\norm{z_i}$
for every $z_i\in V$.
Applying \eqref{eq:product-bound} to $\mathbf z$
gives a constant $C>0$ such that
\[
 m_{\ell-1}(\Simpl(\mathbf z))^{d-\ell+2}
 \prod_{i=1}^{\ell}(1+\norm{y+z_i})^{-d-\nu}
 \leq C\prod_{i=1}^{\ell}
 (1+\norm{z_i})^{-(\ell-1)-(\nu-1)}.
\]
The integral over $V^\ell$ of the right-hand side does not
depend on $V$ and is finite, since $\dim V=\ell-1$ and $\nu>1$.
The affine integration formula \eqref{eq:BP} also gives
\[
 \int_{\Gr(\ell-1,d)}\int_{V^\perp}\widetilde b(V,y)\,dy\,d\vartheta(V)
 =\int_{(\R^d)^\ell}\frac{m_{\ell-1}(\Simpl(\mathbf x))}\ell
       \prod_{i=1}^{\ell}(1+\norm{x_i})^{-d-\nu}\,d\mathbf x.
\]
After applying \eqref{eq:product-bound} and increasing $C$ if necessary,
we obtain
\begin{equation}\label{eq:a-tilde}
\begin{gathered}
 \sup_{V\in\Gr(\ell-1,d)}\sup_{y\in V^\perp}\widetilde b(V,y)<\infty,\\
 \int_{\Gr(\ell-1,d)}\int_{V^\perp}\widetilde b(V,y)\,dy\,d\vartheta(V)
 \leq C\left(\int_{\R^d}(1+\norm x)^{-d-\nu+1}\,dx\right)^\ell<\infty.
\end{gathered}
\end{equation}
In particular, $b_0$ is uniformly bounded.

By Assumption~\ref{ass:density} and the bound on $\nabla f$
in the statement of the proposition, there exists a constant
$C>0$ such that, for every $V\in\Gr(\ell-1,d)$,
$y,v\in V^\perp$ and $\mathbf z\in V^\ell$,
\[
\begin{aligned}
 \abs{\prod_{i=1}^{\ell}f(y+v+z_i)
       -\prod_{i=1}^{\ell}f(y+z_i)}
 &=\abs{\int_0^1\sum_{j=1}^{\ell}
       \scal{\nabla f(y+sv+z_j)}{v}
       \prod_{i\ne j}f(y+sv+z_i)\,ds}\\
 &\leq C\norm v\int_0^1\prod_{i=1}^{\ell}
       (1+\norm{y+sv+z_i})^{-d-\nu}\,ds.
\end{aligned}
\]
Formula~\eqref{eq:a-density} and Tonelli's theorem give
\begin{equation}\label{eq:a-translation}
 |b_0(V,y+v)-b_0(V,y)|
 \leq C\norm v\int_0^1\widetilde b(V,y+sv)\,ds.
\end{equation}
By \eqref{eq:a-tilde}, there exists a constant $L_{b_0}\geq0$
such that, for every $V\in\Gr(\ell-1,d)$,
$b_0(V,\cdot)$ is $L_{b_0}$-Lipschitz on $V^\perp$.

\smallskip\noindent\emph{The Hessian when $\ell<d$.}\enspace
For $\ell<d$, Proposition~\ref{prop:hessian-nonoja} and the
affine integration formula~\eqref{eq:BP} give
\begin{equation}\label{eq:H-integral}
 \nabla^2G\left(u\right)=\int\!\!\int
 \frac{b_0\left(V,y\right)}{\norm{\Gamma_{V^\perp}u-y}}
 \left(\Gamma_{V^\perp}-\frac{\left(\Gamma_{V^\perp}u-y\right)\left(\Gamma_{V^\perp}u-y\right)^\top}
 {\norm{\Gamma_{V^\perp}u-y}^2}\right)\,dy\,d\vartheta\left(V\right).
\end{equation}
For $y\ne\Gamma_{V^\perp}u$, the difference of matrices
in parentheses in \eqref{eq:H-integral} is the orthogonal
projection onto $V^\perp\cap(\Gamma_{V^\perp}u-y)^\perp$,
so its operator norm is at most one.
For $u,w\in\R^d$, the change of variables
$y\mapsto y+\Gamma_{V^\perp}w$ in $\nabla^2G(u+w)$ and the bound
\eqref{eq:a-translation} give
\[
\begin{aligned}
 \norm{\nabla^2G(u+w)-\nabla^2G(u)}_{\op}
 &\leq\int\!\!\int
  \frac{|b_0(V,y+\Gamma_{V^\perp}w)-b_0(V,y)|}{\norm{\Gamma_{V^\perp}u-y}}
  \,dy\,d\vartheta(V)\\
 &\leq C\norm w\int_0^1\int\!\!\int
  \frac{\widetilde b(V,y+s\Gamma_{V^\perp}w)}{\norm{\Gamma_{V^\perp}u-y}}
  \,dy\,d\vartheta(V)\,ds\\
 &=C\norm w\int_0^1\int\!\!\int
  \frac{\widetilde b(V,z)}{\norm{\Gamma_{V^\perp}(u+sw)-z}}
  \,dz\,d\vartheta(V)\,ds,
\end{aligned}
\]
where the last equality follows from the change of variables
$z=y+s\Gamma_{V^\perp}w$.
The double integral over $V$ and $z$ in the last line is
bounded uniformly in $u,w$ and $s\in[0,1]$.
Indeed, on the set $\norm{\Gamma_{V^\perp}(u+sw)-z}\leq1$,
the uniform bound in \eqref{eq:a-tilde} applies, and
$\norm z^{-1}$ is integrable on the unit ball in dimension
$d-\ell+1\geq2$.
On the set $\norm{\Gamma_{V^\perp}(u+sw)-z}>1$,
the integrand is bounded by $\widetilde b(V,z)$,
so the second bound in \eqref{eq:a-tilde} applies.

Consequently, there exists a constant $C>0$ such that, for every
$u,w\in\R^d$,
\[
 \norm{\nabla^2G(u+w)-\nabla^2G(u)}_{\op}\leq C\norm w.
\]
Thus $\nabla^2G$ is Lipschitz, and $G$ is twice
continuously differentiable on $\R^d$.

\smallskip\noindent\emph{The Hessian when $\ell=d$.}\enspace
Let $n:\Gr(d-1,d)\to\R^d$ be a Borel measurable map
assigning to each $V$ a unit normal $n(V)\perp V$.
Applying \eqref{eq:BP} to the gradient formula in
Proposition~\ref{prop:population}, with the change of variables
$y=tn(V)$, $t\in\R$, in the integral over $V^\perp$, we obtain
\[
 \nabla G\left(u\right)=\int_{\Gr\left(d-1,d\right)}n\left(V\right)\int_\R
 \sgn\left(n\left(V\right)^\top u-t\right)b_0\left(V,tn\left(V\right)\right)\,dt\,d\vartheta\left(V\right).
\]
For $u,w\in\R^d$, subtracting the inner integral at $u$
from the one at $u+w$ gives
\[
 2\int_{n\left(V\right)^\top u}^{n\left(V\right)^\top\left(u+w\right)}b_0\left(V,tn\left(V\right)\right)\,dt.
\]
Using the Lipschitz constant $L_{b_0}$, we obtain, for every
$u,w\in\R^d$,
\[
\begin{aligned}
 \norm{\nabla G(u+w)-\nabla G(u)
       -2\int b_0(V,\Gamma_{V^\perp}u)\Gamma_{V^\perp}w\,d\vartheta(V)}
 &\leq2L_{b_0}\int\int_0^{|n(V)^\top w|}s\,ds\,d\vartheta(V)\\
 &\leq L_{b_0}\norm w^2.
\end{aligned}
\]
Thus $\nabla G$ is differentiable, with
\begin{equation}\label{eq:H-oja}
 \nabla^2G(u)=2\int_{\Gr(d-1,d)}
              b_0(V,\Gamma_{V^\perp}u)\Gamma_{V^\perp}\,d\vartheta(V).
\end{equation}
Using \eqref{eq:H-oja}, we also obtain, for every $u,w\in\R^d$,
\[
 \norm{\nabla^2G(u+w)-\nabla^2G(u)}_{\op}
 \leq2\int
 \abs{b_0(V,\Gamma_{V^\perp}(u+w))-b_0(V,\Gamma_{V^\perp}u)}
 \,d\vartheta(V)
 \leq2L_{b_0}\norm w.
\]
Thus $\nabla^2G$ is Lipschitz.

\smallskip\noindent\emph{Positive definiteness.}\enspace
Let $u\in\R^d$ and $v\in\R^d\setminus\{0\}$.
Suppose first that $\ell<d$.
The points $u,u+v,X_1,\ldots,X_\ell$ are affinely independent
almost surely, as established in the proof of
Proposition~\ref{prop:population}.
Thus the vectors
$a\coloneqq\Gamma(\mathbf X)(u-X_\ell)$ and
$\Gamma(\mathbf X)v$ are linearly independent almost surely.
Equality in the Cauchy--Schwarz inequality is therefore
impossible for these two vectors. The formula in
Proposition~\ref{prop:hessian-nonoja} gives
\[
 v^\top\nabla^2G(u)v
 =\frac1\ell\E\left[
 \frac{m_{\ell-1}(\Simpl(\mathbf X))}{\norm a^3}
 \left(\norm{\Gamma(\mathbf X)v}^2\norm a^2
       -\scal{\Gamma(\mathbf X)v}{a}^2\right)
 \right]>0.
\]

Suppose now that $\ell=d$.
Since $f$ is continuous and integrates to one, there exists
an open ball $B(x_0,\varepsilon)$ on which $f>0$.
Let $V\in\Gr(d-1,d)$ and $y\in V^\perp$.
If the hyperplane $y+V$ meets this ball, their intersection
is a nonempty open ball in $y+V$ and contains $d$ affinely
independent points $y+z_1,\ldots,y+z_d$.
Hence $m_{d-1}(\Simpl(\mathbf z))>0$ and $f(y+z_i)>0$
for every $1\leq i\leq d$.
By continuity, these inequalities hold on an open neighborhood
of $\mathbf z$ in $V^d$.
The definition \eqref{eq:a-density} therefore gives
\begin{equation*}
 b_0\left(V,y\right)>0\quad\text{whenever }\norm{y-\Gamma_{V^\perp}x_0}<\varepsilon.
\end{equation*}
Let $\sigma$ be the uniform probability measure on the unit sphere
$\mathbb S^{d-1}\coloneqq\{e\in\R^d:\norm e=1\}$.
The image of $\sigma$ under $e\mapsto\spanop(e)^\perp$
is a rotation-invariant probability measure on $\Gr(d-1,d)$
and therefore equals $\vartheta$.
For $e\in\mathbb S^{d-1}$, the orthogonal projection onto
$\spanop(e)$ is $ee^\top$, so formula~\eqref{eq:H-oja} gives
\[
 v^\top\nabla^2G(u)v
 =2\int_{\mathbb S^{d-1}}
 b_0\bigl(\spanop(e)^\perp,ee^\top u\bigr)
 (e^\top v)^2\,d\sigma(e).
\]
The set
\[
 \mathcal U\coloneqq
 \{e\in\mathbb S^{d-1}:|e^\top(u-x_0)|<\varepsilon\}
\]
is open in $\mathbb S^{d-1}$.
Since $d\geq2$, there exists a unit vector orthogonal to $u-x_0$,
so $\mathcal U$ is nonempty and $\sigma(\mathcal U)>0$.
Since $v\ne0$, we have
$\sigma\bigl(\mathbb S^{d-1}\cap\spanop(v)^\perp\bigr)=0$,
and therefore
\[
 \sigma\bigl(\mathcal U\setminus\spanop(v)^\perp\bigr)
 =\sigma(\mathcal U)>0.
\]
For every $e\in\mathcal U\setminus\spanop(v)^\perp$,
the relations
$\norm{ee^\top(u-x_0)}=|e^\top(u-x_0)|<\varepsilon$
and $e^\top v\ne0$ imply that the integrand is strictly positive.
Consequently, $v^\top\nabla^2G(u)v>0$.

In both cases, $\nabla^2G(u)$ is positive definite for every
$u\in\R^d$.

By Proposition~\ref{prop:population}, $G$ has a unique
minimizer $\mu^*$.
With $H\coloneqq\nabla^2G\left(\mu^*\right)$, we have
\[
 \nabla G\left(\mu^*+u\right)-Hu=\int_0^1\left(\nabla^2G\left(\mu^*+tu\right)-H\right)u\,dt.
\]
By the Lipschitz property of the Hessian, there exists a finite constant $C>0$, independent of $u$, such that
$\norm{\nabla G\left(\mu^*+u\right)-Hu}\leq\left(C/2\right)\norm u^2$. This proves Assumption~\ref{ass:taylor} with $\alpha=1$.
\end{proof}

\begin{remark}
When $\ell=d$, the base--height formula shows that the Hessian
of $g(\mathbf x;\cdot)$ is zero outside $\aff(\mathbf x)$
for every $\mathbf x$.
Under the assumptions of Proposition~\ref{prop:smooth-density},
for every $u\in\R^d$, the Hessian of $g(\mathbf X;\cdot)$ at $u$
is therefore zero almost surely, whereas $\nabla^2G(u)$ in
\eqref{eq:H-oja} is positive definite.
Thus $\nabla^2G(u)$ cannot be obtained by taking the expectation
of the Hessian of $g(\mathbf X;\cdot)$ at $u$.
\end{remark}

\begin{proof}[Verification of the examples following Proposition~\ref{prop:smooth-density}]
For a nondegenerate Gaussian density $f$, there exist constants
$C,\eta>0$ such that, for every $x\in\R^d$,
\[
 f(x)+\norm{\nabla f(x)}\leq Ce^{-\eta\norm{x}^2}.
\]
Hence $f$ satisfies the hypotheses of
Proposition~\ref{prop:smooth-density} for every $\nu>2$.
Finite mixtures of such densities satisfy the same hypotheses.
If $f\in C^1(\R^d)$ has compact support, then $f$ and $\nabla f$
are bounded and vanish outside a compact set, so the same
conclusion holds.

\begin{samepage}
Let $f$ be a Student density of the form
\[
 f\left(x\right)\coloneqq c\left(1+
 \frac{\left(x-a\right)^\top S^{-1}\left(x-a\right)}{\nu}
 \right)^{-\left(d+\nu\right)/2},
\]
where $a\in\R^d$, $\nu>2$, $S$ is a symmetric positive definite
$d\times d$ matrix, and $c>0$ is the normalizing constant.
\end{samepage}
Since $S$ is positive definite, there exists a constant $C>0$
such that, for every $x\in\R^d$,
\[
 f(x)
 \leq C\left(1+\norm{x}^2\right)^{-(d+\nu)/2}
 \leq 2^{(d+\nu)/2}C\left(1+\norm{x}\right)^{-d-\nu}.
\]
Differentiation and this bound give a constant $C'>0$ such that,
for every $x\in\R^d$,
\[
 \norm{\nabla f(x)}
 =\frac{(d+\nu)\norm{S^{-1}(x-a)}}
 {\nu+(x-a)^\top S^{-1}(x-a)}f(x)
 \leq
 \frac{(d+\nu)\norm{S^{-1}}_{\op}}
 {2\sqrt{\nu\lambda_{\min}(S^{-1})}}f(x)
 \leq C'\left(1+\norm{x}\right)^{-d-\nu}.
\]
This proves the hypotheses of
Proposition~\ref{prop:smooth-density} with the same $\nu$.
\end{proof}

\section{Computation of the gradients}\label{sec:gradient-computation}

We describe the computation of the gradients in the
sliding-window recursion and establish the operation counts
stated in Section~\ref{sec:computation}.

We first consider the case $\ell\geq2$, under
Assumption~\ref{ass:density} or, when $\ell<d$,
under Assumption~\ref{ass:moment-alternative}.
Under either assumption, $P$ assigns probability zero to every
affine subspace of dimension at most $\ell-1$.
Successive conditioning on the independent observations in each
window then shows that, with probability one, for every
$n\geq\ell$, the observations $X_{n-\ell+1},\ldots,X_n$
are linearly independent, and hence affinely independent.

\smallskip\noindent\textbf{Evaluation from the Gram matrix.}\enspace
Let $\mathbf x\coloneqq(x_1,\ldots,x_\ell)$ be a tuple of
affinely independent points in $\R^d$.
We recall the subgradient selection in~\eqref{eq:little_g_grad}:
for every $\mu\in\R^d$,
\[
 \nabla_\mu g(\mathbf x;\mu)=
 \begin{cases}
 \displaystyle\frac{m_{\ell-1}(\Simpl(\mathbf x))}{\ell}
 \frac{\Gamma(\mathbf x)(\mu-x_\ell)}
 {\norm{\Gamma(\mathbf x)(\mu-x_\ell)}},
 &\text{if }\Gamma(\mathbf x)(\mu-x_\ell)\ne0,\\[0.8em]
 0,&\text{if }\Gamma(\mathbf x)(\mu-x_\ell)=0.
 \end{cases}
\]

We define the $d\times(\ell-1)$ matrix $B$ of edge vectors
from $x_\ell$ and its Gram matrix $S$ by
\[
 B\coloneqq\left(x_1-x_\ell,\ldots,x_{\ell-1}-x_\ell\right),
 \quad\text{and}\quad S\coloneqq B^\top B.
\]
The affine independence of $x_1,\ldots,x_\ell$ implies that
the columns of $B$ are linearly independent, so $S$ is positive
definite.
For every $\mu\in\R^d$, we have
\[
 m_{\ell-1}\left(\Simpl\left(\mathbf x\right)\right)
 =\frac{\sqrt{\det S}}{\left(\ell-1\right)!},
 \quad\text{and}\quad
 \Gamma\left(\mathbf x\right)\left(\mu-x_\ell\right)
 =\mu-x_\ell
 -B\left(S^{-1}\left(B^\top\left(\mu-x_\ell\right)\right)\right).
\]
Once $B$, $S^{-1}$ and $\det S$ are available, computing
$\Gamma(\mathbf x)(\mu-x_\ell)$ by successive matrix--vector
products and then evaluating the selected subgradient requires
$O(d\ell+\ell^2)=O(d\ell)$ operations.

\smallskip\noindent\textbf{Updating the window.}\enspace
After $n\geq\ell$ observations, only the current window
$X_{n-\ell+1},\ldots,X_n$ is stored.
Since permuting these observations leaves the simplex volume
and the selected subgradient unchanged, we replace
$X_{n-\ell+1}$ by $X_{n+1}$ without moving the other
observations.

We now derive formulas for updating $B$, $S^{-1}$ and $\det S$
when one observation is replaced.
We denote by $j\in\{1,\ldots,\ell\}$ the index of the observation
being replaced and by $x_j'\in\R^d$ the new observation.
We set $u\coloneqq x_j'-x_j$.
We assume that the observations remain affinely independent
after replacement, as they do almost surely for the successive
windows.
The updated matrix $B'$ is defined by the same formula as $B$,
with $x_j$ replaced by $x_j'$.
Thus
\[
 B'=B+uv^\top,
 \qquad
 v\coloneqq
 \begin{cases}
 e_j,&\text{if }j<\ell,\\
 -\mathbf1_{\ell-1},&\text{if }j=\ell,
 \end{cases}
\]
where $e_1,\ldots,e_{\ell-1}$ are the coordinate vectors of
$\R^{\ell-1}$ and
$\mathbf1_{\ell-1}\coloneqq(1,\ldots,1)^\top\in\R^{\ell-1}$.
With $z\coloneqq B^\top u+(\norm u^2/2)v$, the updated
Gram matrix satisfies
\begin{align*}
 S'\coloneqq(B')^\top B'
 &=S+B^\top uv^\top+vu^\top B+\norm u^2vv^\top\\
 &=S+vz^\top+zv^\top=S+(v,z)(z,v)^\top.
\end{align*}
The change from $S$ to $S'$ therefore has rank at most two.
The affine independence of the observations after replacement
implies that $S'$ is positive definite. The matrix determinant lemma gives
\[
 \det S'=\det S\,\det\left(I_2+(z,v)^\top S^{-1}(v,z)\right).
\]
Hence $I_2+(z,v)^\top S^{-1}(v,z)$ is nonsingular, and the
Woodbury formula gives
\[
 (S')^{-1}
 =S^{-1}-\bigl(S^{-1}(v,z)\bigr)
 \left(I_2+(z,v)^\top S^{-1}(v,z)\right)^{-1}
 \bigl((z,v)^\top S^{-1}\bigr).
\]

\smallskip\noindent\textbf{Operation counts.}\enspace
For the first window $X_1,\ldots,X_\ell$, computing $B$ and
$S=B^\top B$ costs $O(d\ell^2)$ operations.
Computing $S^{-1}$ and $\det S$ by factorization costs
$O(\ell^3)$ operations, so initialization costs
$O(d\ell^2)$ since $\ell\leq d$.

At each replacement, computing $u$, $z$ and $B'$ by the
formulas above costs $O(d\ell)$ operations.
Computing $S^{-1}(v,z)$ and $(z,v)^\top S^{-1}$, and then
evaluating the formula for $(S')^{-1}$, costs $O(\ell^2)$
operations.
The $2\times2$ inverse and the determinant update cost $O(1)$.
Including the gradient evaluation, the cost per replacement
is $O(d\ell)$ since $\ell\leq d$.
For $n\geq\ell$, the total cost of the sliding-window
recursion for one pass over $X_1,\ldots,X_n$, including
initialization and the $n-\ell+1$ gradient evaluations, is
therefore
\begin{equation*}
 O\left(d\ell^2+\left(n-\ell+1\right)d\ell\right)
 =O\left(nd\ell\right).
\end{equation*}
Storing the observations, $B$ and $S^{-1}$ requires
$O(d\ell)$ memory.

For the disjoint-block recursion, computing a Gram matrix
and a gradient separately for each of the $\lfloor n/\ell\rfloor$
complete blocks costs $O(d\ell^2)$ operations per block.
The total cost of the disjoint-block recursion for one pass
over $X_1,\ldots,X_n$ is therefore
$O(\lfloor n/\ell\rfloor d\ell^2)=O(nd\ell)$ operations.

\smallskip\noindent\textbf{The case $\ell=1$.}\enspace
For $\ell=1$, the conventions $\Gamma(x_1)=I_d$ and
$m_0(\Simpl(x_1))=1$ give, by~\eqref{eq:little_g_grad},
\[
 \nabla_\mu g(x_1;\mu)=
 \begin{cases}
 \displaystyle\frac{\mu-x_1}{\norm{\mu-x_1}},
 &\text{if }\mu\ne x_1,\\[0.8em]
 0,&\text{if }\mu=x_1.
 \end{cases}
\]
Each gradient evaluation costs $O(d)$ operations.
The total cost of the sliding-window recursion for one pass
over $X_1,\ldots,X_n$ is therefore $O(nd)$ operations,
with $O(d)$ memory.

\smallskip\phantomsection\label{sec:oja-implementation}
\noindent\textbf{An alternative calculation for the Oja median.}\enspace
For $\ell=d$, the preceding Gram-matrix calculation therefore
gives a total cost of $O(nd^2)$ operations for the
sliding-window recursion over $X_1,\ldots,X_n$, with
$O(d^2)$ memory.
The implementation used in Section~\ref{sec::simus} is based
on the square matrix of observations.
Replacing one observation changes this matrix by a matrix
of rank at most one, and the determinant formula below gives
the selected subgradient without computing an orthogonal
projection.

Let $\mathbf x\coloneqq(x_1,\ldots,x_d)$ be a tuple of
linearly independent points in $\R^d$, and let $A$ be the
matrix with columns $x_1,\ldots,x_d$.
We set $\mathbf1_d\coloneqq(1,\ldots,1)^\top\in\R^d$.
The simplex volume formula~\eqref{eq:simplex-gram-volume}
and the matrix determinant lemma give, for every $\mu\in\R^d$,
\[
 g(\mathbf x;\mu)
 =\frac{\abs{\det(A-\mu\mathbf1_d^\top)}}{d!}
 =\frac{\abs{\det A}}{d!}
  \abs{1-\mathbf1_d^\top A^{-1}\mu}.
\]
For $\mu\in\R^d$ such that
$\mathbf1_d^\top A^{-1}\mu\ne1$, differentiation gives
\[
 \nabla_\mu g(\mathbf x;\mu)
 =-\frac{\abs{\det A}}{d!}
 \sgn\left(1-\mathbf1_d^\top A^{-1}\mu\right)
 (A^{-1})^\top\mathbf1_d.
\]
Since $\mathbf1_d^\top A^{-1}\mu=1$ is equivalent to
$\mu\in\aff(\mathbf x)$, the same formula, with
$\sgn(0)\coloneqq0$, gives the selected subgradient
in~\eqref{eq:little_g_grad} for every $\mu\in\R^d$.

For the replacement of $x_j$ by $x_j'$ considered above,
with $u=x_j'-x_j$, we denote by $A'$ the updated matrix
and assume that $A'$ is nonsingular.
Let $e_j$ be the $j$th coordinate vector of $\R^d$.
Then $A'=A+ue_j^\top$, and the matrix determinant lemma gives
\[
 \det A'=\det A\left(1+e_j^\top A^{-1}u\right).
\]
Since $A$ and $A'$ are nonsingular,
$1+e_j^\top A^{-1}u\ne0$, and the Sherman--Morrison formula
for rank-one updates \citep{sherman1950adjustment} gives
\[
 (A')^{-1}
 =A^{-1}-\frac{(A^{-1}u)(e_j^\top A^{-1})}
 {1+e_j^\top A^{-1}u}.
\]
Computing $A^{-1}u$, updating the inverse and determinant,
and evaluating the selected subgradient costs $O(d^2)$
operations per replacement, after an initial $O(d^3)$
factorization.
For $n\geq d$, the total cost of the sliding-window recursion
for one pass over $X_1,\ldots,X_n$, including initialization,
is therefore $O(nd^2)$ operations, with $O(d^2)$ memory.
This gives a simpler implementation than the Gram-matrix
calculation, with a rank-one inverse update and fewer
matrix--vector products for the gradient, while retaining
the same order of complexity.

\section{Gaussian calibration of the step sizes}\label{sec:gaussian-calibration}

We justify the choice $c_\gamma=\sqrt{d!}$ used for the Gaussian
simulations in Section~\ref{sec::simus}.

Suppose that $\ell=d$ and that $X_1,\ldots,X_d$ are independent
with distribution $\mathcal N_d(0,I_d)$. Then $\mu^*=0$.
To compute $\Gamma_0$, for $u\in\R^d$ we define
$D(u)\coloneqq\det\left(X_1-u,\ldots,X_d-u\right)$.
In the multilinear expansion of $D(u)$, all terms containing $u$
in at least two columns vanish. Hence there exists a random vector
$a$ such that $D(u)=D(0)-a^\top u$ for every $u\in\R^d$.
Proposition~2.1 of~\citet{durre_affine-equivariant_2022},
applied with both functions equal to $x\mapsto x-u$, gives
\begin{equation*}
 \E\left[D(u)^2\right]
 =d!\det\left(\E\left[\left(X_1-u\right)\left(X_1-u\right)^\top\right]\right)
 =d!\det\left(I_d+uu^\top\right)
 =d!\left(1+\norm{u}^2\right).
\end{equation*}
Identifying the quadratic terms in $u$ gives, for every
$u\in\R^d$, $u^\top\E\left[aa^\top\right]u=d!\norm{u}^2$,
and hence $\E\left[aa^\top\right]=d!I_d$.
Since $D(0)\ne0$ almost surely, differentiating
$g\left(\mathbf X;u\right)=|D(u)|/d!$ at $u=0$ gives
\begin{equation*}
 \nabla_\mu g\left(\mathbf X;0\right)
 =-\frac{\sgn\left(D(0)\right)}{d!}a.
\end{equation*}
The gradient is centred by symmetry, so
\begin{equation*}
 \Gamma_0=\frac{\E\left[aa^\top\right]}{(d!)^2}
 =\frac1{d!}I_d.
\end{equation*}
Thus the choice $c_\gamma\coloneqq\sqrt{d!}$ used in
Section~\ref{sec::simus} gives $c_\gamma^2\Gamma_0=I_d$, so
each coordinate of $c_\gamma\nabla_\mu g\left(\mathbf X;0\right)$
has variance one.


\begin{thebibliography}{99}
\small
\emergencystretch=3em
\setlength{\itemsep}{0.3em}

\bibitem[Agarwal et~al.(2020)Agarwal, Anil, Koren, Talwar, and Zhang]{agarwal2020stochastic}
Agarwal, N., Anil, R., Koren, T., Talwar, K. and Zhang, C. (2020).
\href{https://proceedings.neurips.cc/paper/2020/hash/74dbd1111727a31a2b825d615d80b2e7-Abstract.html}{Stochastic optimization with laggard data pipelines.}
In \emph{Advances in Neural Information Processing Systems} \textbf{33}, 10282--10293.

\bibitem[Andersen and Gill(1982)]{AndersenGill1982}
Andersen, P. K. and Gill, R. D. (1982).
Cox's regression model for counting processes: A large sample study.
\emph{The Annals of Statistics} \textbf{10}, 1100--1120.
\doi{10.1214/aos/1176345976}.

\bibitem[Benveniste et~al.(1990)Benveniste, M\'etivier, and Priouret]{benveniste-book90}
Benveniste, A., M\'etivier, M. and Priouret, P. (1990).
\emph{Adaptive Algorithms and Stochastic Approximations}.
Applications of Mathematics \textbf{22}. Springer-Verlag, Berlin.
\doi{10.1007/978-3-642-75894-2}.

\bibitem[Brown(1983)]{brown1983statistical}
Brown, B. M. (1983).
Statistical uses of the spatial median.
\emph{Journal of the Royal Statistical Society, Series B} \textbf{45}, 25--30.
\doi{10.1111/j.2517-6161.1983.tb01226.x}.

\bibitem[Cardot et~al.(2017)Cardot, C\'enac, and Godichon-Baggioni]{CCG2017}
Cardot, H., C\'enac, P. and Godichon-Baggioni, A. (2017).
Online estimation of the geometric median in Hilbert spaces: Nonasymptotic confidence balls.
\emph{The Annals of Statistics} \textbf{45}, 591--614.
\doi{10.1214/16-AOS1460}.

\bibitem[Cardot et~al.(2013)Cardot, C\'enac, and Zitt]{CCZ2013}
Cardot, H., C\'enac, P. and Zitt, P.-A. (2013).
Efficient and fast estimation of the geometric median in Hilbert spaces with an averaged stochastic gradient algorithm.
\emph{Bernoulli} \textbf{19}, 18--43.
\doi{10.3150/11-BEJ390}.

\bibitem[Duflo(1997)]{Duflo}
Duflo, M. (1997).
\emph{Random Iterative Models}.
Applications of Mathematics \textbf{34}. Springer-Verlag, Berlin.
Translated by S. S. Wilson.
\doi{10.1007/978-3-662-12880-0}.

\bibitem[D\"urre and Paindaveine(2022a)]{durre_affine-equivariant_2022}
D\"urre, A. and Paindaveine, D. (2022a).
Affine-equivariant inference for multivariate location under $L_p$ loss functions.
\emph{The Annals of Statistics} \textbf{50}, 2616--2640.
\doi{10.1214/22-AOS2199}.

\bibitem[D\"urre and Paindaveine(2022b)]{DPsupp}
D\"urre, A. and Paindaveine, D. (2022b).
Supplement to ``Affine-equivariant inference for multivariate location under $L_p$ loss functions.''
\doi{10.1214/22-AOS2199SUPP}.



\bibitem[Fischer et~al.(2020)Fischer, Mosler, M\"ott\"onen, Nordhausen, Pokotylo, and Vogel]{Fischer2020}
Fischer, D., Mosler, K., M\"ott\"onen, J., Nordhausen, K., Pokotylo, O. and Vogel, D. (2020).
Computing the Oja median in R: The package OjaNP.
\emph{Journal of Statistical Software} \textbf{92}(8), 1--36.
\doi{10.18637/jss.v092.i08}.

\bibitem[Godichon-Baggioni(2019)]{Godichon2019variance}
Godichon-Baggioni, A. (2019).
Online estimation of the asymptotic variance for averaged stochastic gradient algorithms.
\emph{Journal of Statistical Planning and Inference} \textbf{203}, 1--19.
\doi{10.1016/j.jspi.2019.01.001}.

\bibitem[Godichon-Baggioni and Lu(2024)]{GodichonLu2024}
Godichon-Baggioni, A. and Lu, W. (2024).
Online stochastic Newton methods for estimating the geometric median and applications.
\emph{Journal of Multivariate Analysis} \textbf{202}, 105313.
\doi{10.1016/j.jmva.2024.105313}.

\bibitem[Godichon-Baggioni and Werge(2025)]{GodichonWerge2025}
Godichon-Baggioni, A. and Werge, N. (2025).
On adaptive stochastic optimization for streaming data: A Newton's method with $O\left(dN\right)$ operations.
\emph{Journal of Machine Learning Research} \textbf{26}(59), 1--49.
\url{https://www.jmlr.org/papers/v26/23-1565.html}.

\bibitem[Haldane(1948)]{Hal48}
Haldane, J. B. S. (1948).
Note on the median of a multivariate distribution.
\emph{Biometrika} \textbf{35}, 414--417.
\doi{10.1093/biomet/35.3-4.414}.

\bibitem[Hoeffding(1948)]{Hoeffding1948}
Hoeffding, W. (1948).
A class of statistics with asymptotically normal distribution.
\emph{The Annals of Mathematical Statistics} \textbf{19}, 293--325.
\doi{10.1214/aoms/1177730196}.

\bibitem[Hoeffding(1961)]{Hoeffding1961}
Hoeffding, W. (1961).
The strong law of large numbers for $U$-statistics.
Institute of Statistics Mimeo Series No.~302.
University of North Carolina, Chapel Hill.
\url{https://repository.lib.ncsu.edu/handle/1840.4/2128}.

\bibitem[Kemperman(1987)]{Kem87}
Kemperman, J. H. B. (1987).
The median of a finite measure on a Banach space.
In \emph{Statistical Data Analysis Based on the $L_1$-Norm and Related Methods} (Y. Dodge, ed.), 217--230.
North-Holland, Amsterdam.

\bibitem[Kushner and Yang(1993)]{KushnerYang1993}
Kushner, H. J. and Yang, J. (1993).
Stochastic approximation with averaging of the iterates: Optimal asymptotic rate of convergence for general processes.
\emph{SIAM Journal on Control and Optimization} \textbf{31}, 1045--1062.
\doi{10.1137/0331047}.

\bibitem[Kushner and Yin(2003)]{Kushner2003stochastic}
Kushner, H. J. and Yin, G. G. (2003).
\emph{Stochastic Approximation and Recursive Algorithms and Applications}, second edition.
Applications of Mathematics \textbf{35}. Springer-Verlag, New York.
\doi{10.1007/b97441}.



\bibitem[Miles(1971)]{Miles1971}
Miles, R. E. (1971).
Isotropic random simplices.
\emph{Advances in Applied Probability} \textbf{3}, 353--382.
\doi{10.2307/1426176}.

\bibitem[Oja(1983)]{oja1983descriptive}
Oja, H. (1983).
Descriptive statistics for multivariate distributions.
\emph{Statistics \& Probability Letters} \textbf{1}, 327--332.
\doi{10.1016/0167-7152(83)90054-8}.

\bibitem[Oja and Niinimaa(1985)]{oja1985asymptotic}
Oja, H. and Niinimaa, A. (1985).
Asymptotic properties of the generalized median in the case of multivariate normality.
\emph{Journal of the Royal Statistical Society, Series B} \textbf{47}, 372--377.
\doi{10.1111/j.2517-6161.1985.tb01366.x}.

\bibitem[Pelletier(2000)]{pelletier2000asymptotic}
Pelletier, M. (2000).
Asymptotic almost sure efficiency of averaged stochastic algorithms.
\emph{SIAM Journal on Control and Optimization} \textbf{39}, 49--72.
\doi{10.1137/S0363012998308169}.

\bibitem[Polyak and Juditsky(1992)]{PolyakJud92}
Polyak, B. T. and Juditsky, A. B. (1992).
Acceleration of stochastic approximation by averaging.
\emph{SIAM Journal on Control and Optimization} \textbf{30}, 838--855.
\doi{10.1137/0330046}.

\bibitem[Robbins and Monro(1951)]{robbins1951}
Robbins, H. and Monro, S. (1951).
A stochastic approximation method.
\emph{The Annals of Mathematical Statistics} \textbf{22}, 400--407.
\doi{10.1214/aoms/1177729586}.

\bibitem[Robbins and Siegmund(1971)]{RobbinsSiegmund1971}
Robbins, H. and Siegmund, D. (1971).
A convergence theorem for non negative almost supermartingales and some applications.
In \emph{Optimizing Methods in Statistics} (J. S. Rustagi, ed.), 233--257.
Academic Press, New York.
\doi{10.1016/B978-0-12-604550-5.50015-8}.

\bibitem[Rockafellar and Wets(1998)]{rockafellar_wets_1998}
Rockafellar, R. T. and Wets, R. J.-B. (1998).
\emph{Variational Analysis}.
Grundlehren der mathematischen Wissenschaften \textbf{317}.
Springer-Verlag, Berlin.
\doi{10.1007/978-3-642-02431-3}.

\bibitem[Ros\'en(1967)]{Rosen1967}
Ros\'en, B. (1967).
On the central limit theorem for sums of dependent random variables.
\emph{Zeitschrift f\"ur Wahrscheinlichkeitstheorie und Verwandte Gebiete} \textbf{7}, 48--82.
\doi{10.1007/BF00532097}.

\bibitem[Ruppert(1988)]{ruppert1988efficient}
Ruppert, D. (1988).
Efficient estimators from a slowly convergent Robbins--Monro process.
Technical Report 781, School of Operations Research and Industrial Engineering, Cornell University, Ithaca, NY.
\url{https://hdl.handle.net/1813/8664}.

\bibitem[Schneider and Weil(2008)]{SchneiderWeil2008}
Schneider, R. and Weil, W. (2008).
\emph{Stochastic and Integral Geometry}.
Probability and Its Applications. Springer-Verlag, Berlin.
\doi{10.1007/978-3-540-78859-1}.

\bibitem[Sherman and Morrison(1950)]{sherman1950adjustment}
Sherman, J. and Morrison, W. J. (1950).
Adjustment of an inverse matrix corresponding to a change in one element of a given matrix.
\emph{The Annals of Mathematical Statistics} \textbf{21}, 124--127.
\doi{10.1214/aoms/1177729893}.

\end{thebibliography}
\end{document}